\documentclass[preprint,11pt, 3p]{elsarticle}
\biboptions{sort&compress}

\usepackage{amssymb}
\usepackage{amsmath}
\usepackage{amsthm}

\newdefinition{rmk}{Remark}
\usepackage{stmaryrd}
\usepackage{mathsemantics}
\usepackage{subcaption}
\usepackage{placeins}
\usepackage{enumitem}
\usepackage{comment}
\usepackage{tikz}
\usetikzlibrary{hobby}

\newcommand{\tang}[1]{T_{\bp}(#1)}
\newcommand{\normal}[1]{\bn_{#1}(\bp)}

\newcommand{\matder}{\partial^\bullet}
\newcommand{\aleder}{\partial^\cA}

\newcommand{\normder}{\partial^\circ}
\newcommand{\Id}{\mathbf{Id}}

\newcommand{\autofigref}[1]{Figure \ref{#1}}

\usepackage{lineno}

\journal{Computer Methods in Applied Mechanics and Engineering}

\begin{document}


\begin{frontmatter}



\title{A finite element framework for solving coupled multiphysics problem with moving boundaries in cell biophysics\tnoteref{label1}}




%

\tnotetext[label1]{Submitted to the editors.}


\author[NTNU]{Alessandro Contri\corref{cor1}}
\cortext[cor1]{Corresponding author.}
\ead{alessandro.contri@ntnu.no}
\ead[url]{https://www.ntnu.no/ansatte/alessandro.contri}
\author[NTNU]{Andr\'e Massing}
\ead{andre.massing@ntnu.no}
\ead[url]{https://www.ntnu.no/ansatte/andre.massing}
\author[UCSD]{Padmini Rangamani}
\ead{prangamani@health.ucsd.edu}
\ead[url]{https://rangamani.ucsd.edu}

\affiliation[NTNU]{organization={Department of Mathematics, NTNU},
            city={Trondheim},
            country={Norway}}
\affiliation[UCSD]{organization={Department of Pharmacology, UCSD School of Medicine},
            city={San Diego},
            state={CA},
            country={USA}}


\begin{abstract}
Cellular morphodynamics requires solving systems of coupled partial
differential equations on moving bulk and surface domains, where
advection-dominant transport, structure preservation, and severe mesh
distortions make robust simulation difficult.  
We present a holistic finite element framework that jointly addresses these obstacles for
biophysical applications by combining model-agnostic
structure-preserving postprocessing, ALE-based mesh redistribution
strategies driven by surface-tangential velocities, and stabilized
discretization for advection-diffusion-reaction problems tailored to
evolving domains.  
The methodology is modular and applies to
advection-diffusion-reaction systems, Cahn--Hilliard phase separation,
Helfrich-type geometric flows, as well as their staggered and potentially mixed-dimensional couplings.
We provide a concise notation for evolving bulk and surface geometries, extend positivity-, bound-, and mass-preserving projections to moving meshes, and develop a two-step redistribution procedure that maintains element quality without remeshing.
Convergence studies, manufactured solutions, and biologically motivated test cases ---including tumor-growth surrogates and phase segregation on deformable membranes--- demonstrate accuracy, stability, and versatility across the problem classes considered. 
\end{abstract}



\begin{highlights}
\item Fast and accurate structure preserving scheme for moving bulk and surface PDEs
\item Flexible  mesh-redistribution algorithm for biological moving boundary problems
\item Stabilized advection-dominant transport equations on moving surfaces
\item Spacetime convergence for structure preserving PDEs on moving surfaces
\item Application to state-of-the-art coupled problems involving deformable membranes
\end{highlights}


\end{frontmatter}



\section{Introduction}
\label{sec:intro}
\subsection{Background and state of the art}
\label{subsec:background_and_soa}

Many biophysical phenomena in cell biology can be characterized by the feedback between form and function.
Cells routinely change their shape during cell motility, wound healing, and development.
Such shape changes are a result of coordination of multiple biochemical and mechanical processes within the cell.
Not surprisingly, these phenomena have been studied extensively using models based on partial differential equations (PDEs) with different levels of biophysical complexity \cite{LomakinLeeHanEtAl2015, AuddyaZhangGulatiEtAl2021}.
In concert with experiments, such models have provided insight into the mechanisms underlying the change in cell shape \cite{DayelAkinLanderyouEtAl2009, Sens2020, KerenPincusAllenEtAl2008}.
However, an open challenge in the field of computational cell
biophysics is the development of robust computational tools for
numerical solutions of the governing PDEs. 
For example, cellular membranes are mechanical boundaries that are characterized by lateral dimensions that are large compared to their thickness, appearing in a variety of structures in cells \cite{Steigmann_book}.
It has been proven that the properties of these surfaces are crucial for the genesis and maintenance of cellular membrane systems~\cite{BaumgartHessWebb2003,Mullins_book}. 
These restructuring phenomena strongly depend on the components of the membranes themselves and their surroundings \cite{BrangwynneEckmannCoursonEtAl2009, BananiLeeHymanEtAl2017,ShinBrangwynne2017, BeutelMaraspiniPombo-GarciaEtAl2019, LeeCatheyWuEtAl2020, ZhaoZhang2020, AlbertiHyman2021, AmbroggioCostaNavarroPerezSocasEtAl2021, SneadJalihalGerbichEtAl2022}, leading to PDE systems coupling complex multiphysics to free surfaces problems. The resulting systems typically exhibit highly non-linear dynamics and a general numerical framework to simulate these rich dynamics in a unified and stable setting is still missing.
The goal of this article is to present a finite element based computational framework which addresses some of the key challenges arising in continuum-mechanics based computational models of cells. 

Numerous numerical methods have been proposed to deal with the challenges posed by a moving deformable membrane \cite{IyerGompperFedosov2023, OngLai2020, AlandEgererLowengrubEtAl2014, LaadhariSaramitoMisbah2014, RosolenPecoArroyo2013, KrugerVarnikRaabe2011,Noguchi2004,MokbelMokbelLieseEtAl2024,BachiniKrauseNitschkeEtAl2023}, including: 
parametric finite element methods, immersed boundary methods,
level-set method, mesh-free methods, particle methods, and the
phase-field method. We refer to~\cite{AndreaBonito2020}[Chapter 4-6] for a detailed review of some of these methods.
However, even the most advanced schemes often do not simultaneously address the following features
which are frequently crucial in the context of biophysical cell modeling.

First, it is not rare, in the targeted applications, to incur \emph{advection-dominant} equations possibly posed on a \emph{moving surface}. Characteristic examples can be the transport of surfactants in two-phase flow~\cite{GrossReusken2011,GanesanHahnSimonEtAl2017} or the interaction between actin's barbed ends and the cellular membrane~\cite{Bonilla-QuintanaRangamani2024,DoubrovinskiKruse2011}. There exists a large bulk of literature for parabolic advection-diffusion-reaction (ADR) PDEs  on moving surfaces 
\cite{DziukElliott2007, DziukLubichMansour2012, BonitoKyzaNochetto2013,LubichMansourVenkataraman2013, DziukElliott2013, KovacsLiLubichEtAl2017, KovacsPowerGuerra2018, ElliottRanner2021, CaetanoElliottGrasselliEtAl2023}, but primarily with a focus on the diffusion-dominant regime.
For what concerns advection-dominant equations, the bulk setting has been extensively explored and the surface case has seen recent developments ~\cite{OlshanskiiReuskenXu2014,BurmanHansboLarsonEtAl2020, SimonTobiska2019, BurmanHansboLarsonEtAl2018b, DednerMadhavan2015,UlfsbyMassingSticko2023,XiaoZhaoFeng2020}. Surprisingly, most aforementioned approaches only consider the case of
stationary surfaces, except for \cite{HansboLarsonZahedi2015},
where a characteristic cut finite element method on moving surfaces is proposed.

Second, the major part of the existing literature on surface PDEs pays little attention to
structure preservation of physical quantities.
In contrast, computational cell biophysics models often require the preservation of physical constraints ---including mass conservation, density positivity, phase bounds, and energy decay--- to maintain both physical interpretability and numerical stability of the computed solution.
While there is a vast and diversified literature on bulk PDEs
dedicated to preserve positivity, bounds and mass~\cite{LuHuangVanVleck2013,LiYangZhou2020, DuJuLiEtAl2021, LiZhang2020, ZhangShu2011, ZhangShu2010, LiuWangZhou2018, ChenWangWangEtAl2019,HuangShen2021, HuHuang2020} 
this richness of options is unfortunately lacking if we consider PDEs on surfaces, let alone moving surfaces.

Finally, in the context of cell shape dynamics, the domain
evolution depends on a set of coupled PDEs and is not known
\textit{a priori}.  For classical prototypical geometric flows, such as mean curvature flow,
surface diffusion flow, and Helfrich flow, a large number of numerical
methods have been proposed, see, e.g. \cite{Dziuk2008, BarrettGarckeNurnberg2007,
 M.ElliottFritz2017, HuLi2022, DuanLi2024, BaiHuLi2024,
 GarckeNurnbergZhao2025}, assuming mostly idealized geometries.
It is well-known that for most finite-element based algorithms, such geometric flows notoriously drive
surface meshes towards bad conditioning. 
Consequently, it is crucial to design numerical methods that
can maintain good mesh quality during the shape evolution
also in the presence of more realistic geometries
resulting e.g. from imaging techniques
\cite{LeeLaughlinAnglivieldeLaBeaumelleEtAl2020,
YangVenkataramanStylesEtAl2016, GulatiRudraraju2024},
or when coupled to other bulk-surface PDEs~\cite{FrancisLaughlinDokkenEtAl2024}.


\subsection{New contributions and outline}
\label{subsec:new_contributions}
The present work aims to present a unified finite-element based computational framework which focuses 
on the three challenges outlined above: advection-dominant transport
on moving surfaces, structure-preservation, and handling of large mesh
deformations in problems where
multi-physics, potential multi-dimensional problems are nonlinearly coupled to the shape evolution.
It is meant to bridge the complexity gap between numerical accuracy
and biophysical complexity in moving domains scenario.
While there exist a plethora of numerical methods which each target a
specific problem class or particular issue, we here employ and extend
some of the most recent state-of-the-art methods and combine them into
a holistic simulation framework which can handle coupled multiphysics
and shape deformation problems arising in computational cell biology.

We start in Section \ref{sec:model} by providing the notation
necessary for the formulation and discretization of PDEs posed on
evolving domains.  
For purposes of generality, the notation is kept general to cover both evolving surface and
bulk domains, only distinguishing between the two as needed. 
In Section \ref{sec:structure_pres}, we adapt the structure-preserving
scheme in \cite{ChengShen2022} to moving bulk and surface problems as
a general approach to design positivity, bound and mass preserving
schemes.
The model-agnostic, pointwise approach of \cite{ChengShen2022} is fit to the first-order, piecewise linear finite element framework used here. 
The result is a fast, generalizable and flexible approach that can be
readily adapted to multiple PDEs posed on moving bulk and surfaces.
In our context, it is used to bound PDE solution to
physically-relevant values while at the same time preserving the total
mass. 

Then, in Section \ref{sec:mesh_redistribution}, we present a mesh
redistribution algorithm based on the surface tangential motion of
\cite{DuanLi2024} and an Arbitrary Lagrangian Eulerian (ALE) approach.
The two-step scheme in \cite{DuanLi2024} is decoupled from the
gradient flow dynamics and applied to general motions. 
This allows us to greatly alleviate the need of remeshing when simulating domains
characterized by large deformations. At the same time, the scheme does
not interfere with the shape evolution, maintaining accuracy. The main advantage of the
schemes presented in Section \ref{sec:structure_pres} and Section
\ref{sec:mesh_redistribution} is that they can be realized as a post-processing step
and thus represent tunable features of the simulation framework.
In Section \ref{sec:adr} we augment an ADR scheme based on the evolving surface finite element method
with a continuous interior penalty (CIP) stabilization to handle
advection-dominant problems.
We demonstrate that stabilized ADR discretization, potentially combined with a structure preservation postprocessing step, yields optimal convergent schemes for the problem class.
The flexibility of the setting is further
verified with an additional real-case example. 
Furthermore,
we apply structure preservation to a Cahn-Hilliard
type equation on moving surfaces in Section \ref{sec:phase_sep}.
Numerical tests are presented that
highlight how the algorithm in Section \ref{sec:structure_pres} is
readily applicable to moving surfaces and preserves the underlying
physics.  
In Section \ref{sec:geom_flows}, we combine the efficient mesh redistribution scheme of Section \ref{sec:mesh_redistribution} with the parametric method of Barrett, Garke and N\"{u}rnberg \cite{BarrettGarckeNurnberg2017a} for Helfrich flow. 
The postprocessing step is extensively tested with state-of-the-art examples revealing good convergence properties and long-term stability. 
Afterwards, in Section \ref{sec:coupling}, 
we couple the solvers from Section \ref{sec:adr} through \ref{sec:geom_flows} together in a staggered approach. 
Convergence is tested with manufactured solutions~\cite{KovacsLiLubichEtAl2017} where the ADR algorithm and a simplified version of the Helfrich algorithm are coupled. In the same context a PDE system modeling tumor growth is also reproduced. The applicability of our approach to relevant cell biophysics is further explored simulating phase separation on a deformable membrane, where the presented schemes for Cahn-Hilliard and Helfrich are coupled.

\section{Model}
\label{sec:model}
This article explores numerical approaches to biophysical cell phenomena where the dynamics involve moving boundaries, such as cell membranes, that are assumed to have negligible thickness and thus are treated as moving hypersurfaces.
As a result, we begin by outlining the essential notation required for both formulating and discretizing partial differential equations on evolving domains, with specific attention given to surface-bound PDEs \cite{DziukElliott2013, BarrettGarckeNurnberg_bookSection}.

Let $M\subset \bbR^d, \; d=2,3$, be an $m$-dimensional oriented, compact, $C^2$ manifold with boundary $\partial M$. The notation $ \Omega \equiv  M$ when $\dim{(M)}=d$ and $\Gamma \equiv M$ when $\dim{(M)}=d-1$ will also be used. In the case $\dim(M)=d-1$, we denote the tangent space with respect to $\bp\in M$ with $\tang{M}$ and the associated tangent bundle with $TM$. The  \emph{unit normal vector} $\normal{M}$ to $M$ at $\bp$ is defined as the vector $\normal{M} \in \bbR^d$ such that $\normal{M} \perp \tang{M}$, $\norm{\normal{M}} = 1$ and $\normal{M}$ agrees with the orientation given. The norm $\norm{\cdot}$ is the standard Euclidean norm in $\bbR^d$. For $\dim{(M)}=d-1$ we define the \emph{tangential projection} at $\bp$ as 
\begin{equation}
    \bbP_{M}(\bp) = \bbI - \normal{M}\otimes \normal{M},
\end{equation}
where $\bbI$ is the identity matrix. For $\dim{(M)}=d$, we set $\bbP_{M}(\bp) = \bbI$ and $\normal{M}=\mathbf{0}$. This allows us to define the \emph{tangential gradient} of a differentiable function $f:M\to \bbR$ at $\bp$ as
\begin{equation}
    \nabla_{M}f(\bp) = \nabla \bar{f}(\bp) - \bbP_{M}(\bp)\nabla \bar{f}(\bp),
\end{equation}
where $\bar{f}$ is a smooth extension of $f$ to a $d$-dimensional neighborhood of $M$ such that $\restr{\bar{f}}{M} = f$ and $\nabla$ is the Euclidean gradient in $\bbR^d$. Analogously, the \emph{tangential divergence} is defined as $\nabla_{M}\cdot f(\bp) = \mathrm{Tr}(\nabla_{M}f(\bp))$. This leads to the definition of the Laplace-Beltrami operator for a $C^2$-function $f:M\to \bbR$ as 
\begin{equation}
    \Delta_{M}f(\bp) = \nabla_{M}\cdot \nabla_{M}f(\bp).
\end{equation}
For $\bf:M\to \bbR^d$ and $\bbF:M\to \bbR^{d\times d}$ we define $(\nabla_{M}\bf)_{ij} =(\nabla (\bf\cdot\be_i))_j$ and $(\nabla_{M}\cdot \bbF)_i = \nabla_{M}\cdot(\bbF^T\be_i)=\mathrm{Tr}(\nabla_{M}(\bbF^T\be_i))$, where $\be_i, ~i \in \{1,\ldots,d\}$ is the canonical basis in $\bbR^d$. This allows us to define $\Delta_M\bf = \nabla_M\cdot\nabla_M\bf $. 
In order to analyze curvature properties of $M$, we introduce the \emph{extended Weingarten map} $\bbH$ as
\begin{equation}
    \bbH(\bp) = -\nabla_{M}\normal{M}.
\end{equation}
It can be shown that $\bbH$ is symmetric, has an eigenvalue $0$ in the direction of the normal, and restricts to the Weingarten map $\bbW(\bp)$ on the tangent space $T_{\bp}M$. For $\bp\in M$, the \emph{mean curvature} is defined as
\begin{equation}
    \kappa(\bp) = \mathrm{Tr}(\bbH(\bp)).
\end{equation}
We recall that all of the above definitions are independent of the extension $\bar{f}$. A sketch of possible model setups is shown in \autofigref{fig:model}.

\begin{figure}[tbhp!]
  \centering
  \begin{tikzpicture}[scale=2, ,use quick Hobby shortcut]
        \fill[gray!15, thick, draw = black, hobby] plot coordinates
                {(-2, 0)  (-2, 0.1) (-2.1, 1) (-1, 1.7) (0.1, 1.2) (0, 0.1) (0, 0)};
        \node at (-1, 1) {\small $\Gamma$};

        \draw[thick, ->,  blue] (-2,0) -- +(0,-0.4);
        \node[anchor = east] at (-2, -0.3) {\small \textcolor{blue}{$\bn_{\partial\Gamma}$}};

        \draw[thick, ->] (-2.1, 1) -- +(-0.3,0.1);
        \node[anchor = east] at (-2.1, 0.9) {\small $\bn_{\Gamma}$};
        \fill[gray!35, thick, dashed, draw = blue] (0,0) arc[start angle=0,end angle=180,
                                x radius=1,y radius=0.3];

        \fill[gray!35, draw = blue, thick] (0,0) arc[start angle=0,end angle=-180,
                                        x radius=1,y radius=0.3];
        \node at (-1,-0.2) {\small \textcolor{blue}{$\partial\Gamma$}};
        
        \fill[gray!15, thick, draw = black] (0.5, -0.25) rectangle (2.5,1.75);

        \fill[gray!35, thick, draw = blue, closed hobby] plot coordinates
            {(1, 0.75)  (0.9, 1.25)  (1.5, 1.25)
            (2, 1.25) (2, 0.5) (1.5, 0.25) (0.8, 0.25)};

        \node at (1.3,1.1) {\small \textcolor{blue}{$\Gamma$}};
        \draw[thick, ->, blue] (2, 1.25) -- +(0.15,0.3);
        \node at (1.9,1.4) {\small \textcolor{blue}{$\bn_{\Gamma}$}};

        \node at (2,0.0) {\small $\Omega$};
        \node at (2.7,0.0) {\small $\partial\Omega$};
        \draw[thick, ->] (2.5, 0.75) -- +(0.4,0);
        \node at (2.7,0.65) {\small $\bn_{\partial\Omega}$};

    \end{tikzpicture}
  \caption{Sketch of an open and closed geometry for notation purposes. One is a 3D object and the other is a 2D object.}
  \label{fig:model}
\end{figure}

We consider a time interval $J = [0, T], \; T>0$ and a $C^2$-\emph{evolving manifold} $\{M(t)\}_{t\in J}$ in $\bbR^d$. In this article an evolving manifold is modeled by a \emph{reference manifold} $\widehat{M}$ together with a \emph{flow map}
\begin{equation}
    \bPhi: J\times\widehat{M}\to \bbR^d, \quad \bPhi,\bPhi^{-1} \in C^1(J; C^2(\bbR^d; \bbR^d)),
    \label{eq:flow_map}
\end{equation}
such that
\begin{itemize}
    \item denoting $M(t) = \bPhi(t, \widehat{M})$, the map $\bPhi_t:=\bPhi(t,\cdot):\widehat{M} \to M(t) $ is a $C^2$-diffeomorphism with inverse map $ \bPhi_{t}^{-1}: M(t) \to \widehat{M} $,
    \item $\bPhi_0 = \mathbf{Id}_{\widehat{M}}$, where $\mathbf{Id}_{\widehat{M}}$ is the identity map on the reference manifold.
\end{itemize}
The mapping $\bPhi_t$ can be used to define pull-back and push-forward maps of functions \cite{AlphonseElliottStinner2015, AlphonseCaetanoDjurdjevacEtAl2023}
\begin{equation}\begin{aligned}
    \bPhi_{-t}f = f \circ \bPhi_t : \widehat{M} \rightarrow \bbR \text{ for } f : M(t) \rightarrow \bbR ,\\
    \bPhi_t f = f \circ \bPhi_t^{-1} : M(t) \rightarrow \bbR \text{ for } f : \widehat{M} \rightarrow \bbR.
\end{aligned}\end{equation}
Often we will call the space-time set $\cG_T=\cup_{t\in J}(\{t\}\times M(t))$ a $C^2$-\emph{evolving manifold} and identify it with $\{M(t)\}_{t\in J}$.
In addition, 
we assume that there exists a velocity field $\bv:J\times \bbR^d \to \bbR^d$ with $\bv \in C^0(J;C^2(\bbR^d;\bbR^d))$ such that for any $t\in J$ and every $\bp \in \widehat{M}$
\begin{equation}
    \frac{\d}{\d t}\bPhi_t(\bp) = \bv(t, \bPhi_t(\bp)).
\end{equation}
For $\bx\in M(t)$, the velocity $\bv$ can be split into a tangential component $\bv^\top(t, \bx) \in T_{\bx}M(t) $ and a normal one $\bv^\perp(t, \bx) = \bv(t, \bx) - \bv^\top(t, \bx)$. This allows us to define a \emph{normal derivative} of a scalar function $f$ on $M(t)$ as
\begin{equation}
    \normder f(t, \bx) = \partial_t\bar{f}(t, \bx) + \bv^{\perp}(t, \bx)\cdot\nabla\bar{f}(t, \bx),
\end{equation}
where $\bar{f}$ is a smooth extension of $f$ in a space-time neighborhood of $M(t)$. The \emph{material derivative} is defined as 
\begin{equation}
    \matder f(t, \bx) = [\normder f + \bv^{\top}\cdot\nabla_{M(t)}\bar{f}] (t,\bx) = [\partial_t\bar{f} + \bv\cdot\nabla\bar{f}](t,\bx).
\end{equation}
The definition of a flow map is not unique and the same evolving manifold can be described by different maps. The computational domain's flow map can be chosen so to maintain good mesh properties while keeping the original map $\bPhi$ for the PDE evolution. This is achieved by introducing a second flow map $\bPhi^\cA$ called \emph{Arbitrary Lagrangian Eulerian (ALE) map} with corresponding velocity $\bv^{\cA}$ and ALE material derivative $\aleder$. It can be shown that $\bv^\cA(t,\bx) - \bv(t,\bx) \in T_{\bx}M(t)$ and that
\begin{equation}
    \aleder f (t,\bx)= [\matder f + (\bv^\cA - \bv)\nabla_{M(t)}f](t,\bx).
    \label{eq:ale_map}
\end{equation}
It is worth noting that the normal velocity $\bv^{\perp}$ is in any case uniquely determined, i.e. $\bv^{\perp}(t,\cdot) = \bv^{\cA,\perp}(t,\cdot), \; \forall t\in J$, including the normal velocity of the boundary $\partial M(t)$.

\subsection{Space discretization}
The geometries are discretized using piecewise linear elements, where for surfaces we follow the setup of \cite{DziukElliott2013}. An $m$-dimensional manifold $M$ is approximated by a triangulated $m$-dimensional domain denoted by $M_{h}$. The elements of the discretization will be distinguished based on their codimension. 
We suppose $M_h$ is composed by a collection $ \cT_h$ of $m$-simplices which vertices $ \by_i, \; i = 1, \ldots, N$ lie on $M$ and such that
\begin{equation}
     \cT_h=\{ T\}, \quad M_{h} = \bigcup T.
\end{equation}
For each $ T\in  \cT_h$ we denote by $h( T)$ its diameter and define the \emph{mesh-size} as
\begin{equation}
    h = \max_{ T\in  \cT_h}h( T).
\end{equation}
The boundary of every element $ T$ is composed by $m+1$ facets of dimension $m-1$, forming a collection $ \cF=\{ F\}$. The discrete boundary of $M_h$ is defined as the union of those facets which defining points lie on $\partial M$ and is denoted $\partial M_h$. It is useful to distinguish between \emph{internal} and \emph{boundary} facets, defined as
\begin{align*}
     \cF_h^i &= \{ F\in  \cF: \;  F\cap \partial M_h = \emptyset\}, \\
     \cF_h^b &= \{ F\in  \cF: \;  F\cap \partial M_h =  F\}.
\end{align*}
Every element $ T$ has a (constant) normal which we will denote by $\bn_{ T}$. In the case $\mathrm{dim}(M_h)\equiv d$ we define $\bn_{ T} \equiv 0$. In the case $\mathrm{dim}(M_h)\equiv d-1$ we define a discrete projection $\bbP_{ T} = \bbI - \bn_{ T}\otimes \bn_{ T}$ and relative tangential gradient, understood in an element-wise sense.
The facet normal, which we denote $\bn_{ F}$, is well-defined once we consider it on the boundary of an element $ T$. For this reason for each $ F\in \partial T$, $\bn_{ F}$ is  the unit co-normal vector perpendicular to both $\bn_{ T}$ and $ F$ and directed outward with respect to the element $ T$. We will also use the notation $\bn_F^{\pm}$, where given as $T^+$ and $T^-$ the neighboring elements to $F$, then $ \bn_F^+$ is the normal as seen from $T^+$ and $ \bn_F^-$ is the normal as seen from $T^-$. In the bulk case we simply have that $\bn_F^- = - \bn_F^+$, instead for surface cases the two normals are not co-planar. An illustration of these objects is given in \autofigref{fig:discrete_normals}.
\begin{figure}[tbhp!]
  \centering
  \begin{tikzpicture}[scale=2, ,use quick Hobby shortcut]

        \fill[gray!35, thick, draw=black] (0,0) -- (2,1) -- (2.2,0.3) -- cycle;
        \draw[thick, ->] (1.2, 0.4) -- (0.92,1);
        \node at (0.87,0.65) {\small $\bn_{T^+}$};
        \node at (1.5,0.35) {\small $T^{+}$};
        
        \fill[gray!15, thick, draw=blue] (2,1) -- (2.2,0.3) -- (4,0) -- cycle;
        \draw[thick, ->, draw = blue] (3, 0.4) -- (3.3,1);
        \node at (3.4,0.65) {\small \textcolor{blue}{$\bn_{T^-}$}};
        \draw[thick, ->, draw = blue] (3, 0.4) -- (3.3,1);
        \node at (2.8,0.35) {\small \textcolor{blue}{$T^{-}$}};

        \draw[thick, ->] (2.1, 0.65) -- (2.8,0.85);
        \node at (2.5,0.6) {\small $\bn_{F^+}$};
        \draw[thick, ->, draw = blue] (2.1, 0.65) -- (1.45,0.9);
        \node at (1.85,0.65) {\small \textcolor{blue}{$\bn_{F^-}$}};

    \end{tikzpicture}
  \caption{Sketch of an discrete discrete mesh entities of a discrete surface $M_h$.}
  \label{fig:discrete_normals}
\end{figure}
We will again use the notation $\Omega_h$ and $\Gamma_h$ when we need to distinguish between $d$-dimensional and $(d-1)$-dimensional discrete manifolds, respectively.

Throughout this work, we use a first-order fitted finite element method to discretize the PDEs in space \cite{ErnGuermond_booka}. In particular, for surfaces, we employ the surface finite element method (SFEM) as presented in \cite{ DziukElliott2013}. The corresponding finite element space is defined as
\begin{equation}
    \label{eq:fem_space}
    V_h(M_h) := \{v_h \in C(M_h) \;|\; \restr{v_h}{ T} \text{ is linear for all }  T \in  \cT_h\},
\end{equation}
which can be described as the span of $N$ piecewise linear continuous basis functions $\phi_i, \; i = 1, \ldots, N$ such that $\phi_i(\by_j) = \delta_{ij}$. Occasionally, we will also need the space
\begin{equation}
    \label{eq:fem_space_0}
    V_{h0}(M_h) := \{v_h \in C(M_h) \;|\; \restr{v_h}{ T} \text{ is linear for all }  T \in  \cT_h \text{ and } v_h|_{\partial M_h}=0\}.
\end{equation}
Integration over mesh entities follows the same notation as for the continuous one, where the summation over every element of the collection is implicit, i.e. we write
\begin{equation}
    \inner{\cdot}{\cdot}_{M_h} = \sum_{T\in M_h}\inner{\cdot}{\cdot}_T, \quad   \inner{\cdot}{\cdot}_{\cF_h^i} = \sum_{F\in \cF_h^i}\inner{\cdot}{\cdot}_F,
\end{equation}
where $\inner{\cdot}{\cdot}$ is a certain inner product. The same element-wise summation is used in case of discrete norms.

\begin{rmk}
    Together with the fitted finite element discretization explained above, an unfitted approach can also be adopted. We refer to \cite{BurmanClausHansboEtAl2015, PalzhanovZhiliakovQuainiEtAl2021, WangPalzhanovQuainiEtAl2022, OlshanskiiPalzhanovQuaini2022, NeivaTurlier2025, GargOlshanskii2025} for advances in this direction.
\end{rmk}

\subsection{Time discretization}
\label{subsec:time_discretization}

The evolving manifold $M(t)$ is discretized using the reference discretization $\widehat{M}_h$ and transport it using a suitable ALE map $\bPhi^\cA$ \cite{HirtAmsdenCook1974, HughesLiuZimmermann1981, Nobile2001, Gastaldi2001, FormaggiaNobile2004, BadiaCodina2006, BonitoKyzaNochetto2013, KovacsPowerGuerra2018, AlphonseCaetanoDjurdjevacEtAl2023}. The reference points $ \widehat{\by}_i\in \widehat{M}_h$ defining the simplices of the reference triangulation $ \cT_h$ are evolved as $\by_i(t) = \bPhi_t^\cA( \widehat{\by}_i)$. This allows to define evolved elements $ T(t)$ and facets $F(t)$ such that
\begin{equation}\begin{aligned}
    \cT_h(t)=\{T(t)\}, \quad \cF_h(t)=\{F(t)\}, \quad M_{h}(t) = \bigcup T(t),\\
    \cF_h^i(t) = \{F(t)\in \cF_h(t): \; F(t)\cap \partial M_h(t) = \emptyset\}, \\
    \cF_h^b(t) = \{F(t)\in \cF_h(t): \; F(t)\cap \partial M_h(t) = F(t)\}.
\end{aligned}\end{equation}
For the discretization of the function spaces,  we start with the finite element space defined in \eqref{eq:fem_space}, then transport the basis functions using the flow map $\phi_i(t) = \bPhi_t^\cA(\phi_i)$ and define the evolved finite element space as their span
\begin{equation}
    V_h(t) = \{ \mathrm{span}\{\phi_i(t)\}: i = 1,\ldots, N. \; \phi_i \in V_h(\widehat{M}_h)\}.
    \label{eq:evolved_fes}
\end{equation}
The space $V_{h0}(t)$ is defined analogously. It follows from the definitions that $V_h(0) \equiv V_h(\widehat{M}_h)$. Discrete material velocities are defined by
\begin{equation}
    \bv_h(t, \cdot) = \sum_{i=1}^N \bv(t, \by_i(t))\,\phi_i(t), \quad
    \bv^\cA_h(t, \cdot) = \sum_{i=1}^N \bv^\cA(t, \by_i(t))\,\phi_i(t),
\end{equation}
with corresponding derivatives
\begin{equation}
    \matder_h \phi_i = \partial_t\bar{\phi_i} + \bv_h\cdot\nabla\bar{\phi_i}, \quad \aleder_h \phi_i = \partial_t\bar{\phi_i} + \bv^\cA_h\cdot\nabla\bar{\phi_i}.
\end{equation}
Using an ALE map, which describes the motion of the computational domain, it holds that $\aleder_h \phi_i = 0$.

In this work, time-dependent PDEs are discretized in a collection of time points $\{t_0,\ldots, t_N\}$. We denote quantities on discrete time points by a superscript, where $(\cdot)^{n}$ refers to the previous time step and $(\cdot)^{n+1}$ refers to the current time step. The time step size is defined as $\tau^{n+1} = t_{n+1}-t_n$ and for the sake of simplicity we restrict in what follows to constant time step size $\tau$, with the understanding that every first-order-in-time scheme is amenable to adaptive time stepping.

\begin{rmk}
    From now on, we will often omit the explicit dependence on time, indicating for example $M_h(t)$ as simply $M_h$ and leaving it to the context to indicate the considered time.
\end{rmk}

\section{Structure preservation}
\label{sec:structure_pres}
In this section we describe how to deal with the second one of the challenges highlighted in subsection \ref{subsec:background_and_soa}, i.e. positivity, bounds and mass preservation in the interest of interpretability.  The method recently proposed in \cite{ChengShen2022, ChengShen2022a} provides a flexible approach to achieve such goal for a wide variety of finite element formulations on stationary bulk domains. Here, we extend those advances to moving surfaces and moving bulk domains. The main advantages of the employed method are:
\begin{itemize}
    \item Bounds preservation is achieved in such a way that, in its simplest form, it reduces to the widely used cut-off approach.
    \item Accuracy is theoretically retained for higher order methods in space and time, granted some underlying hypotheses are satisfied. This is not the case, for example, for discrete maximum principle preserving schemes.
    \item The approach can easily be incorporated as a post-processing step without requiring explicit time integration. While explicit time integration as reviewed in \cite{ZhangShu2011} can be effective for hyperbolic equations, here we are also interested in parabolic equations.
    \item The simplicity of the algorithm makes it suitable for implementation in legacy codes with negligible overhead computational time with respect to the original non-preserving discretization method. Nevertheless, the method enjoys good stability properties and even allows for error analysis for some specific cases, see \cite{ChengShen2022a}.
    \item No \emph{ad hoc} problem reformulation is needed, but the scheme only acts on the nodal values of the solution through a Lagrange multiplier approach.
\end{itemize}
Following the presentation in \cite{ChengShen2022a}, the derivation of the method starts from a nonlinear PDE in the form
\begin{equation}
    u_t+\cL u + \cN(u) = 0,
\end{equation}
where $\cL$ is a linear or non-negative operator and $\cN(u)$ is a semilinear or quasi-linear operator. The operators $\cL$ and $\cN(u)$ are acting on functions defined on the space-time manifold $\cG_T$. As underlying assumption, the solution to the continuum problem lies in the interval $[a,b]$, i.e. given $a\leq u(0, \bp)\leq b$ for all $\bp\in \widehat{M}$ then $a\leq u(t, \bx)\leq b$ for all $(t, \bx) \in \cG_T$.
The corresponding generic spatial discretization is
\begin{equation}
    \partial_t u_h + \cL_h u_h + \cN_h (u_h) = 0,
\end{equation}
where we assume $u_h \in V_h$. To make the scheme bound-preserving, a Lagrange multiplier $\lambda_h$ is introduced together with the quadratic function $g(u) = (b-u)(u-a)$. The problem is then reformulated as
\begin{subequations}
\label{eq:bp_discrete_general}
\begin{align}
    \partial_t u_h + \cL_h u_h + \cN_h (u_h) = \lambda_h g'(u_h), \\
    \lambda_h\geq 0, \quad g(u_h)\geq 0, \quad \lambda_h g(u_h) =0,
\end{align}
\end{subequations}
where the second line represents the usual Karush-Kuhn-Tucker conditions for constrained optimization \cite{Karush_book, KuhnTucker_bookSection}. A core assumption of the scheme presented in \cite{ChengShen2022a} is that \eqref{eq:bp_discrete_general} is satisfied \emph{pointwise} in a set of points $\Sigma_h$ that can be both mesh points or collocation points for $M_h$. It must then hold:
\begin{subequations}
\label{eq:bp_discrete_general_poitwise}
\begin{align}
    \partial_t u_h(\bp) + \cL_h u_h(\bp) + \cN_h (u_h(\bp)) = \lambda_h(\bp) g'(u_h(\bp)),\; \forall \bp\in \Sigma_h, \\
    \lambda_h(\bp)\geq 0, \quad g(u_h(\bp))\geq 0, \quad \lambda_h(\bp) g(u_h(\bp)) =0, \; \forall \bp\in \Sigma_h.
\end{align}
\label{eq:struct_pres_kkt}
\end{subequations}
It is important to notice that $\Sigma_h$ does not include points where essential boundary conditions are applied.
The next step is to discretize \eqref{eq:struct_pres_kkt} in time, where we here focus on the case of
Backward-Euler (BE or BDF-1) time integration in the interval $[t^n, t^{n+1}]$. Higher order extensions are presented in \cite{ChengShen2022a} and can be adapted following similar arguments. The key idea is to apply an operator splitting approach to \eqref{eq:bp_discrete_general} dividing the problem in two steps.

In the \emph{predictor} step we solve the unconstrained problem 
\begin{equation}
    \frac{\tilde{u}_h^{n+1}(\bp)-u_h^{n}(\bp)}{\tau} + \cL_h \tilde{u}_h^{n+1}(\bp) + \cN_h (\tilde{u}_h^{n+1}(\bp)) = 0, \; \forall \bp\in \Sigma_h.
\end{equation}

In the \emph{corrector} step we solve the constrained problem
\begin{subequations}
\label{eq:bp_discrete_corrector}
    \begin{align}
        \frac{u_h^{n+1}(\bp)-\tilde{u}_h^{n+1}(\bp)}{\tau} = \lambda_h^{n+1}(\bp) g'(u_h^{n+1}(\bp)),\; \forall \bp\in \Sigma_h, \\
    \lambda_h^{n+1}(\bp)\geq 0, \quad g(u_h^{n+1}(\bp))\geq 0, \quad \lambda_h^{n+1}(\bp) g(u_h^{n+1}(\bp)) =0, \; \forall \bp\in \Sigma_h.
    \end{align}
\end{subequations}
If one is only interested in bound preservation, it turns out that the solution to \eqref{eq:bp_discrete_corrector} is the cutoff function
\begin{equation}
    u_h^{n+1} = \text{cutoff}[\tilde{u}_h^{n+1}]  =
    \begin{cases}
        \tilde{u}_h^{n+1}(\bp), &\text{ if } a<\tilde{u}_h^{n+1}(\bp)<b,\\
        a, &\text{ if } \tilde{u}_h^{n+1}(\bp)<a, \\
        b, &\text{ if } \tilde{u}_h^{n+1}(\bp)>b.
    \end{cases}
\end{equation}
The above step is unfortunately not mass-preserving even if the solution is. To fix this, the authors in \cite{ChengShen2022,ChengShen2022a} introduce an additional global space-independent Lagrange multiplier $\xi_h^{n+1}$. The corrector step \eqref{eq:bp_discrete_corrector} is modified to include the mass preservation constraint,
\begin{subequations}
\label{eq:bp_discrete_corrector_mp}
    \begin{align}
        \frac{1}{\tau}u_h^{n+1}(\bp) = \lambda_h^{n+1}(\bp) g'(u_h^{n+1}(\bp)) + \xi_h^{n+1} + \frac{1}{\tau}\tilde{u}_h^{n+1}(\bp),\; \forall \bp\in \Sigma_h, \label{eq:bp_discrete_corrector_mp_1}\\
    \lambda_h^{n+1}(\bp)\geq 0, \quad g(u_h^{n+1}(\bp))\geq 0, \quad \lambda_h^{n+1}(\bp) g(u_h^{n+1}(\bp)) =0, \; \forall \bp\in \Sigma_h, \\
    \inner{u_h^{n+1}}{1}^h = \inner{u_h^{n}}{1}^h,
    \end{align}
\end{subequations}
where $\inner{\cdot}{\cdot}^h$ is a suitably chosen discrete inner product. In our case this discrete inner product, that has to be expressed pointwise as $\inner{u}{v}^h=\sum_{\bp\in \Sigma_h}\omega_\bp u(\bp)v(\bp)$, is assumed to be equivalent to the classical $L^2$ inner product on $V_h$. For this reason, whenever required,  $\inner{\cdot}{\cdot}^h$ is chosen to be the mass-lumped inner product.  Rearranging \eqref{eq:bp_discrete_corrector_mp_1} as 
\begin{align}
    \frac{1}{\tau}(u_h^{n+1}(\bp) - (\tilde{u}_h^{n+1}(\bp) + \tau \xi_h^{n+1} )) = \lambda_h^{n+1}(\bp) g'(u_h^{n+1}(\bp)),\; \forall \bp\in \Sigma_h,
\end{align}
the solution is then $u_h^{n+1} = \text{cutoff}[\tilde{u}_h^{n+1} + \tau \xi_h^{n+1}]$ where $\xi_h^{n+1}$ is the solution to the nonlinear equation
\begin{equation}
    F^{n+1}(\xi_h^{n+1}) = \inner{\text{cutoff}[\tilde{u}_h^{n+1}+\tau \xi_h^{n+1}]}{1}^h - \inner{u_h^{n}}{1}^h = 0.
    \label{eq:sturct_pres_nonlin_eq}
\end{equation}
Bounds on the energy norm are provided for a class of non-linear PDEs in \cite{ChengShen2022a} together with error analysis of a restricted class of problems.
\begin{rmk}
    Since $(F^{n+1})'$ might not exist, the authors of the scheme in \cite{ChengShen2022,ChengShen2022a} suggest the use of a secant method to find the solution of $F^{n+1}(\xi) = 0$ in \eqref{eq:sturct_pres_nonlin_eq} using
\begin{equation}
    \xi_{k+1} = \xi_{k}-\frac{F^{n+1}(\xi_k)(\xi_k - \xi_{k-1})}{F^{n+1}(\xi_k)-F^{n+1}(\xi_{k-1})},
\end{equation}
together with initial guesses $\xi_0=0$ and $\xi_1 = \cO(\tau)$. From a practical standpoint, the algorithm converged in less than four iterations for each timestep for the case at hand.
\end{rmk}

The technique introduced is clearly directly applicable to a wide variety of equations. The authors themselves in \cite{ChengShen2022,ChengShen2022a} have proposed experiments for Allen-Cahn, Cahn-Hilliard with variable mobility and Fokker-Planck equations on a fixed bulk domain. Inspired by the biophysical mechanisms driving cell reshaping, we will demonstrate their accuracy when modeling complex advection-diffusion-reaction equations and phase-separation phenomena. The structure preservation properties will allow to simultaneously maintain stability, interpretability and physical complexity.

\section{Tangential grid control and ALE discretization}
\label{sec:mesh_redistribution}
When dealing with evolving domains in a fitted framework, it is not uncommon to incur deformations that lead to highly distorted meshes, as mentioned in subsection \ref{subsec:background_and_soa}. This is especially true when gradient flows are involved, since they only prescribe the velocity and hence the displacement in normal direction. Following \eqref{eq:ale_map}, we introduce a particular ALE map tailored for gradient flow dynamics that works as follows.

Consider the time interval $[t_n, t_{n+1}]$, and suppose the surface $\Gamma_h^n$ would be displaced, following the physics of the problem, to an updated surface $\widetilde{\Gamma}_h^{n+1}$ at time $t_{n+1}$ which is highly distorted. We define an artificial tangential motion aimed at maintaining the surface shape (and consequently energy) but redistributing the nodes more favorably to reduce mesh distortions. We choose the two-stage algorithm proposed in \cite{DuanLi2024}  for the following reasons:
\begin{itemize}
    \item The method has proven to be very effective for gradient flows such as the mean curvature flow and the surface diffusion flow.
    \item Its two-stage nature makes it a tunable feature of the code, that can be turned on and off at will.
    \item Being a two-stage process and not embedded in the gradient-flow algorithm, it can be applied on the mesh evolution independently of the presence of the gradient flow.
    \item Beyond the scope of the current work, the method can potentially be applied to higher order in time integration schemes.
\end{itemize}
The algorithm presented in \cite{DuanLi2024} is based on the assumption that the starting mesh is regular enough. The idea is then to impose an artificial tangential velocity that requires $\bPhi^\cA(t, \cdot):\widehat{\Gamma}\to \Gamma$ to be an harmonic map, i.e.
\begin{equation}
    \begin{cases}
        \frac{\partial\bPhi^\cA}{\partial t}\cdot (\bn_\Gamma\circ\bPhi^\cA) = (\bv\cdot\bn)\circ \bPhi^\cA, \\
        -\Delta_{\,\widehat{\Gamma}}\bPhi^\cA = (\varkappa\bn)\circ\bPhi^\cA.
    \end{cases}
\end{equation}
Equivalently, $\bPhi^\cA$ is required to minimize the energy $\int_{\,\widehat{\Gamma}}|\nabla_{\widehat{M}}\bPhi^\cA(t, \cdot)|^2$, under the constraint $\bv\cdot\bn_\Gamma=0$ on $\Gamma$,  which is imposed using the scalar-valued Lagrange multiplier $\varkappa$. 
The fully discrete version of the scheme reads: given a surface $\widetilde{\Gamma}_h^{n+1}$ evolved from an initial reference surface $\widehat{\Gamma}_h$, find $\bw_h^{n+1}\in [V_{h0}(\widetilde{\Gamma}^{n+1}_h)]^d$ and $\varkappa_h^{n+1}\in V_{h0}(\widetilde{\Gamma}^{n+1}_h)$ such that
\begin{subequations}
    \begin{align}
        &\inner{\nabla_{\widehat{\Gamma}_h}(\bw_h^{n+1}\circ \bPhi^\cA_{t_{n+1}})}{\nabla_{\widehat{\Gamma}_h}(\bphi_h\circ \bPhi^\cA_{t_{n+1}})}_{\widehat{\Gamma}_h} - \inner{\varkappa_h^{n+1}\bn_{\widetilde{\Gamma}_h^{n+1}}}{\bphi_h}_{\widetilde{\Gamma}_h^{n+1}} \nonumber\\
        &\qquad=-\inner{\nabla_{\widehat{\Gamma}_h}(\Id_{\widetilde{\Gamma}_h^{n+1}}\circ \bPhi^\cA_{t_{n+1}})}{\nabla_{\widehat{\Gamma}_h}(\bphi_h\circ \bPhi^\cA_{t_{n+1}})}_{\widehat{\Gamma}_h}\\
        &- \inner{\bw_h^{n+1}\cdot\bn_{\widetilde{\Gamma}_h^{n+1}}}{\psi_h}_{\widetilde{\Gamma}_h^{n+1}} = 0,
    \end{align}
    \label{eq:duanli_discrete}
\end{subequations}
for all $\bphi_h\in [V_{h0}(\widetilde{\Gamma}^{n+1}_h)]^d$ and $\psi_h\in V_{h0}(\widetilde{\Gamma}^{n+1}_h)$. 
The function $\bw_h^{n+1}$ represents the displacement on $\widetilde{\Gamma}_h^{n+1}$ that moves tangentially the nodes maintaining good mesh properties. From this we have that $(\bv_h-\bv_h^\cA)^{n+1} = -\bw_h^{n+1}/\tau$.

\begin{rmk}
    The algorithm proposed in \eqref{eq:duanli_discrete} does not constitute the only option available. Since \eqref{eq:duanli_discrete} is a post-processing step, alternative schemes can be used as long as the mesh redistribution is purely tangential  \cite{Sauer2014, HuLi2022}. To highlight the flexibility of our framework we further introduce the scheme in \cite{HuLi2022}, which closely resembles the one presented in \cite{DuanLi2024}.
    The fully discrete version of the modified scheme reads: given a surface $\widetilde{\Gamma}_h^{n+1}$ evolved from a surface $\Gamma^{n}_h$ at the previous timestep, find $\bw_h^{n+1}\in [V_{h0}(\widetilde{\Gamma}^{n+1}_h)]^d$ and $\varkappa_h^{n+1}\in V_{h0}(\widetilde{\Gamma}^{n+1}_h)$ such that
    \begin{subequations}
        \begin{align}
            \frac{1}{\tau}\inner{\nabla_{\Gamma^{n}_h}(\bw_h^{n+1}\circ \bPhi^\cA_{\tau})}{\nabla_{\Gamma^{n}_h}(\bphi_h\circ \bPhi^\cA_{\tau})}_{\Gamma^{n}_h} - \inner{\varkappa_h^{n+1}\bn_{\widetilde{\Gamma}_h^{n+1}}}{\bphi_h}_{\widetilde{\Gamma}_h^{n+1}} &=0\\
            - \inner{\bw_h^{n+1}\cdot\bn_{\widetilde{\Gamma}_h^{n+1}}}{\psi_h}_{\widetilde{\Gamma}_h^{n+1}} &= 0,
        \end{align}
    \label{eq:huli_discrete}
    \end{subequations}
    for all $\bphi_h\in [V_{h0}(\widetilde{\Gamma}^{n+1}_h)]^d$ and $\psi_h\in V_{h0}(\widetilde{\Gamma}^{n+1}_h)$, where $\bPhi^\cA_{\tau}=\bPhi^\cA_{t_{n+1}}\circ\bPhi^\cA_{-t_{n}}$.
\end{rmk}

After the surface has been displaced, a second step might be needed to extend the ALE motion in the bulk. The bulk mesh is advected through a continuous harmonic extension of $\bw_h^{n+1}$. The Laplace problem reads correspondingly
\begin{align}
    \Delta \be_{h}^{n+1} &= 0,  \qquad \text{ in } \widetilde{\Omega}^{n+1}_h, \\
    \be_{h}^{n+1} &= \bw_h^{n+1},  \quad\text{ on } \partial\widetilde{\Omega}^{n+1}_h,
\end{align}
and we have that $(\bv_h-\bv_h^\cA)^{n+1} = -\be_{h}^{n+1}/\tau$.

\section{Advection-diffusion-reaction system}
\label{sec:adr}
The model equation we consider reads
\begin{equation}
\label{eq:adr_moving}
\begin{cases}
    \partial_tu+\nabla_M\cdot(\bb_M u - \bbA_M \nabla_M u) + c_Mu=f\; &\text{ on }M(t), \\
    u(0, \bx) = u_0 \; &\text{ on } \widehat{M},
\end{cases}
\end{equation}
where $u: \cG_T\to \bbR$ is the concentration, $\bb_M$ is the advective velocity field, $\bbA_M$ is the diffusion matrix and $c_M>0$ is the reaction constant \cite{Gastaldi2001, DziukElliott2013, ElliottRanner2021}.
As in \cite[p.304]{DziukElliott2013} it is assumed that $(\bbA_M)_{ij}, \;(\bb_M)_i,\;\nabla_M\cdot\bb_M,\; c_M \in L^\infty(M),$
together with $\bb_M(\bx)\in T_{\bx}M, \; \forall \bx\in M$ and the requirement for $\bbA_M$ to be symmetric and to map the tangent space $T_{\bx}M$ onto itself.

Extensive literature exists on finite element discretizations of such equations for what concerns the bulk case, see \cite{ErnGuermond_booka,ErnGuermond_bookb}. For \emph{closed surfaces} the SFEM discretization of parabolic PDEs in the form of \eqref{eq:adr_moving} has been comprehensively reviewed in \cite{DziukElliott2013}. Convergence estimates for its continuous-in-time formulation using general order polynomials are found in e.g. \cite{ElliottRanner2021}. For the sake of generality,  here we focus on \emph{moving surfaces with boundary}. The problem reads: Find $u_h\in V_h(\Gamma_h)$ such that
\begin{equation}\begin{aligned}
\label{eq:adr_weak_boundary}
    \frac{\d}{\d t}m_h\inner{u_h}{v_h} + \overline{a}_h\inner{u_h}{v_h} + \inner{(\bb_{rel}u_h-\bbA_{\Gamma_h}\nabla_{\Gamma_h} u_h)\cdot\bn_{\partial\Gamma_h}}{v_h} _{\partial\Gamma_h}\\
     =\overline{l}_h(v_h)\;, \forall v_h \in V_h(\Gamma_h).
\end{aligned}\end{equation}
where
\begin{equation}\begin{aligned}
    m_h(u_h, v_h) =& ~\inner{u_h}{v_h}_{\Gamma_h}, \\
    \overline{a}_h(u_h, v_h) =& ~\inner{\bbA_{M_h}\nabla_{\Gamma_h} u_h}{\nabla_{\Gamma_h} v_h}_{\Gamma_h} + \inner{-\bb_{rel} u_h}{\nabla_{\Gamma_h} v_h}_{\Gamma_h}\\
    &+ \inner{c_{M_h}u_h}{v_h}_{\Gamma_h}, \\
    \overline{l}_h(v_h) =&~ \inner{f}{v_h}_{\Gamma_h}.
\end{aligned}\end{equation}
The coefficients $\bbA_{\Gamma_h},\; \bb_{\Gamma_h},\; c_{\Gamma_h}$ are required to satisfy the necessary smoothness conditions element-wise on the discrete surface $\Gamma_h$. Recalling that  $\bv_h^{\cA}$ is the discrete velocity associated with the ALE mapping, the \emph{relative velocity} is defined as  $\bb_{rel} = \bb_{\Gamma_h} -\bv^\cA_h$.
All boundary conditions are imposed weakly and the boundary is divided based on the type of boundary condition applied. We define as $\partial\Gamma_{h,D}$ the part of the boundary where Dirichlet conditions are applied and as $\partial\Gamma_{h,N}$ the part where Neumann boundary conditions are applied. It is required that
\begin{equation}
    \partial\Gamma_{h,D}\cap\partial\Gamma_{h,N} = \emptyset, \quad \partial\Gamma_{h,D}\cup\partial\Gamma_{h,N}=\partial\Gamma_{h}.
\end{equation}
Nitsche's boundary penalty method is used in order to impose the inhomogeneous Dirichlet boundary condition $u_h(t, \bx)|_{\partial\Gamma_{h,D}}=u(t,\bx)$ weakly for the diffusive term \cite{Nitsche1971, ErnGuermond_bookb}. For convenience, we define $v^+ = \max\{v, 0\}$ and $ v^- = \min\{v, 0\} $. The resulting bilinear forms are given by
\begin{equation}\begin{aligned}
    a_h\inner{u_h}{v_h} =&~ \overline{a}_h\inner{u_h}{v_h} - \inner{(\bbA_{\Gamma_h}\nabla_{\Gamma_h} u_h)\cdot\bn_{\partial\Gamma_h}}{v_h}_{\partial\Gamma_{h,D}}\\
    & - \inner{u_h}{(\bbA_{\Gamma_h}\nabla_{\Gamma_h} v_h)\cdot\bn_{\partial\Gamma_h}}_{\partial\Gamma_{h,D}}+ \inner{\gamma_{\bbA }h^{-1}_F u_h}{v_h}_{\partial\Gamma_{h,D}}\\
    &+ \inner{(\bb_{rel}\cdot \bn_{\partial\Gamma_h})^+u_h}{v_h}_{\partial\Gamma_{h}}\\
    l_h(v_h) =&~ \overline{l}_h(v_h)- \inner{u}{(\bbA_{\Gamma_h}\nabla_{\Gamma_h} v_h)\cdot\bn_{\partial\Gamma_h}}_{\partial\Gamma_{h,D}}+  \inner{\gamma_{\bbA}h^{-1}_F u}{v_h}_{\partial\Gamma_{h,D}}\\
    &+ \inner{(\bbA_{\Gamma_h}\nabla_{\Gamma_{h,N}} u)\cdot\bn_{\partial\Gamma_h}}{v_h}_{\partial\Gamma_{h,N}} - \inner{(\bb_{rel}\cdot \bn_{\partial\Gamma_h})^-u}{v_h}_{\partial\Gamma_{h}},
\end{aligned}\end{equation}
where $\gamma_\bbA>0$ is the \emph{penalty parameter} enforcing $u_h=u$ on $\partial\Gamma_{h,D}$.
Overall the problem for open surfaces reads: find $u_h\in V_h$ such that
\begin{equation}
\label{eq:adr_weak_nostab}
    \frac{\d}{\d t}m_h\inner{u_h}{v_h} + a_h\inner{u_h}{v_h}=l_h(v_h),\; \forall v_h \in V_h(\Gamma_h).
\end{equation}

\subsection{Numerical methods for advection-dominant problems}
\label{subsec:adr_adv_dom}

In biophysical cell models, as introduced in subsection \ref{subsec:background_and_soa}, it is not uncommon to incur advection-dominant ADR equations. In our case, even purely parabolic equations might become advection-dominant due to the introduction of an ALE velocity $\bv^\cA$. 
In the present work we employ and extend the continuous interior penalty (CIP) stabilization to cope with the possibly dominant advection regime. The CIP method was proposed and analyzed in~\cite{BurmanHansbo2004,BurmanFernandez2009a} and then later extended to the case of stationary surfaces in~\cite{BurmanHansboLarsonEtAl2018b}. Moreover, the CutFEM approach for moving surfaces developed in \cite{HansboLarsonZahedi2015} applies a CIP-type stabilization for advection dominant problems on moving surfaces. We chose CIP for the following reasons:
\begin{itemize}
    \item It is a widely known, easily implementable technique whose implementation tools are usually shipped in classical finite element packages;
    \item Although only weakly consistent, it has proven to be very successful in time-dependent problems given its commutativity with the time derivative \cite{ErnGuermond_book};
    \item Therefore, it has been shown to lead to convergent algorithms for discretizations of bulk, surface and implicitly described moving surfaces problems.
\end{itemize}
To formulate the CIP stabilization, we need to define averages and jumps of functions across edges and faces. 
For a piecewise discontinuous function $f$ defined on
a surface or bulk mesh $\mathcal{T}_h(M_h)$,
we define its average and jump over an interior facet $F \in \mathcal{F}_h^i(M_h)$ by
\begin{equation}
\{f\}|_F  =\frac{1}{2}\left(f_F^{+}+f_F^{-}\right), \quad
{\jump{f}|_F }  =f_F^{+}-f_F^{-},
\end{equation}
respectively, where 
$
f_F^{\pm}(\bx) = \lim_{\delta \to 0^+} f(\bx_F^{\pm} - \delta \bn_F^{\pm})
$.

The main idea of the CIP approach is to penalize the jump of the streamline derivative across element interfaces. 
To this end, the CIP stabilization form is defined as follows,
\begin{equation}
    s_{h}^{b}\left({u}_{h}, {v}_{h}\right)
    :=\gamma_{b} \sum_{F \in \mathcal{F}_{h}^{i}}
    \phi_{F} h_F\left(\jump{ {\boldsymbol{b}_{M_h}} \cdot \nabla_{M_h} {u}_{h} }, \jump{\boldsymbol{b}_{M_h} \cdot \nabla_{M_h} {v}_{h} }\right)_{F},
    \label{eq:cip_stab_b-I}
\end{equation}
where $\gamma_b>0$ is a dimensionless stabilization constant and
$\phi_{F}$ denotes a stabilization parameter defined by
\begin{equation}
    \label{eq:phi_b}
    \phi_F = \max(\phi_{T^+}, \phi_{T^-}), \quad \phi_{T^{\pm}} =\min(\beta_{T^{\pm}}^{-1}h_{T^{\pm}}, \mu_{T^{\pm}}^{-1}), 
    \quad T^{\pm}\in \cT_h. 
\end{equation}
In the above $\beta_T$ is a local velocity scale and $\mu_T$ is the reciprocal of a time \cite{ErnGuermond_book}. If $\bb_{M_h}\in L^{\infty}(M_h)\cap C^0(M_h)$, then \eqref{eq:cip_stab_b-I} can be simplified to
\begin{equation}
    s_{h}^{b}({u}_{h}, {v}_{h})
    :=\gamma_{b} \sum_{F \in \mathcal{F}_{h}^{i}}
    \beta_Fh_F^2\left(\jump{ \nabla_{M_h} {u}_{h} }, \jump{\nabla_{M_h} {v}_{h} }\right)_{F}.
    \label{eq:cip_stab_b-II}
\end{equation}
with $\beta_F=\max (\beta_{T^+}, \beta_{T^-})$. We augment the system in \eqref{eq:adr_weak_nostab} further to handle advection-dominant problems, which now reads: find $u_h\in V_h$ such that
\begin{equation}
\label{eq:adr_weak_cip}
    \frac{\d}{\d t}m_h\inner{u_h}{v_h} + a_h\inner{u_h}{v_h} +  s_h^{b}\inner{u_h}{v_h}=l_h(v_h),\; \forall v_h \in V_h(\Gamma_h).
\end{equation}

\subsection{Numerical results for the ADR system}
Given the novelty of the methods introduced, little to no theoretical convergence analysis is available in the literature and is left for future work. We thus proceed in reporting experimental convergence studies using select examples of increasing complexity to testify the accuracy of the framework.

\begin{figure}[tbhp!]
    \centering

    \begin{subfigure}[tbhp!]{\textwidth}
        \centering
        \includegraphics[width=\textwidth, clip, trim= {0cm 7.4cm 0cm 13.5cm}, page = 2]{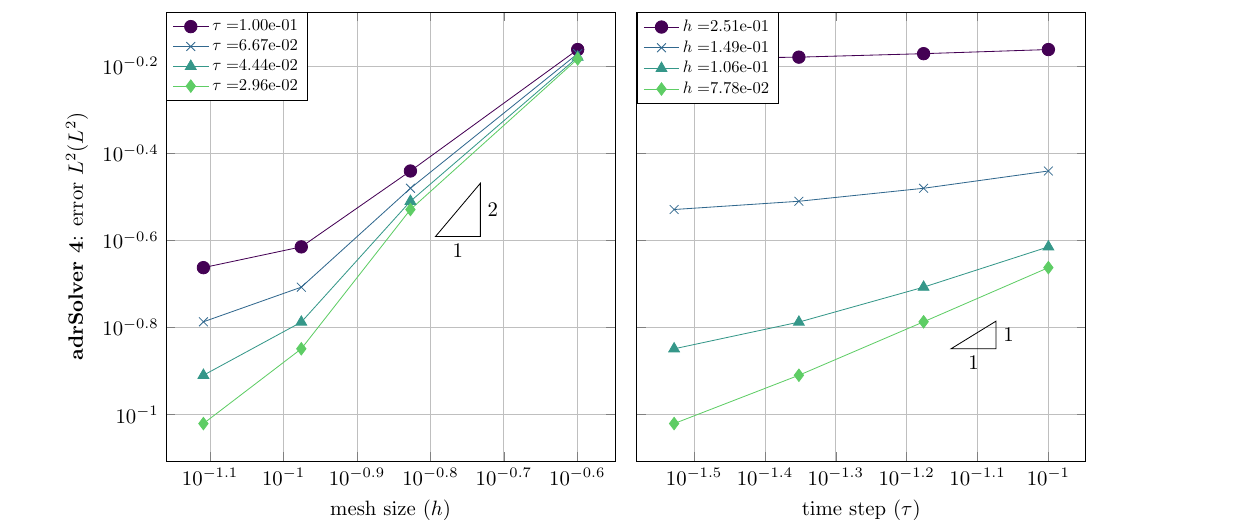}
    \end{subfigure}%

    \caption{Convergence studies for the solver \textbf{adrSolver 4} in List \ref{item:ad_solvers}. As expected, first and second order convergence are achieved in time and space, respectively.}
    \label{fig:adr_conv_open}
\end{figure}

We begin by considering a domain in 3D that moves under the linear transformation $\bPhi(t, \bp) = \bPhi^\cA(t, \bp) = \bbA(t)\bp + B(t)$ where
\begin{equation}
    \bbA(t) = 
    \begin{bmatrix}
        \cos(t) & -\sin(t) & 0 \\
        \sin(t) & \cos(t) & 0 \\
        0 & 0 & 1
    \end{bmatrix},
    \quad B(t)=
    \begin{bmatrix}
        0 \\
        0 \\
        0
    \end{bmatrix}.
    \label{eq:nr_transformation}
\end{equation}
We have that $\bPhi^{-1}(t, \bx) = [\bbA(t)]^{-1}(\bx - B(t))$ and that the domain velocity is $\bv(t, \bx) = \bbA(t)'\bPhi^{-1}(t, \bx) + B(t)' \equiv \bv^{\cA}$. We construct manufactured solution for the problem
\begin{equation}
    \label{eq:adv_dom_manufactured}
    \partial_t u + \nabla_{\Gamma}\cdot(\bb_{\Gamma} u) + c_{\Gamma}u = f,
\end{equation}
with
\begin{equation}\begin{aligned}
    u_{ex} & = \cos(2\pi x_1)\cos(t),\\
    c_\Gamma &= 1+t^2, \\
    \bb_\Gamma &= \bv^\cA +\bb_{rel} = \bv^\cA + \bbP_{\Gamma}(2, -x_1,0),
\end{aligned}\end{equation}
where $\bx = (x_1, x_2, x_3)$. The right-hand side $f$ has been chosen so to satisfy the imposed solution. For $u_{ex}$ we have that $\int_{\Gamma} u(t) =0$ and that $u\in [-1, 1]$. The geometry chosen is a half-sphere with unit radius that was cut along the $x_1$-$x_2$ plane. Under the above rotation the half-sphere has normal
$\bn(t, \bx) = \bx/\norm{\bx}$
.
We test convergence properties for the following solvers:
\begin{enumerate}
\label{item:ad_solvers} 
    \item \label{item:ad_solver} \textbf{adrSolver 1}: CIP stabilized solver \eqref{eq:adr_weak_cip}.
    \item \label{item:ad_solver_BP} \textbf{\textbf{adrSolver \ref{item:ad_solver_BP}}}: CIP stabilized solver with imposed bound preservation in the interval $[-1, 1]$.
    \item \label{item:ad_solver_MP} \textbf{\textbf{adrSolver \ref{item:ad_solver_MP}}}: CIP stabilized solver with imposed mass preservation using the lumped-mass inner product $\inner{\cdot}{\cdot}^h$.
    \item \label{item:ad_solver_BP_MP} \textbf{\textbf{adrSolver \ref{item:ad_solver_BP_MP}}}: CIP stabilized solver with both imposed bound preservation in the interval $[-1, 1]$ and mass preservation.
\end{enumerate}
All solvers reveal the expected convergence in time and space.
The results for \textbf{adrSolver 4} with BDF-1 time stepping are shown in \autofigref{fig:adr_conv_open}. Analogous results are obtained for the other solvers.

\begin{figure}[tbhp!]
    \centering
    \begin{subfigure}[tbhp!]{0.45\textwidth}
        \centering
        \includegraphics[clip, trim = {0cm 10cm 0cm 15cm}, width = \textwidth]{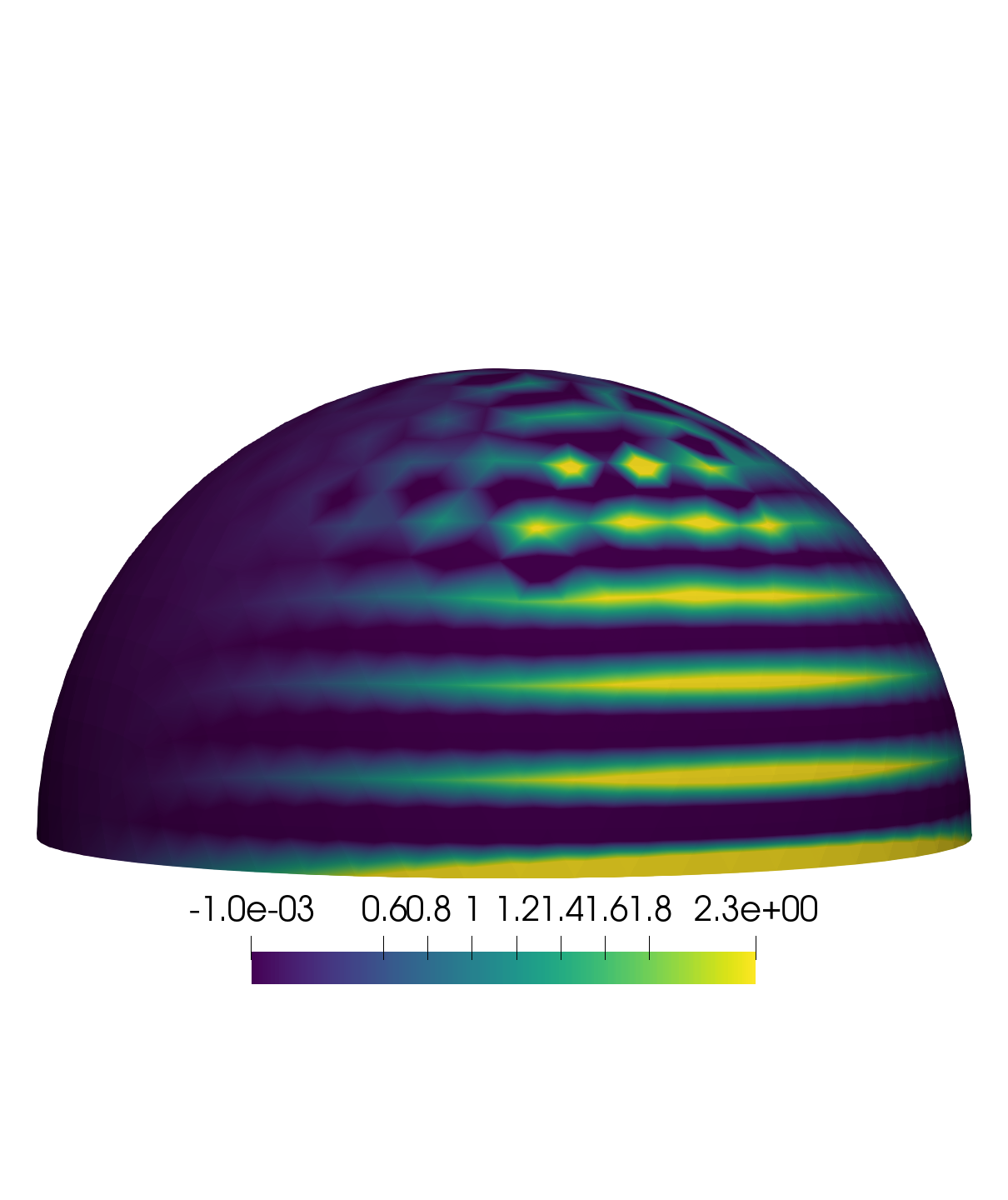}
        \caption{Non-stabilized solver}
        \label{fig:illposed0}
    \end{subfigure}
    ~ 
    \begin{subfigure}[tbhp!]{0.45\textwidth}
        \centering
        \includegraphics[clip, trim = {0cm 10cm 0cm 15cm}, width = \textwidth]{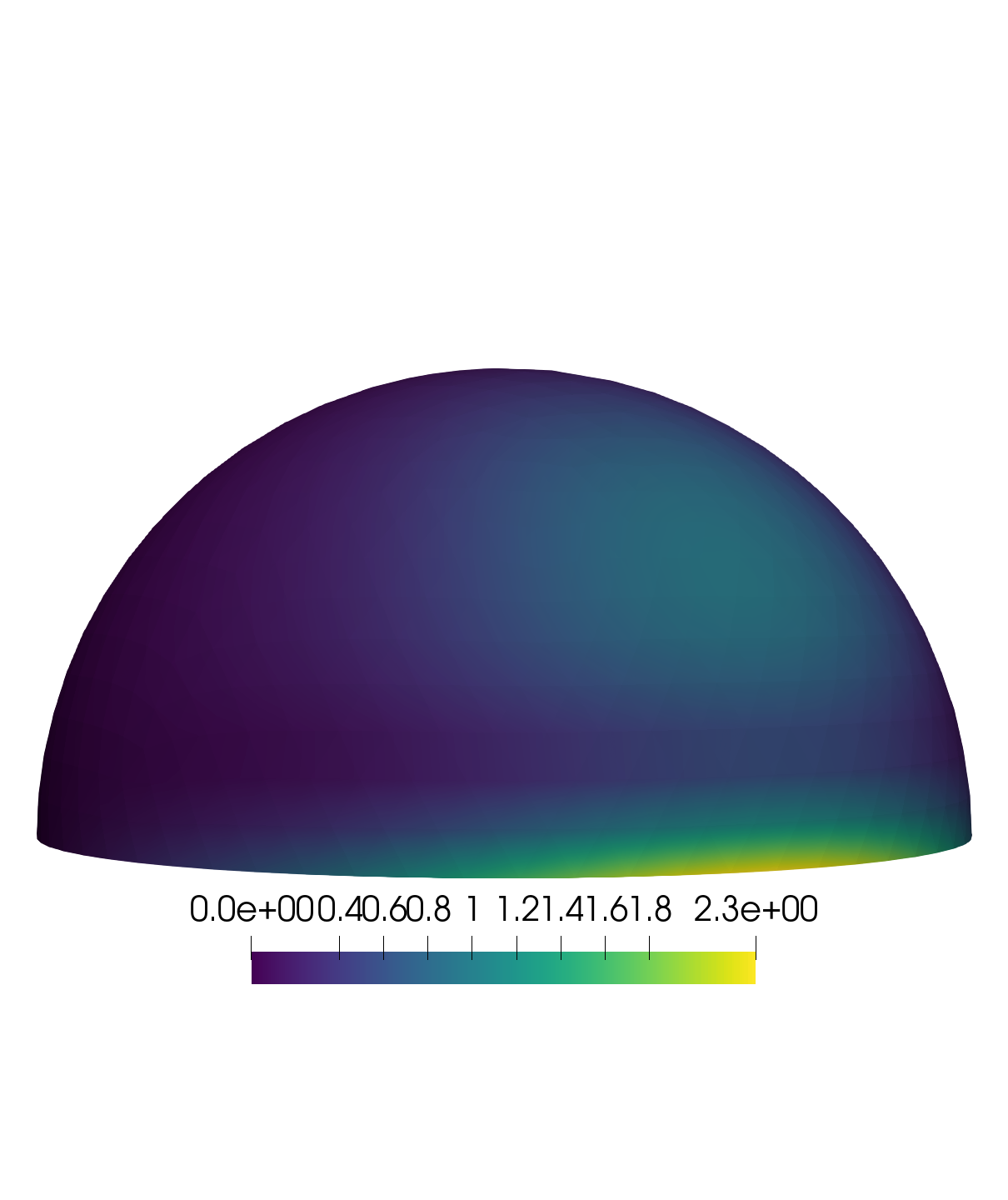}
        \caption{\textbf{adrSolver \ref{item:ad_solver_BP_MP}}}
        \label{fig:illposed4}
    \end{subfigure}
    \caption{(a) Oscillating final solution of non-stabilized ADR solver at $t=1$. (b) Final solution for \textbf{adrSolver \ref{item:ad_solver_BP_MP}} at $t=1$.}
    \label{fig:adr_figs}

    \begin{subfigure}[tbhp!]{\textwidth}
        \centering
        \includegraphics[width=\textwidth, clip, trim=0cm 18.1cm 0cm 2.5cm]{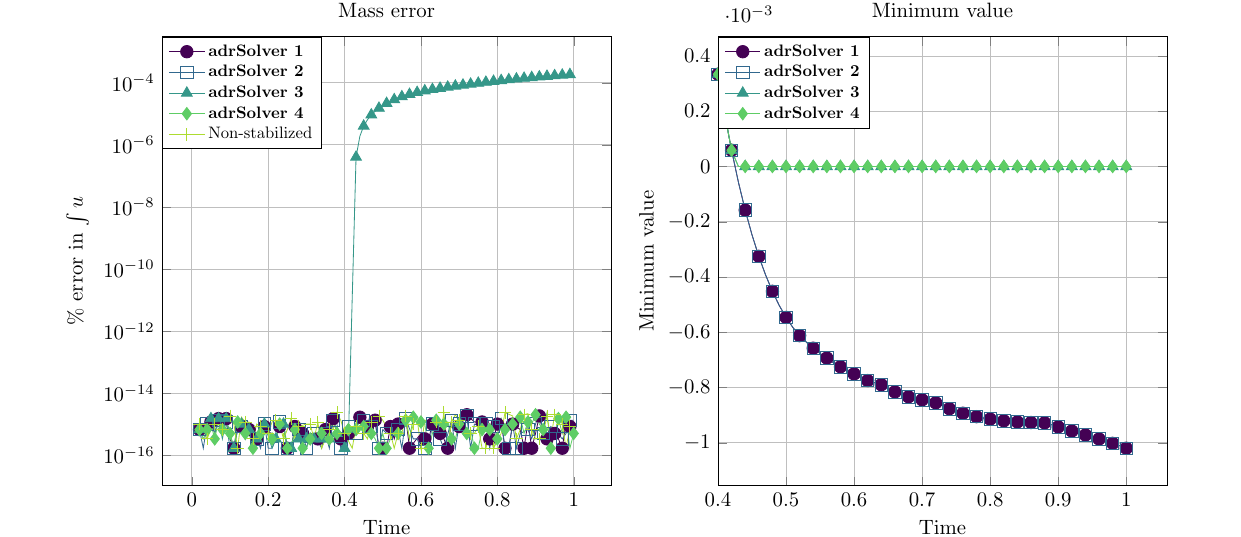}
    \end{subfigure}%
    \caption{Left: relative absolute error in $\int_{\Gamma(t)} u(t)$ with respect to $\int_{\Gamma(0)} u(0)$ for the solvers in List \ref{item:ad_solvers}. Right: Close-up of minimum value of $u$ for the different solvers in List \ref{item:ad_solvers}. }
    \label{fig:adr_bmp}
\end{figure}

To further test the effectiveness of the algorithm we now consider an ill-posed problem. Maintaining the half-sphere geometry and the transformation as in \eqref{eq:nr_transformation}, we solve the pure transport problem 
\begin{equation}
    \partial_t u + \nabla_{\Gamma}(\bb_{\Gamma} u) =0,
\end{equation}
with 
\begin{equation}
    \bb_{\Gamma} = \bv^{\cA} +\bb_{rel} = \bv^{\cA} + (x_3, 0, -x_1)\cdot (1-e^{-10x_3}), \quad u_0=e^{-3(x_1^2+x_2^2)}.
\end{equation}
By construction the initial concentration is pushed towards the boundary where an incredibly sharp boundary layer forms, since zero flux boundary conditions are enforced by the fact that $\bb_{rel}\cdot \bn_{\partial\Gamma}=0$. Given the zero flux condition, the resulting solution is also mass-preserving. Parameters are set as follows
\vspace{2mm}
\begin{center}
\begin{tabular}{||c |c| c| c| c||} 
 \hline
 $h$ & $nv$ & $ne$ & $T$ & $\tau$\\ [0.5ex] 
 \hline\hline
 0.1 & 780 & 1495  & 1 & 0.01\\ 
 \hline
\end{tabular}
\end{center}
\vspace{2mm}
where $nv$ is the number of vertices and $ne$ is the number of elements. We compare the results of the four different solvers of List \ref{item:ad_solvers} in \autofigref{fig:adr_figs} and \autofigref{fig:adr_bmp}.  The bounds for \textbf{adrSolver \ref{item:ad_solver_BP}} and \textbf{adrSolver \ref{item:ad_solver_BP_MP}} were set to $[0, 10^5]$ in order to guarantee positivity.

It can be seen that the non-stabilized solver is unstable due to the pure advective regime, leading to oscillations that inevitably effect the whole domain, see \autofigref{fig:illposed0}. \textbf{adrSolver \ref{item:ad_solver}} effectively controls the solution's gradient thanks to the stabilization, generating a non-oscillating solution, but fails to maintain its positivity, see \autofigref{fig:adr_bmp} (right plot). This characteristic is lost along $\partial\Gamma$ where the species encounters the steep boundary layer. \textbf{adrSolver \ref{item:ad_solver_BP}} manages, thanks to the cutoff, to keep the solution positive but fails in maintaining mass-preservation as it can be seen in \autofigref{fig:adr_bmp} (left plot), where the relative absolute percentage mass error is plotted. \textbf{adrSolver \ref{item:ad_solver_BP_MP}} successfully manages to deal with the boundary layer, maintain positivity,  and avoid oscillations. The final solution is shown in \autofigref{fig:illposed4}. 

\section{Cahn-Hilliard system as a phase separation model}
\label{sec:phase_sep}
Phase field modeling plays a central role in biophysical systems to investigate phase transition processes. Common phase field models for stationary surfaces are derived from Allen-Cahn and Cahn-Hilliard dynamics \cite{CahnHilliard1958, ErlebacherAzizKarmaEtAl2001, ElliottStinner2010, DuJuTian2011}, where the motion follows the minimization of an energy functional. We will focus here on the functional
\begin{equation}
    \cE_{CH}(u) = \int_{M(t)}\sigma\Bigl(\underbrace{\frac{\varepsilon |\nabla_\Gamma u|^2}{2}}_{I} +\underbrace{\frac{1}{\varepsilon}F(u)}_{II}\Bigr),
    \label{eq:cahn_hillard_energy}
\end{equation}
where $u$ is the phase field, $\varepsilon\in \bbR^+$ is a length scale parameter determining the interface length and $\sigma\in\bbR^+$ describes the interface tension between different phases. Term $I$ encodes the interface energy, penalizing inhomogeneities. Term $II$ encodes the local free energy and drives phase separation through the use of a double-well potential. We consider here two different double-well potentials
\begin{subequations}
\begin{align}
    F_1(u) &:= \frac{1}{4}((1-u) \log(1-u) + (1 + u) \log(1 + u)) + \frac{1-u^2}{2}, \\
    F_2(u) &:= \frac{(1-u^2)^2}{4},
\end{align}
\label{eq:ch_potentials}
\end{subequations}
for which
\begin{equation}
    F_1'(u) = \frac{1}{4}\log\Bigl(\frac{1+u}{1-u}\Bigr)-u, \text{ and } F_2'(u) = u^3-u \text{ respectively}.
\end{equation}
$F_1$ is a Flory-Huggins-type mixing free energy density whose formulation finds its roots in configurational entropy \cite{Huggins1941, Flory1942}. It can be seen that $F_1$ is only well posed for $u\in(-1, 1)$. For this reason, a simpler polynomial potential $F_2$ is often used. The advantage of the latter is that it is well defined for values of $u$ outside the interval $(-1,1)$.



The extension of the Cahn-Hilliard equation to \emph{evolving surfaces} has seen recent developments due to its applicability to biological membranes and alloys \cite{ZhiliakovWangQuainiEtAl2021, EilksElliott2008}. The evolving surface Cahn-Hillard equation \cite{CaetanoElliott2021,CaetanoElliottGrasselliEtAl2023} can be written as
\begin{equation}
\label{eq:cahn_hilliard}
\begin{cases}
    \matder u + (\nabla_\Gamma\cdot\bv) u = \nabla_{\Gamma}(m\nabla_\Gamma w), \\
    w = -\sigma\varepsilon\Delta_\Gamma u + \frac{\sigma}{\varepsilon}F'(u),
\end{cases}
\end{equation}
subject to the initial condition $u(0) = u_0$ for suitable initial data. $w$ represents the chemical potential and $m\in\bbR^+$ is the mobility. As explained in \cite{LiuHuangXiaoEtAl2024}, models for evolving surfaces have similar energy functions to models on stationary surfaces, but important differences arise. Even  if the evolution of the systems is driven by the minimization of \eqref{eq:cahn_hillard_energy}, energy dissipation is no longer satisfied since the velocity is viewed as external force that adds to the system. Analogously, mass preservation depends on the velocity. To keep mass conservation, we must assume $\nabla_\Gamma\cdot \bv = 0$, which is equivalent to assume inextensibility of the cell membrane, and by consequence area preservation. This assumption also allows to consistently extend this one-phase model to $N$-phase models satisfying the hyperlink condition, as described in \cite{LiuHuangXiaoEtAl2024}.

\subsection{Numerical methods for Cahn-Hilliard system}
For the discretization of \eqref{eq:cahn_hilliard}, the evolutionary surface finite element method (ESFEM) is one of the most popular methods adopted and the one we consider here \cite{DziukElliott2013}. Its advantage lies in simple implementation, low memory consumption, ability to handle complex deformations and potential for efficient parallelization. We refer to \cite{LiuHuangXiaoEtAl2024} for a review of alternative methods.  We employ a semi-implicit approach where the fully discrete formulation reads: given $u_h^n \in V_h(\Gamma_h^n)$ find $u^{n+1}_h, w^{n+1}_h\in V_h(\Gamma_h^n)$ such that
\begin{subequations}
\label{eq:cahn_hilliard_fully_discrete}
    \begin{align}
    n_h\inner[auto]{\frac{u^{n+1}_h - u_h^n}{\tau}}{\phi_h} +  \inner{(\bv_h -\bv^\cA_h) u^{n+1}_h}{\nabla_{\Gamma^n_h} \phi_h}_{\Gamma^n_h} + mb_h\inner{w^{n+1}_h}{\phi_h}= 0, \\
    n_h\inner{w^{n+1}_h}{\psi_h} - \sigma\varepsilon b_h\inner{u^{n+1}_h}{\psi_h} = \frac{\sigma}{\varepsilon} n_h\inner{F'(u^n_h)}{\psi_h},
\end{align}
\end{subequations}
for all $\phi_h, \psi_h \in V_h(\Gamma_h^n)$ with
\begin{equation}
    n_h\inner{u_h}{\phi_h} = \inner{u_h}{\phi_h}_{\Gamma_h^n}, \quad b_h\inner{u_h}{\psi_h} = \inner{\nabla_{\Gamma_h^n} u_h}{\nabla_{\Gamma_h^n} \psi_h}_{\Gamma_h^n}.
\end{equation}
Similar discretizations can be found in \cite{ElliottStinner2010,BachiniKrauseNitschkeEtAl2023,MokbelMokbelLieseEtAl2024}, where a Cahn-Hilliard solver is coupled with a domain evolution which is unknown. In those works, the potential $F_2$ was used. In our case, we can take advantage of the scheme presented in Section \ref{sec:structure_pres} to employ both $F_1$ and $F_2$, and apply mass conservation if needed. The scheme proceeds as follows.

In the \emph{predictor} step we solve the unconstrained problem \eqref{eq:cahn_hilliard_fully_discrete}.
Taking the non-linear term explicitly ensures $u_h^{n}$ lays inside the desired bound. 

In the \emph{corrector} step we solve the constrained problem imposing $u^{n+1}_h\in (-1, 1)$. 


\subsection{Numerical results for Cahn-Hilliard system}

\begin{figure}[tbhp!]
    \centering

    \begin{subfigure}[tbhp!]{0.95\textwidth}
        \centering
        \includegraphics[clip, trim = {0cm 18.5cm 0cm 2.9cm}, width = \textwidth, page =2]{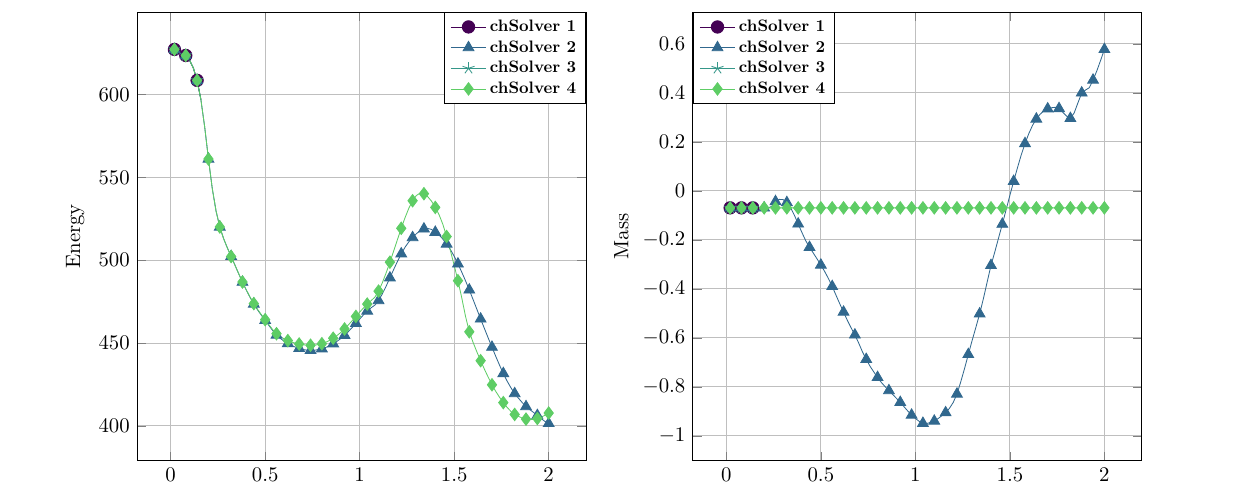}
    \end{subfigure}

    \begin{subfigure}[tbhp!]{0.95\textwidth}
        \centering
        \includegraphics[clip, trim = {0cm 8.7cm 0cm 12.2cm}, width = \textwidth, page = 2]{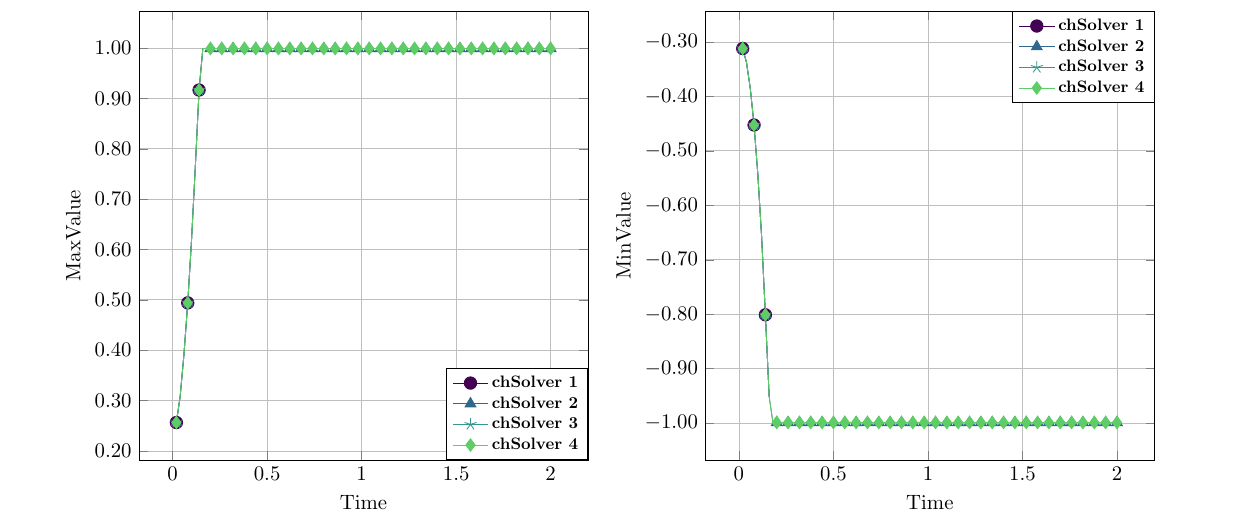}
    \end{subfigure}

    \caption{Experimental studies for potential $F_1(u)=\frac{1}{4}((1-u) \log(1-u) + (1 + u) \log(1 + u)) + \frac{1-u^2}{2}$ and the different solvers in List \ref{item:ch_solvers}.}
    \label{fig:ch_log_studies}
\end{figure}

\begin{figure}[tbhp!]
    \centering
    \begin{subfigure}[tbhp!]{0.95\textwidth}
        \centering
        \includegraphics[clip, trim = {0cm 18.5cm 0cm 2.9cm}, width = \textwidth, page = 2]{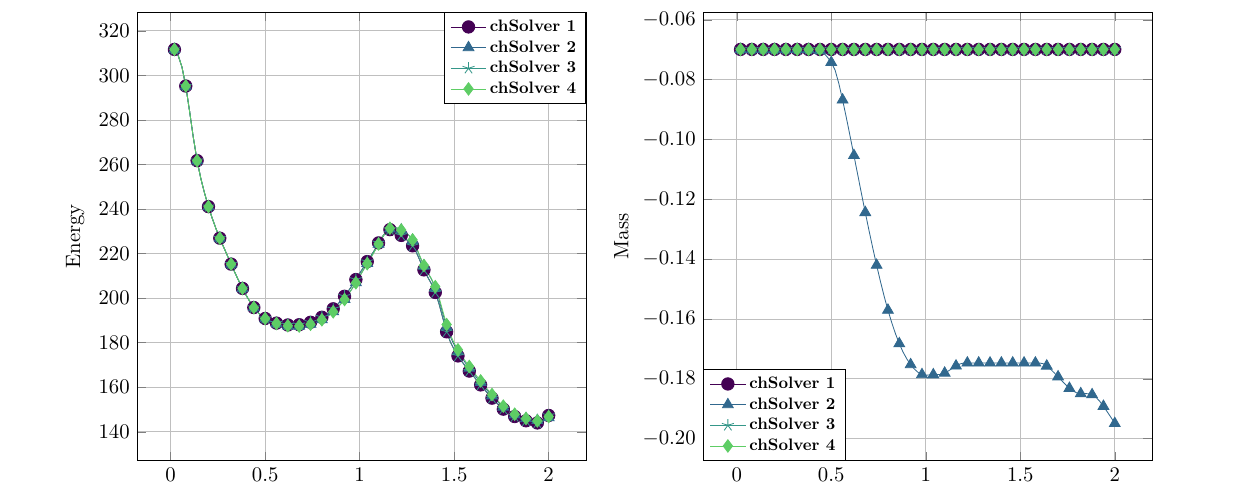}
    \end{subfigure}

    \begin{subfigure}[tbhp!]{0.95\textwidth}
        \centering
        \includegraphics[clip, trim = {0cm 8.7cm 0cm 12.2cm}, width = \textwidth, page = 2]{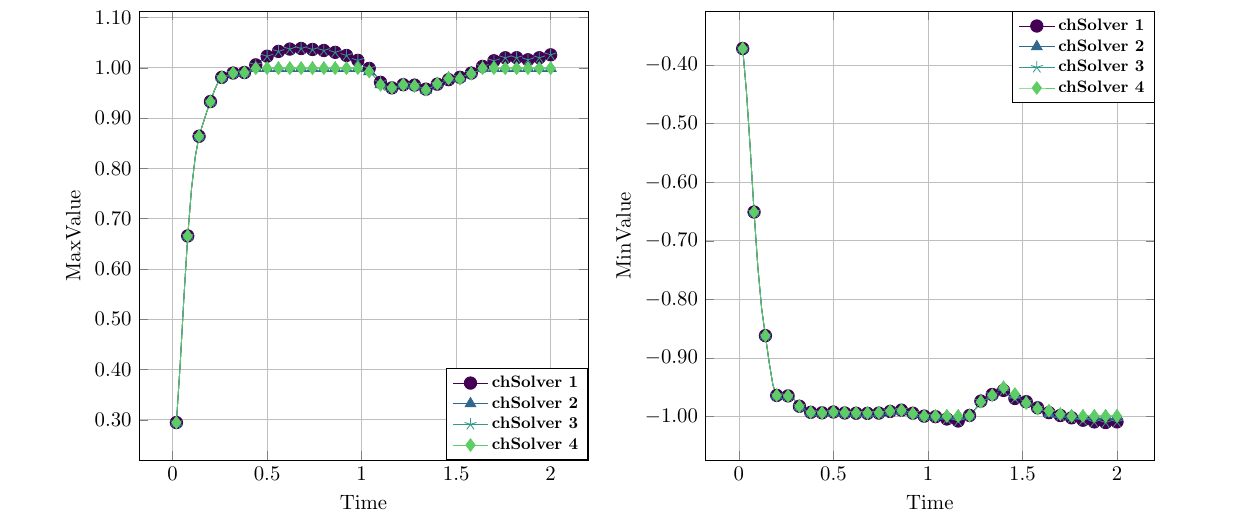}
    \end{subfigure}

    \caption{Experimental studies for potential $F_2(u)=(1-u^2)^2/4$ and the different solvers in List \ref{item:ch_solvers}.}
    \label{fig:ch_pol_studies}
\end{figure}

We proceed in reporting experimental studies to assess the accuracy of the proposed framework.
Consider a cylinder revolving around the $x_3$-axis with radius $R=1$ and height $l=2$. To satisfy the assumption $\nabla_\Gamma\cdot \bv=0$, we transform the initial cylinder through the isometry $\bPhi(t, \bp) = \bPhi^\cA(t, \bp) = \bbA(t)\bp + B(t)$
\begin{equation}
    \bbA(t) = 
    \begin{bmatrix}
        \lambda(t) & 0 & 0 \\
        0 & \lambda(t) & 0 \\
        0 & 0 & 1/\lambda(t)
    \end{bmatrix},
    \quad B(t)=
    \begin{bmatrix}
        0.2t \\
        0.1t \\
        0
    \end{bmatrix},
\end{equation}
where $\lambda(t) = 1+0.5\sin(\pi t)$. 
Convergence studies have already been proposed in the literature and we proceed by testing the quality of the structure preservation technique on realistic initial conditions. The verification is performed for the following solvers:
\begin{enumerate}
\label{item:ch_solvers} 
    \item \label{item:ch_solver} \textbf{chSolver 1}: Solver \eqref{eq:cahn_hilliard_fully_discrete} with no postprocessing.
    \item \label{item:ch_solver_BP} \textbf{chSolver \ref{item:ch_solver_BP}}: Solver \eqref{eq:cahn_hilliard_fully_discrete} with imposed bound preservation in the interval $(-1, 1)$.
    \item \label{item:ch_solver_MP} \textbf{chSolver \ref{item:ch_solver_MP}}: Solver \eqref{eq:cahn_hilliard_fully_discrete} with imposed mass preservation using the lumped-mass inner product $\inner{\cdot}{\cdot}^h$.
    \item \label{item:ch_solver_BP_MP} \textbf{chSolver \ref{item:ch_solver_BP_MP}}: Solver \eqref{eq:cahn_hilliard_fully_discrete} with imposed bound preservation in the interval $(-1, 1)$ and mass preservation.
\end{enumerate}
Zero flux Neumann boundary conditions are imposed for both $u$ and $w$. A random uniform initial distribution in the interval $(-1, 1)$ is chosen for $u_0$. The parameters are set as 
\vspace{2mm}
\begin{center}
\begin{tabular}{||c |c| c| c| c| c||} 
 \hline
 $h$ & $m$ & $\varepsilon$ & $\sigma$ & $T$ & $\tau$\\ [0.5ex] 
 \hline\hline
 0.05 & 0.01 & 0.1 & 10 & 4 & $2\cdot 10^{-3}$\\ 
 \hline
\end{tabular}
\end{center}
\vspace{2mm}
Results for the potential $F_1$ are shown in \autofigref{fig:ch_log_studies} and the same tests for $F_2$ are shown in \autofigref{fig:ch_pol_studies}.
In \autofigref{fig:ch_log_studies}, \textbf{chSolver 1}, \textbf{3} and \textbf{5} break after few steps due to the phase exceeding the bound $(-1, 1)$, while the other solvers correctly reach the end of the simulation. Structure preservation properties are better appreciated in \autofigref{fig:ch_pol_studies}. It can be seen that the energy landscape is not significantly modified by the use of either bound preservation and/or mass preservation. In both mass plots it is evident how the mass preservation property is lost for \textbf{chSolver 2}, that only has bound preservation, leading to a different overall dynamics for \textbf{chSolver 2}. Mass in instead correctly preserved by \textbf{chSolver 4} as expected. Plotting the maximum and the minimum value in the bottom two pictures not only highlights the bound preservation but also visually testify that the phase separation dynamics is respected by the postprocessing step. This is particularly clear in the passages where the bound preservation is deactivated. For example around $t\approx 1$ in \autofigref{fig:ch_pol_studies}, when the bound preservation is not needed, the maximum value of the bounded solvers naturally returns to align with the unbounded ones.


\section{Membrane force system}
\label{sec:geom_flows}

In the realm of mathematical modeling, principal curvatures have long been used for the theory of elastic plates and shells, with the work of Kirchhoff \cite{Kirchhoff1850} representing what may be the most famous case. 
Later in the '70s, the work of Canham \cite{Canham1970} and Helfrich \cite{Helfrich1973} on the link between biological membranes shape and curvature functionals has \emph{de facto} initiated an entirely new field of research focused on the use of geometry of hypersurfaces for biophysical modeling.  


We will here focus on the Canham-Helfrich functional $\cE_{B}(t)$ and its discretization. Considering membranes with boundaries, a commonly used expression for the energy functional is
\begin{equation}
\label{eq:ch_energy}
    \cE_{B}(t)=\int_{\Gamma}\frac{\gamma_W}{2} (\kappa-\kappa_0)^2 + \int_{\Gamma} \gamma_G\kappa_G  +\gamma_{\partial} \int_{\partial\Gamma} 1,
\end{equation}
where $\gamma_W,\; \gamma_G \in \bbR$ are bending rigidities, $\kappa:\cG_T\to\bbR$ is the previously introduced mean curvature, $\kappa_0 \in \bbR$ is the spontaneous curvature, and $\kappa_G: \cG_T\to\bbR$ is the Gaussian curvature $\kappa_G = \det \bbW$. The parameter $\gamma_\partial$ takes into account the possible line energy of $\partial\Gamma$. 
For the sake of simplicity we will consider the case $\gamma_G=constant$, $\gamma_W=constant$,  and assume that no topological change occurs along the evolution. Taken into account these simplifications, the gradient flow of a membrane $\Gamma$ is a family of evolving surfaces with normal velocity $\bv^\perp$ given by \cite{Nitsche1993,BarrettGarckeNurnberg2017a}
\begin{equation}
    \bv^\perp = \gamma_{W}\Bigl(-\Delta_\Gamma\kappa+\frac{1}{2}\kappa(\kappa-\kappa_0)^2-(\kappa-\kappa_0) |\bbW|^2\Bigr)\bn_\Gamma.
    \label{eq:willmore_velocity}
\end{equation}
Note that such velocity is in the normal direction since tangential changes to the surface shape do not modify its energy. Equation \eqref{eq:willmore_velocity} can be considered as a generalized Willmore flow equation. The mean curvature $\kappa$ is directly coupled to the surface $\Gamma$ through a second-order dependency. This can easily be seen in the case of surfaces whose configuration is described by a levelset equation. In that case, $\kappa$ is the trace of the extended Weingarten map $\bbH$ which in turn can be computed as the Hessian of the levelset for $\Gamma$. Equation \eqref{eq:willmore_velocity} thus leads to a fourth-order nonlinear equation posed on a moving manifold.   
\begin{rmk}
    Other geometrical flows like the mean curvature flow and the surface diffusion flow are also of interest in this field. Being somewhat simpler, extensive research is already available for them, while less is known about the Willmore flow, which we focus on in this article. We refer to \cite{DziukElliott2007, DeckelnickDziukElliott2005, BarrettGarckeNurnberg_bookSection} for further details.
\end{rmk}
\begin{rmk}
    Since this is a defining equation for the normal \emph{velocity}, additional forces can be taken into account by introducing a right-hand side to \eqref{eq:willmore_velocity}.
\end{rmk}

\subsection{Numerical methods for Helfrich system}
We choose here the family of BGN-type algorithms as pioneered in \cite{BarrettGarckeNurnberg2007} and later adapted for Willmore flow in \cite{BarrettGarckeNurnberg2016,BarrettGarckeNurnberg2017a}. The schemes pivot around the two following ingredients:
\begin{itemize}
    \item Introduce the \emph{mean curvature vector} $\bkappa(t) = \kappa(t)\bn_\Gamma(t)$ as independent variable and use the defining equation for the mean curvature vector $\bkappa =\Delta_{\Gamma}\mathbf{Id}_{\Gamma}$ in order to reduce the fourth order problem to a coupled second order one \cite{BarrettGarckeNurnberg2007};
    \item Use the weak formulation to rewrite and simplify the term $\kappa |\bbW|^2$ (more precisely $\kappa |\bbH|^2$ in the immersed setting) taking advantage of the newly introduced variable $\bkappa$. This reformulation not only simplifies and linearizes the system but also leads to stable schemes as proven in \cite{Dziuk2008}.
\end{itemize}
Our scheme choice builds on the one derived in \cite{BarrettGarckeNurnberg2017a} for boundary value problems.  Among the different types of boundary conditions considered in \cite{BarrettGarckeNurnberg2017a} we will restrict to the case of \emph{clamped boundary conditions}, for which
\begin{equation}
    \partial\Gamma(t)=\partial\Gamma(0)=\partial \widehat{\Gamma} \text{ and } \bn_{\partial\Gamma}(t)=\bmu(t) \text{ on } \partial\Gamma,
    \label{eq:bgn_bc}
\end{equation}
where $\bmu(t)$ is a user-given function that dictates the evolution of $\bn_{\partial\Gamma}(t)$.
The weak form of \eqref{eq:willmore_velocity} then reads: Find $\bv^\perp, ~\by, ~\bkappa\in [H(t)]^d$ such that
\begin{subequations}
\begin{align}
    &\inner{\bv^\perp}{\bphi} =~ \inner{\nabla_\Gamma \by}{\nabla_\Gamma \bphi} + \inner{\nabla_\Gamma\cdot \by}{\nabla_\Gamma\cdot \bphi}-\inner{(\nabla_\Gamma\by)^T}{\cD(\bphi)\;\bbP_\Gamma} \\
    &\qquad- \kappa_0\inner{\bkappa}{(\nabla_\Gamma\bphi)^T \bn_\Gamma}-\frac{1}{2}\inner{\gamma_W|\bkappa - \kappa_0\bn_\Gamma|^2\;\bbP}{\nabla_\Gamma \bphi} \nonumber\\ 
    &\qquad+ \inner{(\by\cdot\bkappa)\;\bbP_\Gamma}{\nabla_\Gamma \bphi},  \nonumber \\
    &\inner{\by}{\bpsi} =~ \inner{\gamma_W(\bkappa - \kappa_0\bn_\Gamma)}{\bpsi}, \label{eq:willmore_weak_y}\\
    &\inner{\bkappa}{\bxi} =~ -\inner{\bbP_\Gamma}{\nabla_\Gamma \bxi}+ \inner{\bmu}{\bxi}_{\partial\Gamma}, \label{eq:willmore_weak_kappa}
\end{align}
\label{eq:willmore_weak}
\end{subequations}
for all $\bphi, \bpsi, \bxi \in [H(t)]^d$, where $H(t)$ is the Sobolev space  $W^{1,2}(\Gamma(t))$. In the above $\cD(\bphi) = \bbP_{\Gamma}(\nabla_\Gamma\bphi + \nabla_\Gamma\bphi^T)\bbP_\Gamma$ and $\by$ is an auxiliary variable introduced to deal with the spontaneous curvature. 
For the discretization of the normal velocity in \eqref{eq:willmore_weak} we choose to introduce the displacement variable $\bd_h^{n+1}\in[V_h(\Gamma_h^n)]^d$ defined as
\begin{equation}
    \bv_h^\perp(t_n) = \frac{1}{\tau}\sum_i(\by_i(t_{n+1}) - \by_i(t_n))\phi_i(t_n) = \frac{1}{\tau}\bd_h^{n+1},
    \label{eq:discrete_displacement}
\end{equation}
where we recall that $\by_i$ are the mesh points and $\{\phi_i\}$ the associated basis such that $\phi_i(\by_j) = \delta_{ij}$ at all times.
Moreover, \eqref{eq:willmore_weak_kappa} and \eqref{eq:willmore_weak_y} can be collapsed to
\begin{align}
    \frac{1}{\gamma_W}\inner{\by^{n+1}}{\bpsi}_{\Gamma^n} + \inner{\nabla_{\Gamma^n}\bd^{n+1}}{\nabla_{\Gamma^n} \bpsi}_{\Gamma^n} = -\inner{\bbP_{\Gamma^n}}{\nabla_{\Gamma^n} \bpsi}_{\Gamma^n} \nonumber\\ 
    + \inner{\kappa_0\bn_{\Gamma^n}}{\bpsi}_{\Gamma^n}+ \inner{\bmu}{\bpsi}_{\partial\Gamma^n}.
\end{align}
The fully discrete form of \eqref{eq:willmore_weak} reads: Given $\by^{n}_h, \bkappa^{n}_h\in [V_h(\Gamma_h^n)]^d$ find \break
$\bd^{n+1}_h\in [V_{h0}(\Gamma_h^n)]^d$ and $\by^{n+1}_h\in [V_h(\Gamma_h^n)]^d$ such that
\begin{subequations}
\begin{align}
    &\inner[auto]{\frac{\bd_h^{n+1}}{\tau}}{\bphi_h}_{\Gamma_h^n}^h -\inner{\nabla_{\Gamma_h^n} \by^{n+1}_h}{\nabla_{\Gamma_h^n} \bphi_h}_{\Gamma_h^n}=~ \inner{\nabla_{\Gamma_h^n}\cdot \by^n_h}{\nabla_{\Gamma_h^n}\cdot \bphi_h}_{\Gamma_h^n}\\
    &\qquad-\inner{(\nabla_{\Gamma_h^n}\by^n_h)^T}{\cD(\bphi_h)\;\bbP_{\Gamma_h^n}}_{\Gamma_h^n}- \kappa_0\inner{\bkappa^n_h}{(\nabla_{\Gamma_h^n}\bphi_h)^T \bn_{\Gamma_h^n}}_{\Gamma_h^n}^h \nonumber \\
    &\qquad-\frac{1}{2}\inner{\gamma_W|\bkappa^n_h - \kappa_0\bn_{\Gamma_h^n}|^2\;\bbP_{\Gamma_h^n}}{\nabla_{\Gamma_h^n} \bphi_h}_{\Gamma_h^n}^h + \inner{(\by^n_h\cdot\bkappa^n_h)\;\bbP_{\Gamma_h^n}}{\nabla_{\Gamma_h^n} \bphi_h}_{\Gamma_h^n}^h,  \nonumber \\
    &\frac{1}{\gamma_W}\inner{\by^{n+1}_h}{\bpsi_h}_{\Gamma_h^n}^h+ \inner{\nabla_{\Gamma_h^n}\bd_h^{n+1}}{\nabla_{\Gamma_h^n} \bpsi_h}_{\Gamma_h^n} \\
    & \qquad= -\inner{\bbP_{\Gamma_h^n}}{\nabla_{\Gamma_h^n} \bpsi_h}_{\Gamma_h^n}+ \inner{\kappa_0\bn_{\Gamma_h^n}}{\bpsi_h}_{\Gamma_h^n}^h + \inner{\bmu}{\bpsi_h}_{\partial\Gamma^n}, \nonumber\label{eq:auxiliary_mean_curvature}
\end{align}
\label{eq:discrete_willmore}
\end{subequations}
for all $\bphi_h\in [V_{h0}(\Gamma_h^n)]^d,~\bpsi_h\in [V_h(\Gamma_h^n)]^d$, where $\inner{\cdot}{\cdot}^h$ is the mass-lumped inner product.
\begin{rmk}
    The mean curvature $\bkappa_h^{n+1}$ can be recovered a posteriori as $\bkappa_h^{n+1} = \pi^h(\gamma_W^{-1}\by_h^{n+1}+\kappa_0\bn_{\Gamma_h^n})$ where $\pi^h(\cdot)$ is the standard interpolation operator.
\end{rmk}
\begin{rmk}
    Note that \eqref{eq:discrete_willmore} corresponds to the scheme presented in \cite{BarrettGarckeNurnberg2016} for the choice of the parameter $\theta=1$, i.e. without tangential motion control.
\end{rmk}

\subsection{Numerical results for Helfrich system}
\label{subsec:willmore_conv}

To test the accuracy of our method for gradient flows, we start by performing classical studies on Willmore flow such as the ones in \cite{Dziuk2008,BarrettGarckeNurnberg2016,BarrettGarckeNurnberg2017a}. For every test we compare
\begin{enumerate}
\label{item:willmore_solvers}
    \item \label{item:willmore_naive} \textbf{hfSolver \ref{item:willmore_naive}}: algorithm \eqref{eq:discrete_willmore}.
    \item \label{item:willmore_duanli} \textbf{hfSolver \ref{item:willmore_duanli}}: algorithm \eqref{eq:discrete_willmore} together with the mesh redistribution described in \eqref{eq:duanli_discrete}. We highlight that it is of crucial importance that, when the algorithm in Section \ref{sec:mesh_redistribution} is employed, geometrical quantities characteristic of the mesh are updated. In our case, we refer to the mean curvature $\bkappa^n_h$ and the auxiliary variable $\by^n_h$ of the algorithm in \eqref{eq:discrete_willmore}.
\end{enumerate}
We start by testing the convergence  using the example presented in \cite[Section 5.3]{BarrettGarckeNurnberg2016}. It is an evolving sphere with initial radius $R_0$ and spontaneous curvature $\kappa_0$. The radius evolution $R(t)$ is the solution to the ordinary differential equation
\begin{equation}
    \partial_tR = -\frac{\kappa_0}{R}\Bigl(\frac{2}{R}+\kappa_0\Bigr), \quad R(0)=R_0>0.
\end{equation}
The error norm chosen to verify convergence is the following
\begin{equation}
    \norm{\Gamma - \Gamma_h}_{L^\infty(L^\infty)} = \max_{n}\max_{i}|\norm{\by_i(t_n)} - R(t_n)|,
    \label{eq:nodal_linfty_norm}
\end{equation}
where $\by_i$ are the vertices of the mesh simplices.
\begin{figure}
    \centering
    \begin{subfigure}[tbhp!]{0.9\textwidth}
        \centering
        \includegraphics[clip, trim = {0cm 18.2cm 0cm 3.8cm}, width = \textwidth]{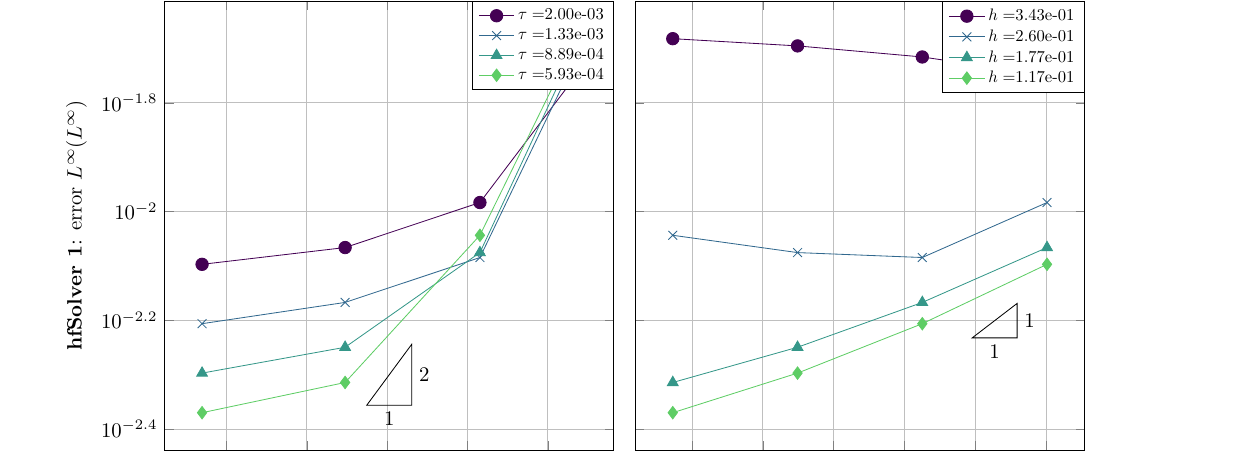}
    \end{subfigure}%

    \begin{subfigure}[tbhp!]{0.9\textwidth}
        \centering
        \includegraphics[clip, trim = {0cm 6.5cm 0cm 14.4cm}, width = \textwidth]{Fig7-test_sphere_spontaneous_mc.pdf}
    \end{subfigure}
    \caption{Convergence studies for Willmore flow of sphere under spontaneous mean curvature  for the solvers in List \ref{item:willmore_solvers}. As expected, first and second order
    convergence are achieved in time and space, respectively.}
    \label{fig:willmore_spontaneous_curv_conv}
    \begin{subfigure}[tbhp!]{0.9\textwidth}
        \centering
        \includegraphics[clip, trim = {0cm, 18cm, 10cm, 2.8cm}, width = 0.6\textwidth]{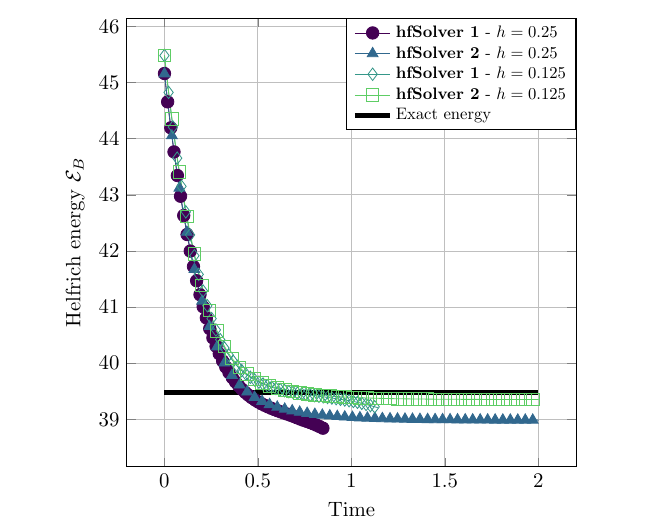}
    \end{subfigure}
    \caption{Energy evolution for torus with $R=2$, $r=1$ for the solvers in List \ref{item:willmore_solvers}.}
    \label{fig:willmore_torus21ngsolve}
\end{figure}
Results are presented in \autofigref{fig:willmore_spontaneous_curv_conv}. The convergence in space and time is clearly visible and in accordance with the optimal expected. Moreover, we can see that the use of the mesh redistribution actually improves the convergence rates. We believe this behavior is due to the automatic mesh adaptivity of \eqref{eq:duanli_discrete} (already mentioned in the original article \cite{DuanLi2024}), which tends to accumulate nodes in certain regions providing better geometry description.

As a second test we check the energy evolution for a torus of major radius $R=2$ and minor radius $r=1$ as done in \cite[Section 5.1]{BarrettGarckeNurnberg2016}. The minimizer of the Willmore energy in this case is the Clifford torus with radius $R=\sqrt{2}, \; r=1$, which has energy $\cE = 4\pi^2\approx 39.4784$. We use both solvers and plot the energy history for two different examples described by the following parameter sets:
\vspace{2mm}
\begin{center}
\begin{tabular}{||c|c|c|c|c||} 
 \hline
 $h$ & $nv$ & $ne$ & $\tau$ & $T$\\ [0.5ex] 
 \hline\hline
 0.2 & 2238 & 4476 & $10^{-3}$ &  2 \\ 
 \hline
\end{tabular}
~
\begin{tabular}{||c|c|c|c|c||} 
 \hline
 $h$ & $nv$ & $ne$ & $\tau$ & $T$\\ [0.5ex] 
 \hline\hline
 0.1 & 8951 & 17902 & $10^{-3}$ & 2 \\ 
 \hline
\end{tabular}
\end{center}
\vspace{2mm}
The results are shown in \autofigref{fig:willmore_torus21ngsolve}. We notice how both simulations using \textbf{hfSolver \ref{item:willmore_naive}} break before reaching the end time. This is due to instabilities that occur at regions of high curvature already observed in \cite[Figure 2]{BarrettGarckeNurnberg2016} for $\theta=1$. Applying the redistribution doesn't impact the energy landscape and allows to proceed with the simulation until the stable state is reached. Moreover, refining the mesh shows how the energy converges towards the expected value.

\subsubsection{Pinching tests for Helfrich system}
\label{subsec:helfrich_pinching}

\begin{figure}[tbhp!]
    \centering
    \begin{subfigure}[tbhp!]{0.32\textwidth}
        \centering
        \includegraphics[width = \textwidth, clip, trim = {0cm 0cm 0cm 0cm}]{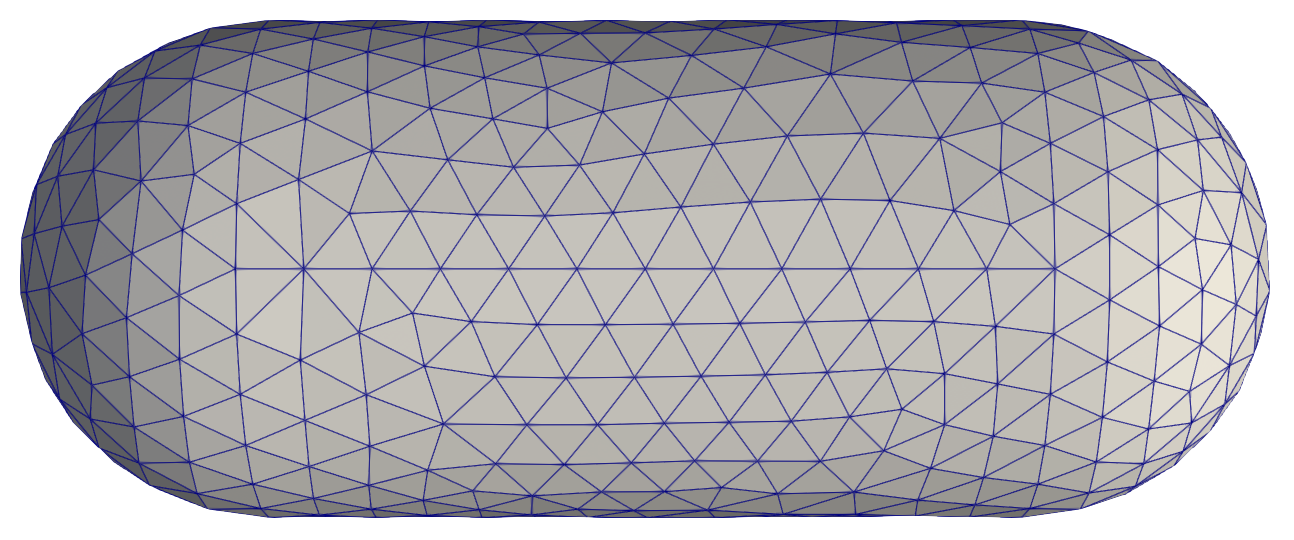}
    \end{subfigure}%
    ~
    \begin{subfigure}[tbhp!]{0.32\textwidth}
        \centering
        \includegraphics[width = \textwidth, clip, trim = {0cm 0cm 0cm 0cm}]{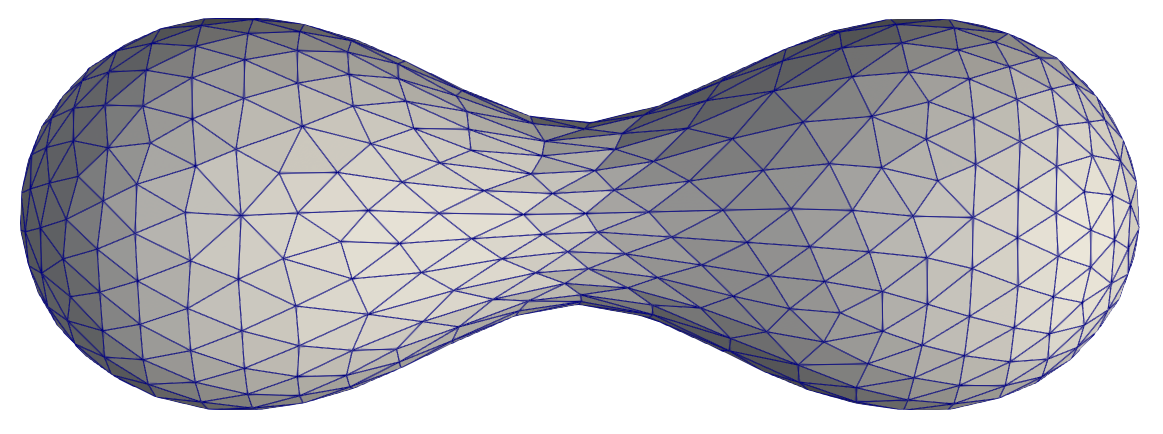}
    \end{subfigure}
    ~
    \begin{subfigure}[tbhp!]{0.32\textwidth}
        \centering
        \includegraphics[width = \textwidth, clip, trim = {0cm 0cm 0cm 0cm}]{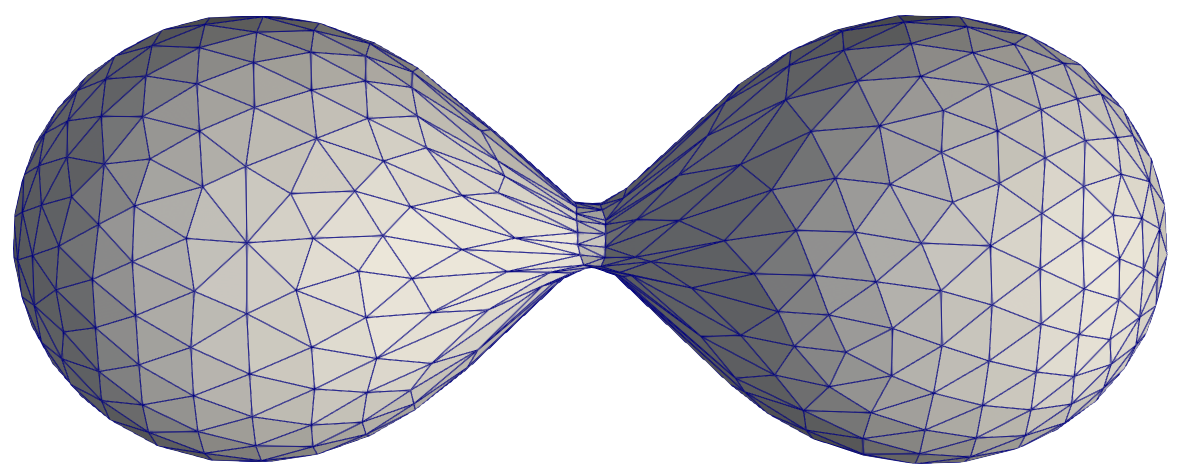}
    \end{subfigure}
    
    \begin{subfigure}[tbhp!]{0.32\textwidth}
        \centering
        \includegraphics[width = \textwidth, clip, trim = {0cm 0cm 0cm 0cm}]{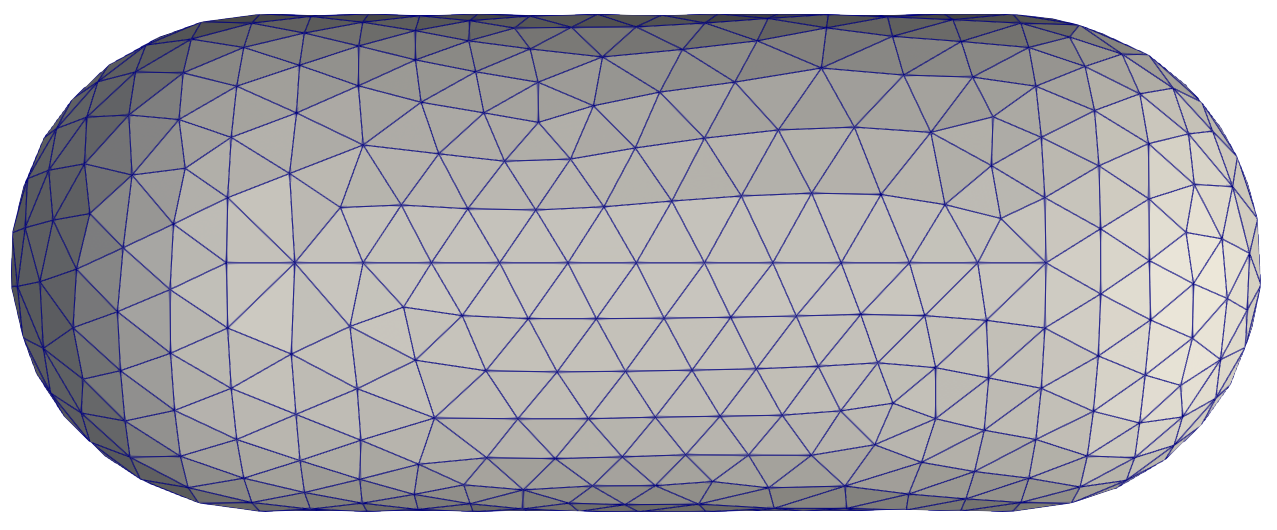}
        \caption{t=0}
        \label{fig:cigar31_0}
    \end{subfigure}%
    ~
    \begin{subfigure}[tbhp!]{0.32\textwidth}
        \centering
        \includegraphics[width = \textwidth, clip, trim = {0cm 0cm 0cm 0cm}]{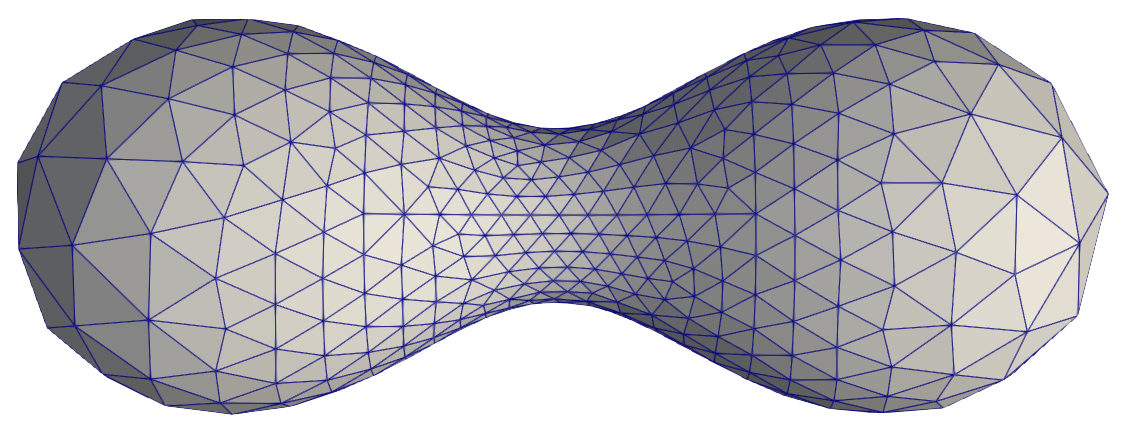}
        \caption{t=0.5}
    \end{subfigure}
    ~
    \begin{subfigure}[tbhp!]{0.32\textwidth}
        \centering
        \includegraphics[width = \textwidth, clip, trim = {0cm 0cm 0cm 0cm}]{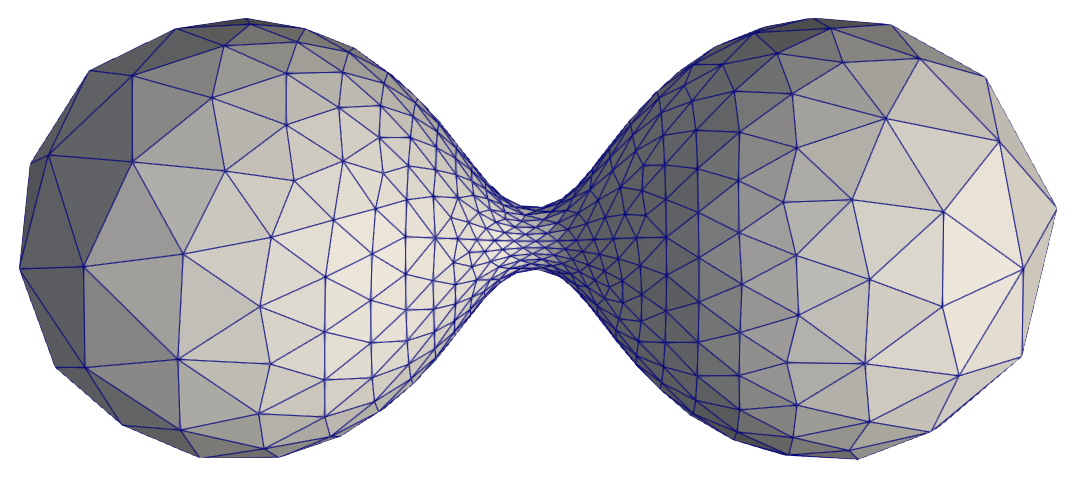}
        \caption{t=1}
        \label{fig:cigar511_duanli_2}
    \end{subfigure}

    \centering
    \begin{subfigure}[tbhp!]{0.32\textwidth}
        \centering
        \includegraphics[width = \textwidth, clip, trim = {0cm 0cm 0cm 0cm}]{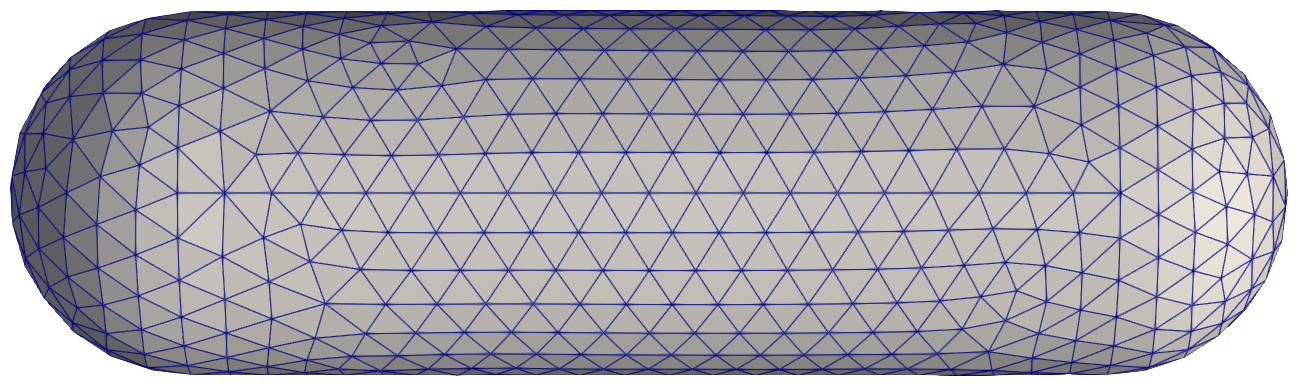}
    \end{subfigure}%
    ~
    \begin{subfigure}[tbhp!]{0.32\textwidth}
        \centering
        \includegraphics[width = \textwidth, clip, trim = {0cm 0cm 0cm 0cm}]{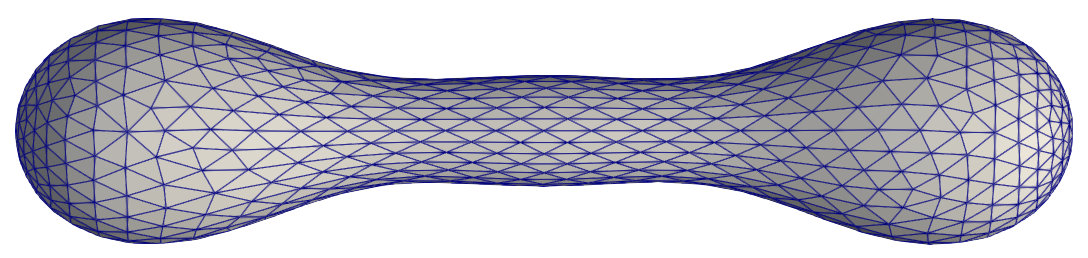}
    \end{subfigure}
    ~
    \begin{subfigure}[tbhp!]{0.32\textwidth}
        \centering
        \includegraphics[width = \textwidth, clip, trim = {0cm 0cm 0cm 0cm}]{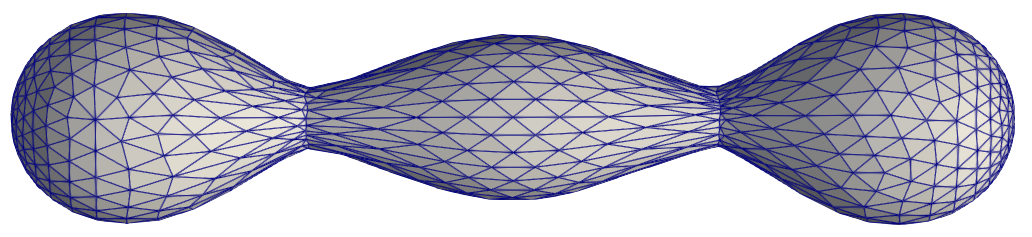}
    \end{subfigure}
    
    \begin{subfigure}[tbhp!]{0.32\textwidth}
        \centering
        \includegraphics[width = \textwidth, clip, trim = {0cm 0cm 0cm 0cm}]{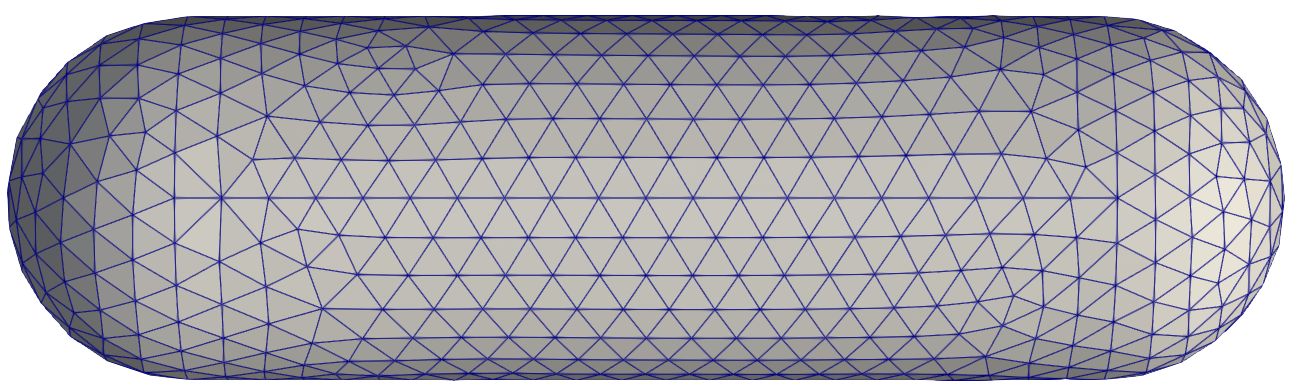}
        \caption{t=0}
        \label{fig:cigar51_0}
    \end{subfigure}%
    ~
    \begin{subfigure}[tbhp!]{0.32\textwidth}
        \centering
        \includegraphics[width = \textwidth, clip, trim = {0cm 0cm 0cm 0cm}]{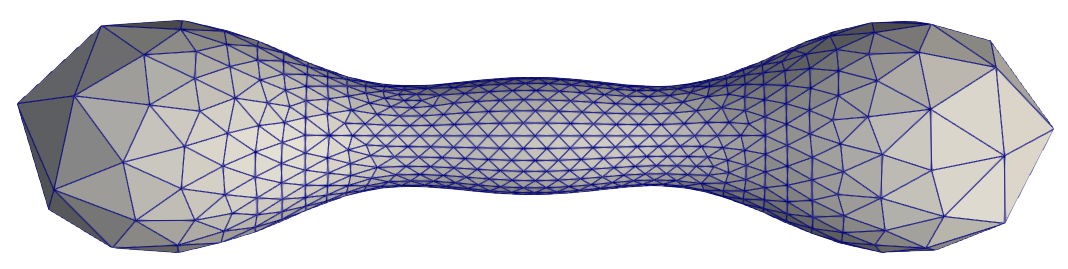}
        \caption{t=0.15}
    \end{subfigure}
    ~
    \begin{subfigure}[tbhp!]{0.32\textwidth}
        \centering
        \includegraphics[width = \textwidth, clip, trim = {0cm 0cm 0cm 0cm}]{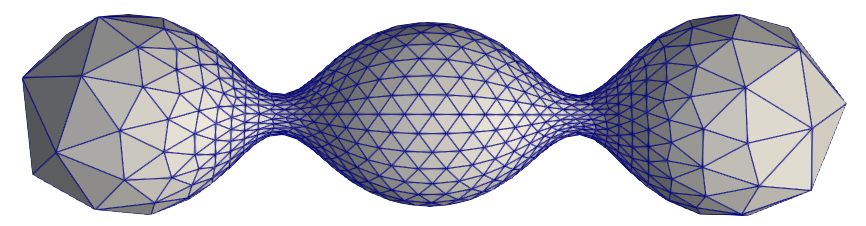}
        \caption{t=0.3}
    \end{subfigure}

    \begin{subfigure}[tbhp!]{0.48\textwidth}
        \centering
        \includegraphics[width = \textwidth, clip, trim = {1cm 18cm 10cm 3cm}]{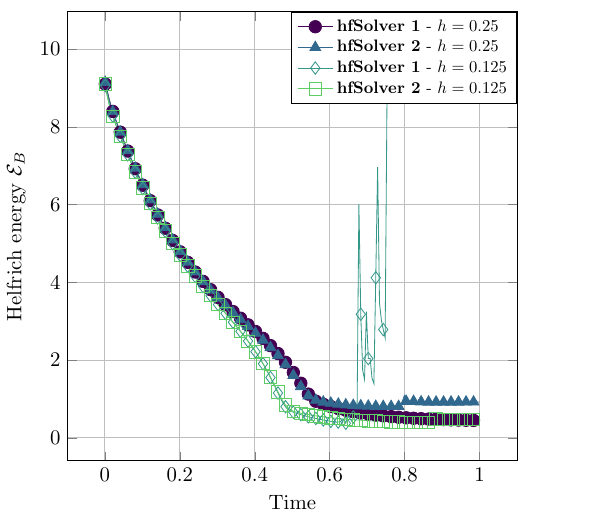}
        \caption{Willmore energy cigar \textbf{hf-Test 1}}
    \end{subfigure}
    ~
    \begin{subfigure}[tbhp!]{0.48\textwidth}
        \centering
        \includegraphics[width = \textwidth, clip, trim = {1cm 18cm 10cm 3cm}]{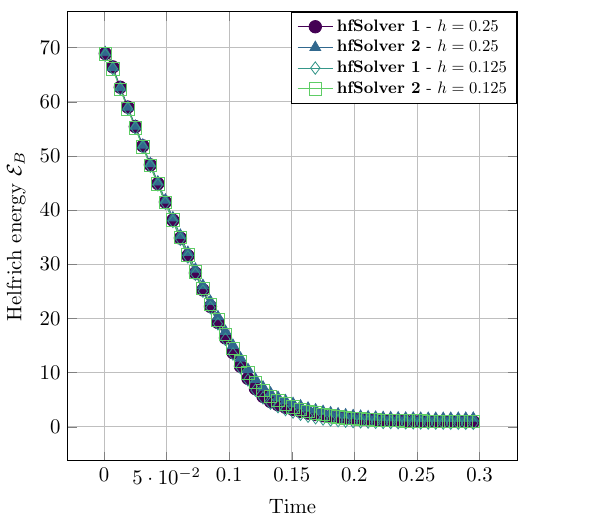}
        \caption{Willmore energy cigar \textbf{hf-Test 2}}
    \end{subfigure}

    \caption{(a)-(c): Mesh evolution for \textbf{hf-Test 1} (coarser mesh). Top: \textbf{hfSolver 1}, bottom \textbf{hfSolver 2}. (d)-(f): Mesh evolution for \textbf{hf-Test 2} (coraser mesh). Top: \textbf{hfSolver 1}, bottom \textbf{hfSolver 2}. (g) Energy evolution for \textbf{hf-Test 1} and (h) Energy evolution for \textbf{hf-Test 2} for the different solvers in List \ref{item:willmore_solvers}. We note that figures are rescaled.}
    \label{fig:cigar}
\end{figure}

We perform consolidated tests for the Willmore flow under the influence of spontaneous mean curvature as presented in \cite[Section 5.3]{BarrettGarckeNurnberg2016} and then repeated in \cite{GarckeNurnbergZhao2025}, where a similar two-step procedure as the one adopted here is used. We consider two different tests for cigar-shaped surfaces:
\begin{enumerate}
    \item \textbf{hf-Test 1} with initial mesh as in \autofigref{fig:cigar31_0} and the following sets of parameters
    \vspace{2mm}
    \begin{center}
    \begin{tabular}{||c|c|c|c|c||} 
    \hline
    $h$ & $nv$ & $ne$ & $T$ & $\kappa_0$\\ [0.5ex] 
    \hline\hline
    0.25 & 601 & 1198 & 1 & -2 \\ 
    \hline
    \end{tabular}
    ~
    \begin{tabular}{||c|c|c|c|c||} 
    \hline
    $h$ & $nv$ & $ne$ & $T$ & $\kappa_0$\\ [0.5ex] 
    \hline\hline
    0.125 & 2321 & 4638 & 1 & -2 \\ 
    \hline
    \end{tabular}
    \end{center}
    \vspace{2mm}
    \item \textbf{hf-Test 2} with initial mesh as in \autofigref{fig:cigar51_0} and the following sets of parameters
    \vspace{2mm}
    \begin{center}
    \begin{tabular}{||c|c|c|c|c||} 
    \hline
    $h$ & $nv$ & $ne$ & $T$ & $\kappa_0$\\ [0.5ex] 
    \hline\hline
    0.25 & 845 & 1686 & 0.3 & -3 \\ 
    \hline
    \end{tabular}
    ~
    \begin{tabular}{||c|c|c|c|c||} 
    \hline
    $h$ & $nv$ & $ne$ & $T$ & $\kappa_0$\\ [0.5ex] 
    \hline\hline
    0.125 & 3282 & 6460 & 0.3 & -3\\ 
    \hline
    \end{tabular}
    \end{center}
    \vspace{2mm}
\end{enumerate}
A common timestep of $\tau=10^{-3}$ has been used for all the experiments. The idea is to test the ability of the algorithm to simulate pinching phenomena, since \textbf{hf-Test 1} develops one neck and \textbf{hf-Test 2} develops two necks. The results are plotted in \autofigref{fig:cigar} for both \textbf{hfSolver \ref{item:willmore_naive}} and \textbf{2}. It can be seen how the mesh adaptivity mentioned in \cite{DuanLi2024} leads to a better resolution of the neck and overall better results for \textbf{hfSolver \ref{item:willmore_duanli}}. By comparison of \autofigref{fig:cigar511_duanli_2} with \cite[Figure 12]{BarrettGarckeNurnberg2016} the shape is qualitatively in accordance with what is expected. Even more, while \cite[Figure 12]{BarrettGarckeNurnberg2016} requires a finer mesh to accurately resolve the three pearls, the mesh adaptivity of \cite{DuanLi2024} alleviates that requirement for similar resolution. \autofigref{fig:cigar} is also in accordance with \cite[Figure 5.14]{GarckeNurnbergZhao2025}, where a similar algorithm is employed for a cigar-shaped surface with different aspect ratio. As it happens in that article, it has to be noted that the refinement of the necks comes to the expense of a coarser mesh at the lobes. Looking at the energy evolution for \textbf{hf-Test 2}, the curve for \textbf{hfSolver \ref{item:willmore_naive}} and \textbf{hfSolver \ref{item:willmore_duanli}} fall in the same place for all experiments, highlighting that the tangential mesh redistribution does not impact the energy landscape of the shape dynamics. For \textbf{hf-Test 1} we have that the dynamics is stable for both examples when $h=0.2$, with a slight difference in energy when necking occurs. Refining the mesh leads to mesh breaking for \textbf{hfSolver 1}, while \textbf{hfSolver 2} maintains its stability, asymptotically tending towards the same energy \textbf{hfSolver 1} was tending to in the coarser case. Overall, the proposed technique appears stable in various regimes and accurate in resolving the energy evolution of the geometry.

\section{Coupling}
\label{sec:coupling}
Simulating actual biological behavior involves complex coupling between different equations. Given we make heavy use of post-processing techniques, we couple the various solvers in a staggered way. We have then the option to choose between explicit and implicit staggering. Both approaches have shown to be effective in the literature \cite{MerckerMarciniak-CzochraRichterEtAl2013, KhanwaleSaurabhIshiiEtAl2023, KhanwaleLofquistSundarEtAl2020,  GarckeNurnbergZhao2025}. In \cite[Appendix E]{BachiniKrauseNitschkeEtAl2023}, where fluid deformable membranes with phase separation are considered, both procedures show very similar results even in the most complex scenarios. Moreover, when implicit staggering is chosen, the authors of \cite{BachiniKrauseNitschkeEtAl2023} claim that convergence is achieved in four steps on average. Strong of these results we proceed in coupling our solvers.

\begin{rmk}
    An alternative option is to use a monolithic approach and couple the solvers in a unique system, avoiding staggering. In this way there is no need for sub-iterations and memory allocation is reduced. On the other side, the solution is usually only available for ad-hoc problems and not flexible. Together with this, eventual non-linearities have to be linearized. One can find an example of such a choice in \cite{MokbelMokbelLieseEtAl2024}, where wetting dynamics of liquid droplets on deformable membranes is simulated.
\end{rmk}

\subsection{Numerical results for coupled mean curvature and ADR equations}
\label{subsec:mcadr_numerical_tests}
\begin{figure}[tbhp!]
    \centering
    \begin{subfigure}[tbhp!]{\textwidth}
        \centering
        \includegraphics[clip, trim = {1cm 17.9cm 0.5cm 3.cm}, width = \textwidth]{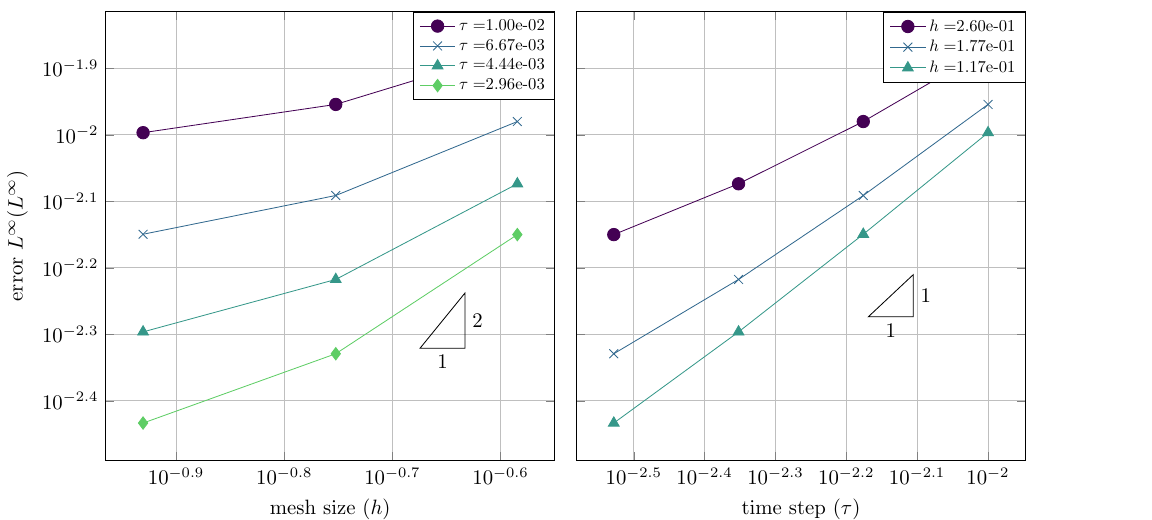}
    \end{subfigure}
    \caption{Convergence studies for \textbf{mcadrSolver 4} and problem \eqref{eq:mc_adr_conv_kovacs}. As expected, first and second order convergence are achieved in time and space, respectively.}
    \label{fig:mc_adr_conv}
\end{figure}

\begin{figure}[tbhp!]
    \begin{subfigure}[tbhp!]{0.48\textwidth}
        \centering
        \includegraphics[clip, trim = {0cm 0cm 0cm 0cm}, width = 0.8\textwidth]{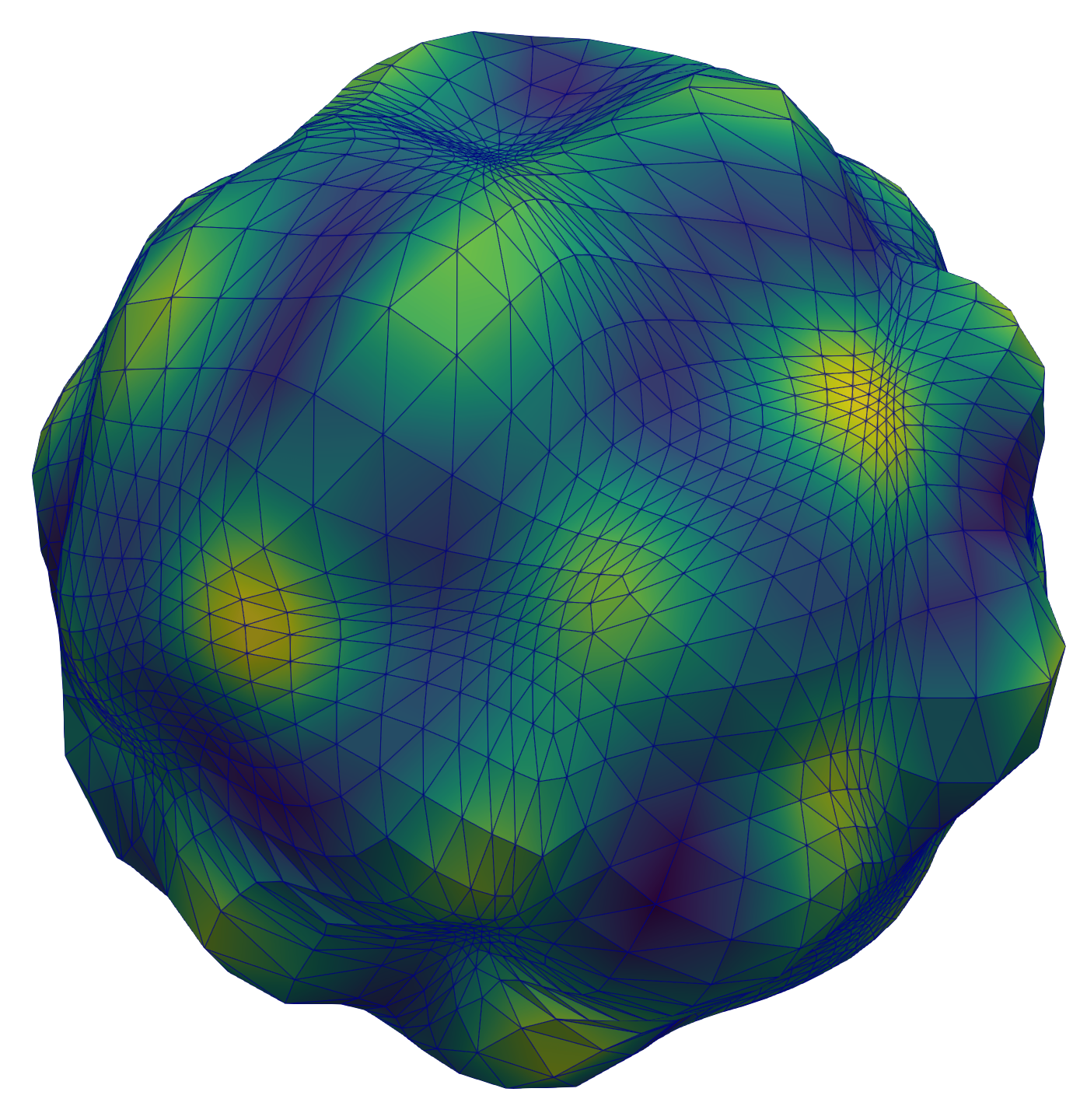}
        \caption{\textbf{mcadrSolver 1}}
    \end{subfigure}
    ~
    \begin{subfigure}[tbhp!]{0.48\textwidth}
        \centering
        \includegraphics[clip, trim = {0cm 0cm 0cm 0cm}, width = 0.8\textwidth]{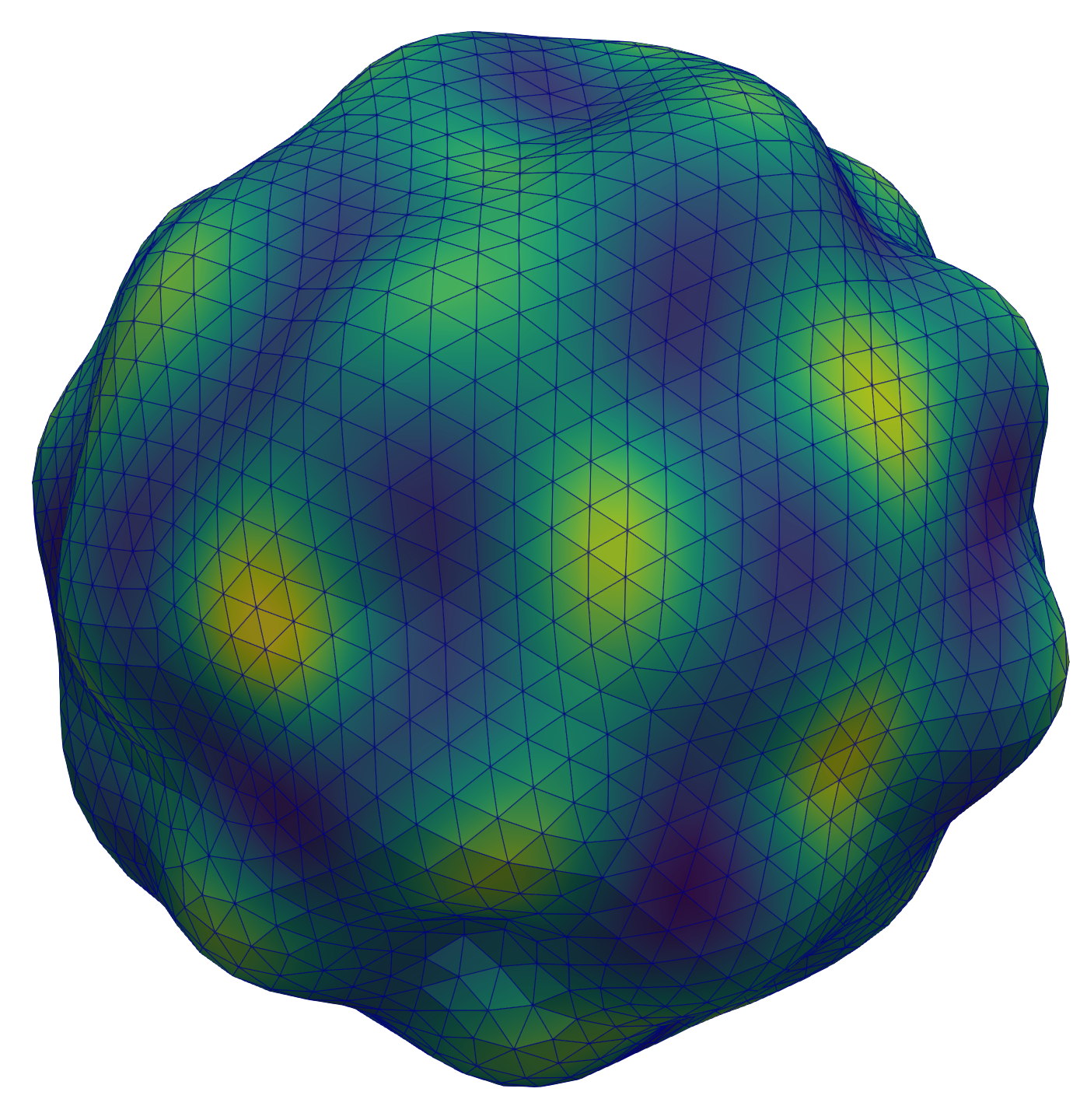}
        \caption{\textbf{mcadrSolver 2}}
    \end{subfigure}
    \caption{Simulation at $t=8$ for the model \eqref{eq:mc_adr_mesh_kovacs} \cite[Sec. 11.2]{KovacsLiLubichEtAl2017} and the implicitly staggered schemes in List \ref{item:mcadr_solvers}. }
    \label{fig:mc_adr_mesh}
\end{figure}

We begin by testing the convergence properties of our coupled solver. To do so we pick the coupling example presented in \cite{BarreiraElliottMadzvamuse2011, KovacsLiLubichEtAl2017}. The system which is solved is the following \cite[p.686]{KovacsLiLubichEtAl2017}:
\begin{subequations}
\begin{align}
    \matder u + u\nabla_\Gamma\cdot \bv -\Delta_\Gamma u &= f(t,\bx) \label{eq:mc_adr_conv_kovacs_1},\\
    \bv - \alpha \Delta_\Gamma\bv - \beta \Delta_\Gamma\Id_\Gamma &= (\delta u + g(t,\bx))\bn_\Gamma \label{eq:mc_adr_conv_kovacs_2},
\end{align}
\label{eq:mc_adr_conv_kovacs}
\end{subequations}
with parameters $\alpha, \beta, \delta \in \bbR^+$ and $\bx = (x_1, x_2, x_3)$. Convergence studies are performed choosing $f, g$ such that the exact solution for $u$ is $u(t, \bx) = x_1x_2e^{-6t}$. The geometry is chosen to be a sphere whose radius evolves following the law
\begin{equation}
    R(t) = \frac{r_0r_K}{r_Ke^{-kt}+r_0(1-e^{-kt})},
\end{equation}
with parameters $r_0,r_K, k \in \bbR^+$. The parameters are chosen as follows
\vspace{2mm}
\begin{center}
\begin{tabular}{||c|c|c|c|c|c|c||} 
\hline
$T$ & $\alpha$ & $\beta$ & $\delta$ & $r_0$ & $r_K$ & $k$\\ [0.5ex] 
\hline\hline
1 & 0 & 1 & 0.4 & 1 & 2 & 0.5 \\ 
\hline
\end{tabular}
\end{center}
\vspace{2mm}
To simulate \eqref{eq:mc_adr_conv_kovacs_1} we can use what has been introduced in Section \ref{sec:adr} and add a right-hand side coupling. Since we chose $\alpha=0$, we can simulate \eqref{eq:mc_adr_conv_kovacs_2} by simplifying \eqref{eq:discrete_willmore}. The result is a weighted mean curvature solver that, in this very case, takes the form
\begin{subequations}
\begin{align}
    \inner[auto]{\frac{\bd_h^{n+1}}{\tau}}{\bphi_h}_{\Gamma_h^n}^h -\beta\inner{ \bkappa^{n+1}_h}{\bphi_h}^h_{\Gamma_h^n}&=\inner{ (\delta u + g(t,\bx))\bn_{\Gamma_h^n}}{\bphi_h}^h_{\Gamma_h^n}, \\
    \inner{\bkappa^{n+1}_h}{\bpsi_h}_{\Gamma_h^n}^h+ \inner{\nabla_{\Gamma_h^n}\bd_h^{n+1}}{\nabla_{\Gamma_h^n} \bpsi_h}_{\Gamma_h^n} &= -\inner{\bbP_{\Gamma_h^n}}{\nabla_{\Gamma_h^n} \bpsi_h}_{\Gamma_h^n}.
\end{align}
\label{eq:mc_kovacs_discrete}
\end{subequations}
We perform the test using implicit staggering. The results of the convergence studies are shown in \autofigref{fig:mc_adr_conv}, where the same norm as in \eqref{eq:nodal_linfty_norm} is used to study the results. One can see how the expected convergence is achieved in both space and time. 
To qualitatively visualize the effect of the mesh redistribution and ADR stabilization, the example in \cite[Sec. 11.2]{KovacsLiLubichEtAl2017} is also reproduced. The simulation is derived from a proposed model for tumor growth. The coupled system is as follows:
\begin{subequations}
\begin{align}
    \matder u + u\nabla_\Gamma\cdot \bv -\Delta_\Gamma u &= f_1(u,w), \label{eq:mc_adr_mesh_kovacs_1}\\
    \matder w + w\nabla_\Gamma\cdot \bv -D_c\Delta_\Gamma w &= f_2(u,w), \label{eq:mc_adr_mesh_kovacs_2}\\
    \bv - \alpha \Delta_\Gamma\bv - \beta \nabla_\Gamma\Id_\Gamma &= \delta u\bn_\Gamma, \label{eq:eq:mc_adr_mesh_kovacs_3}
\end{align}
\label{eq:mc_adr_mesh_kovacs}
\end{subequations}
with the non-linear couplings
\begin{equation}
    f_1(u, w) = \gamma(a-u+u^2w), \quad f_2(u, w) = \gamma(b-u^2w).
    \label{eq:mc_adr_nonlinear_coupling}
\end{equation}
The parameters are set as
\vspace{2mm}
\begin{center}
\label{tab:mc_adr_parameters}
\begin{tabular}{||c|c|c|c|c|c|c|c|c|c||} 
\hline
$T$ & $\tau$ & h & $\alpha$ & $\beta$ & $\delta$ & a & b & $\gamma$ & $D_c$\\ [0.5ex] 
\hline\hline
8 & $10^{-3}$ & 0.07 & 0 & 0.01 & 0.4 & 0.1 & 0.9 & 100& 10 \\ 
\hline
\end{tabular}
\end{center}
\vspace{2mm}
As discussed  in \cite{KovacsLiLubichEtAl2017}, the nonlinear system composed by \eqref{eq:mc_adr_mesh_kovacs_1} and \eqref{eq:mc_adr_mesh_kovacs_2} is solved until $t=5$ without coupling it to the mesh evolution, i.e. ignoring 
\eqref{eq:eq:mc_adr_mesh_kovacs_3}. After that, the full system is evolved until $t=8$. We study the following solvers:
\begin{enumerate}
\label{item:mcadr_solvers}
    \item \textbf{mcadrSolver 1}: Implicit staggering of \eqref{eq:adr_weak_nostab} and \eqref{eq:mc_kovacs_discrete} and no mesh redistribution.
    \item \textbf{mcadrSolver 2}: Implicit staggering of \eqref{eq:adr_weak_cip} and \eqref{eq:mc_kovacs_discrete} and mesh redistribution.
\end{enumerate} 
In \autofigref{fig:mc_adr_mesh} we compare the mesh quality at the final time $t=8$. It is evident that the redistribution step keeps the mesh well-behaved while maintaining the solution accuracy. 

\begin{figure}[tbhp!]
    \begin{subfigure}[tbhp!]{0.48\textwidth}
        \centering
        \includegraphics[clip, trim = {0cm 5cm 0cm 5cm}, width = 0.8\textwidth]{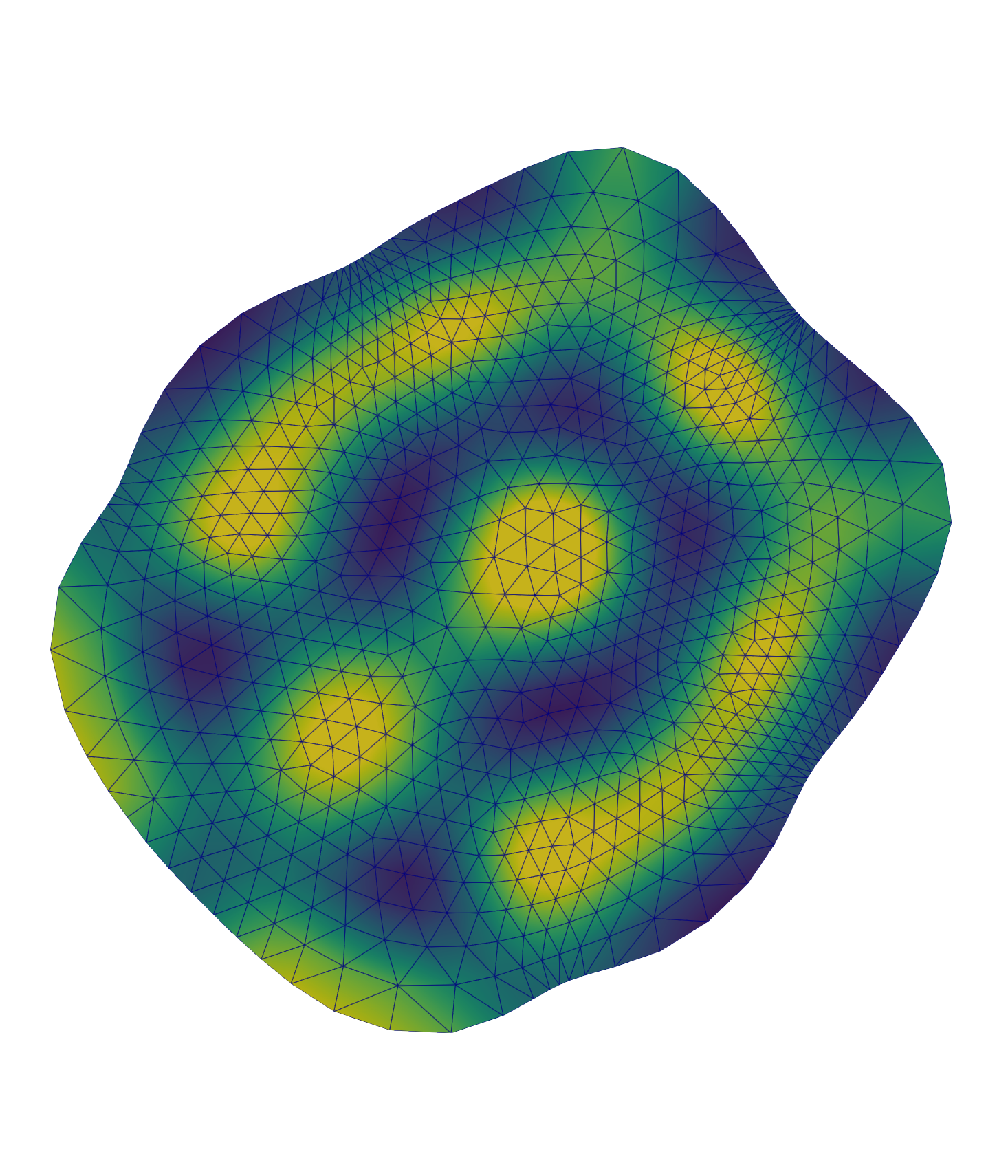}
        \caption{\textbf{mcadrSolver 1}}
    \end{subfigure}
    ~
    \begin{subfigure}[tbhp!]{0.48\textwidth}
        \centering
        \includegraphics[clip, trim = {0cm 5cm 0cm 5cm}, width = 0.8\textwidth]{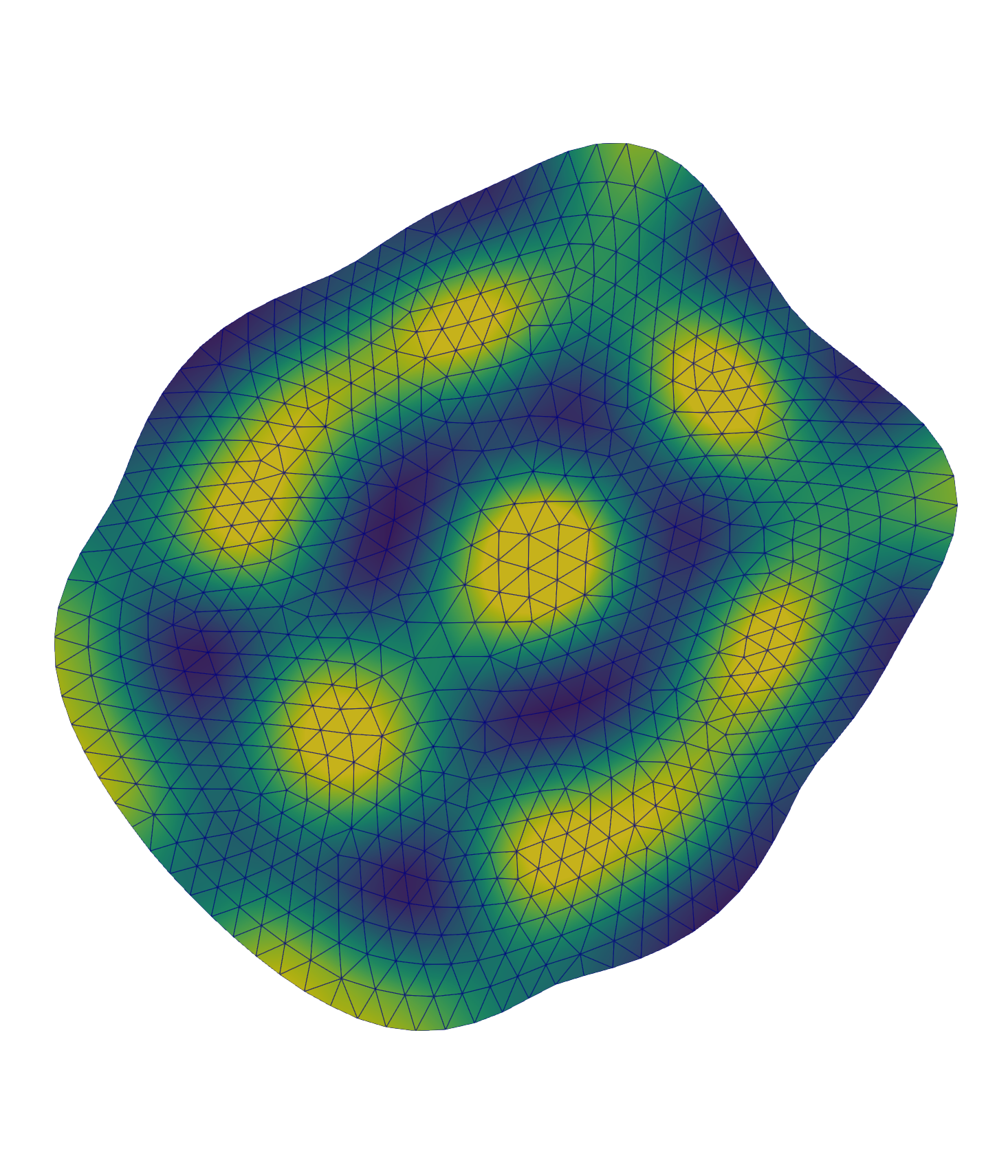}
        \caption{\textbf{mcadrSolver 2}}
    \end{subfigure}
    \caption{Simulation at $t=8$ for the model \eqref{eq:surfbulk_mc_adr_mesh_kovacs} and the implicitly staggered schemes in List \ref{item:mcadr_solvers}. }
    \label{fig:surfbulk_mc_adr_mesh}
\end{figure}

We can construct a bulk-surface problem by modifying the system in \eqref{eq:mc_adr_mesh_kovacs} with the following one
\begin{subequations}
\begin{align}
    \matder u + u\nabla\cdot \bv -\Delta u &= f_1(u,w) \quad\text{ in } \Omega, \\
    \matder w + w\nabla\cdot \bv -D_c\Delta w &= f_2(u,w) \quad\text{ in } \Omega,  \\
    \bv - \alpha \Delta_{\partial\Omega}\bv - \beta \nabla_{\partial\Omega}\Id_{\partial\Omega} &= \delta u\bn_{\partial\Omega}  \quad\text{ on } \partial\Omega . \label{eq:surfbulk_velocity}
\end{align}
\label{eq:surfbulk_mc_adr_mesh_kovacs}
\end{subequations}
The non-linear couplings and the parameters are kept the same as in \eqref{eq:mc_adr_nonlinear_coupling} and Table \ref{tab:mc_adr_parameters}. The initial geometry is a circle of unitary radius and the velocity prescribed by \eqref{eq:surfbulk_velocity} is extended in the bulk harmonically as described in Section \ref{sec:mesh_redistribution}. Results are presented in \autofigref{fig:surfbulk_mc_adr_mesh}, where it can be observed how the physics is correctly resolved while maintaining good mesh nodes distribution. 

\subsection{Numerical results for coupled Helfrich flow and Cahn-Hilliard euqations}
\label{subsec:wmch_numerical_tests}
\begin{figure}[tbhp!]
    \begin{subfigure}[tbhp!]{0.48\textwidth}
        \centering
        \includegraphics[clip, trim = {2cm 9cm 2cm 9cm}, width = 0.78\textwidth]{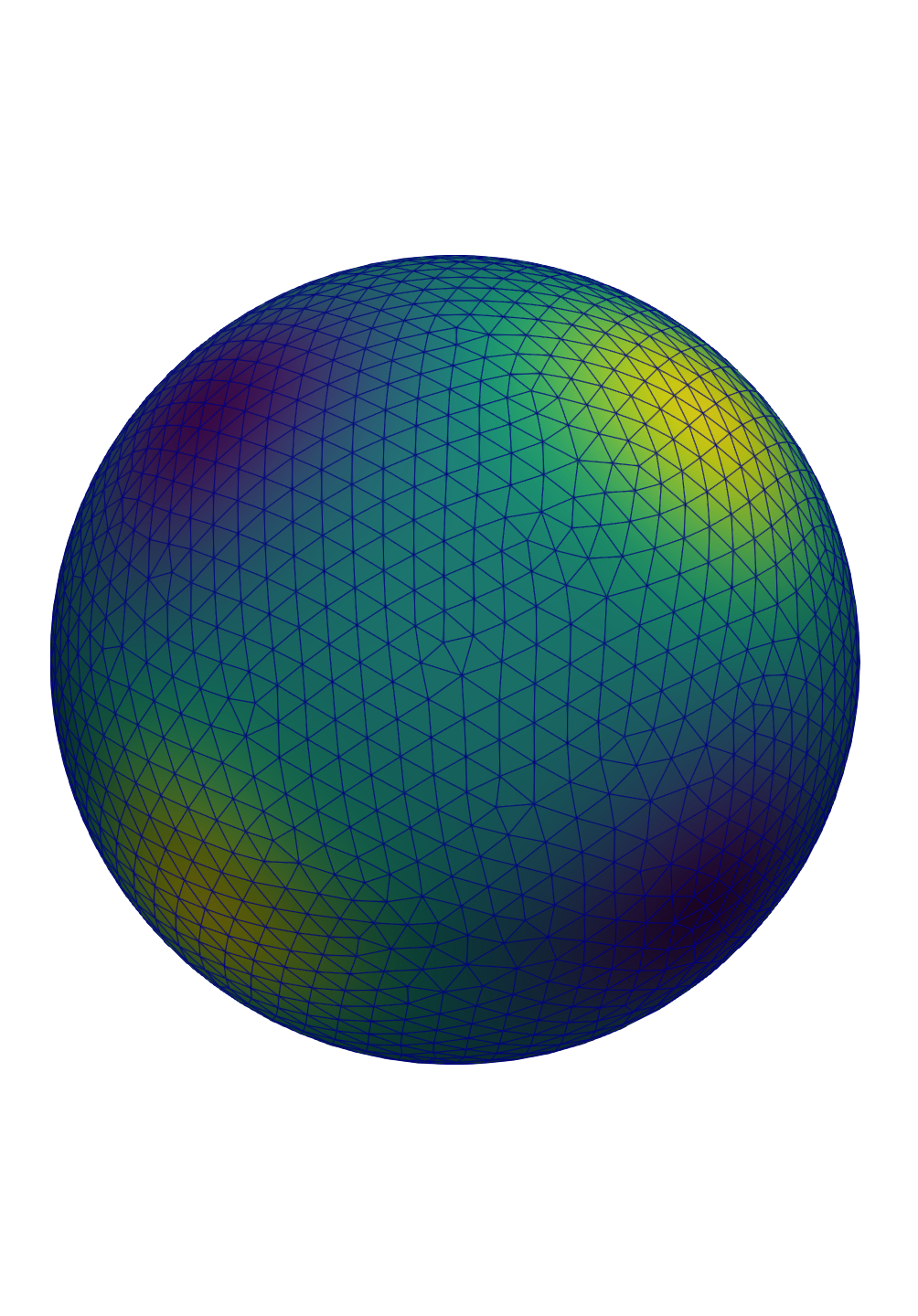}
        \caption{}
    \end{subfigure}%
    ~
    \begin{subfigure}[tbhp!]{0.48\textwidth}
        \centering
        \includegraphics[clip, trim = {2cm 9cm 2cm 9cm}, width = 0.78\textwidth]{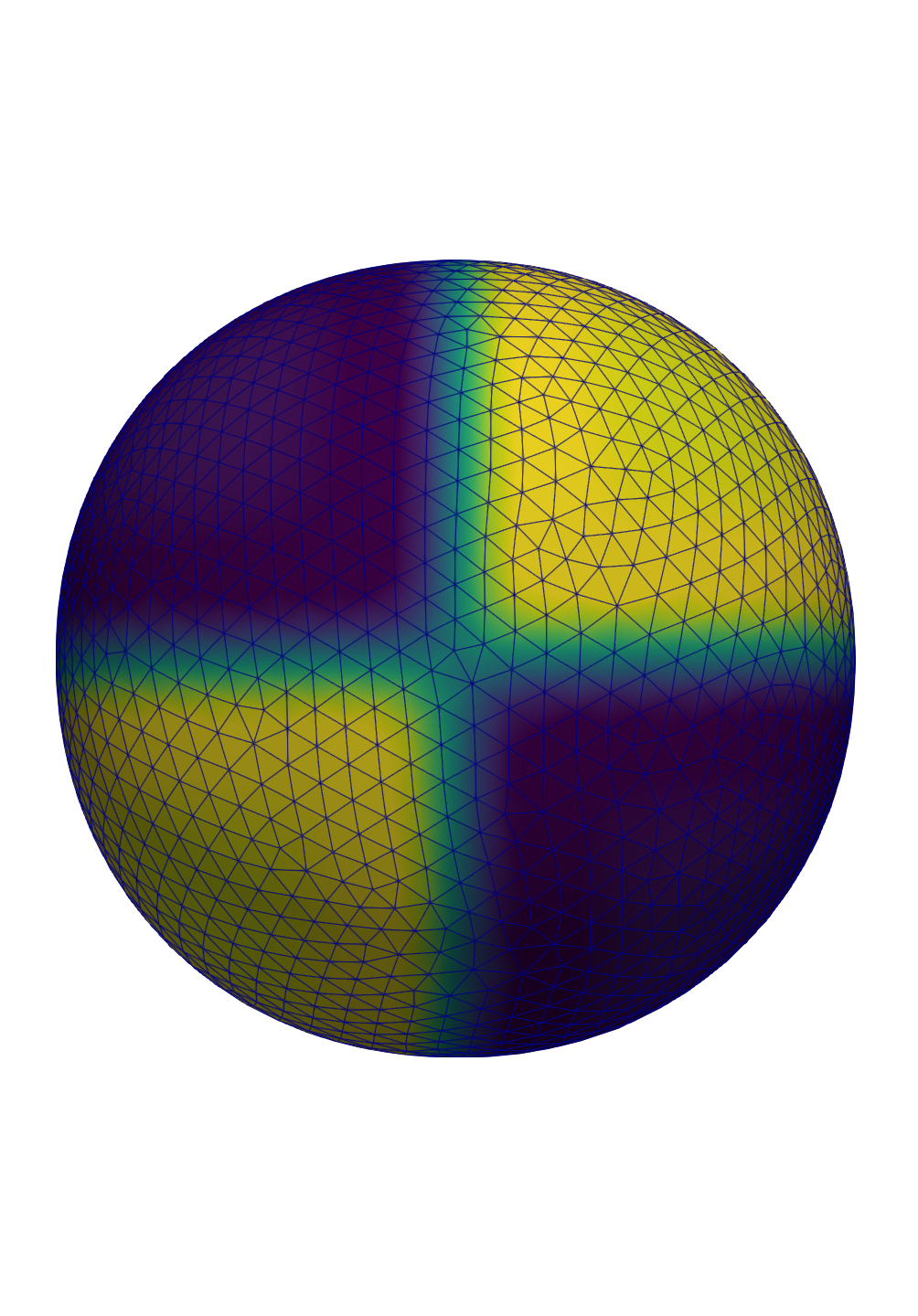}
        \caption{$\gamma_W=0.5, ~t=1$}
    \end{subfigure}

    \begin{subfigure}[tbhp!]{0.48\textwidth}
        \centering
        \includegraphics[clip, trim = {2cm 9cm 2cm 9cm}, width = 0.78\textwidth]{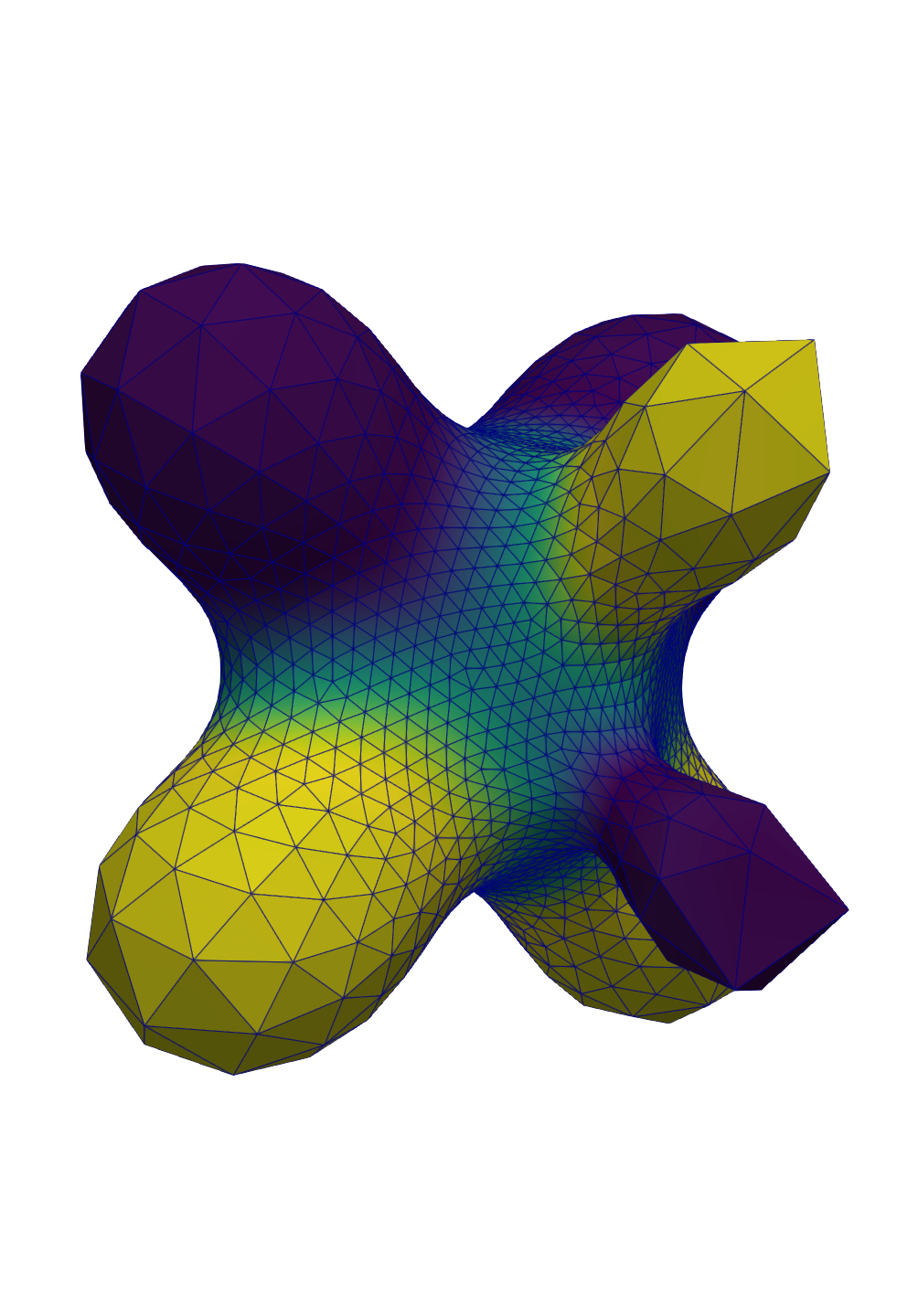}
        \caption{$\gamma_W=0.1, ~t=0.75$}
        \label{fig:bachini_trig_01}
    \end{subfigure}%
    ~
    \begin{subfigure}[tbhp!]{0.48\textwidth}
        \centering
        \includegraphics[clip, trim = {2cm 9cm 2cm 9cm}, width = 0.78\textwidth]{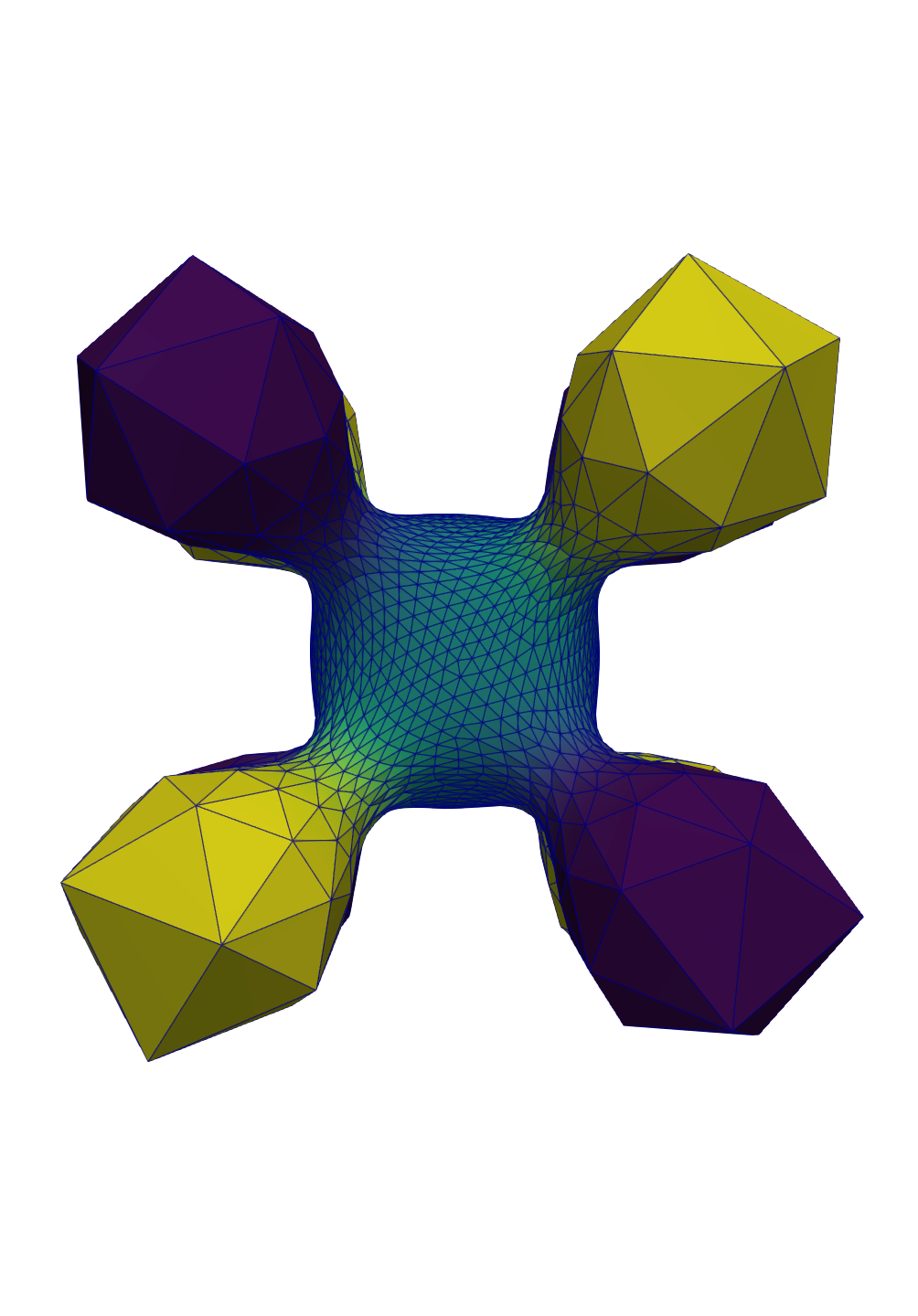}
        \caption{$\gamma_W=0.02, ~t=0.2$}
    \end{subfigure}
    \caption{(a) Initial setup for the problem \ref{eq:wmch_bachini_u0}. (b)-(d) Mesh for different elasticity modulus. The simulation is shown at end time $t=1$ or at last timestep before inextensibility threshold is violated.}
    \label{fig:bachini_trig}
\end{figure}

\begin{figure}[tbhp!]
    \begin{subfigure}[tbhp!]{0.48\textwidth}
        \centering
        \includegraphics[clip, trim = {2cm 7cm 2cm 7cm}, width = 0.78\textwidth]{Fig26-bachini_trig_initial.png}
        \caption{}
    \end{subfigure}%
    ~
    \begin{subfigure}[tbhp!]{0.48\textwidth}
        \centering
        \includegraphics[clip, trim = {2cm 6cm 2cm 6cm}, width = 0.78\textwidth]{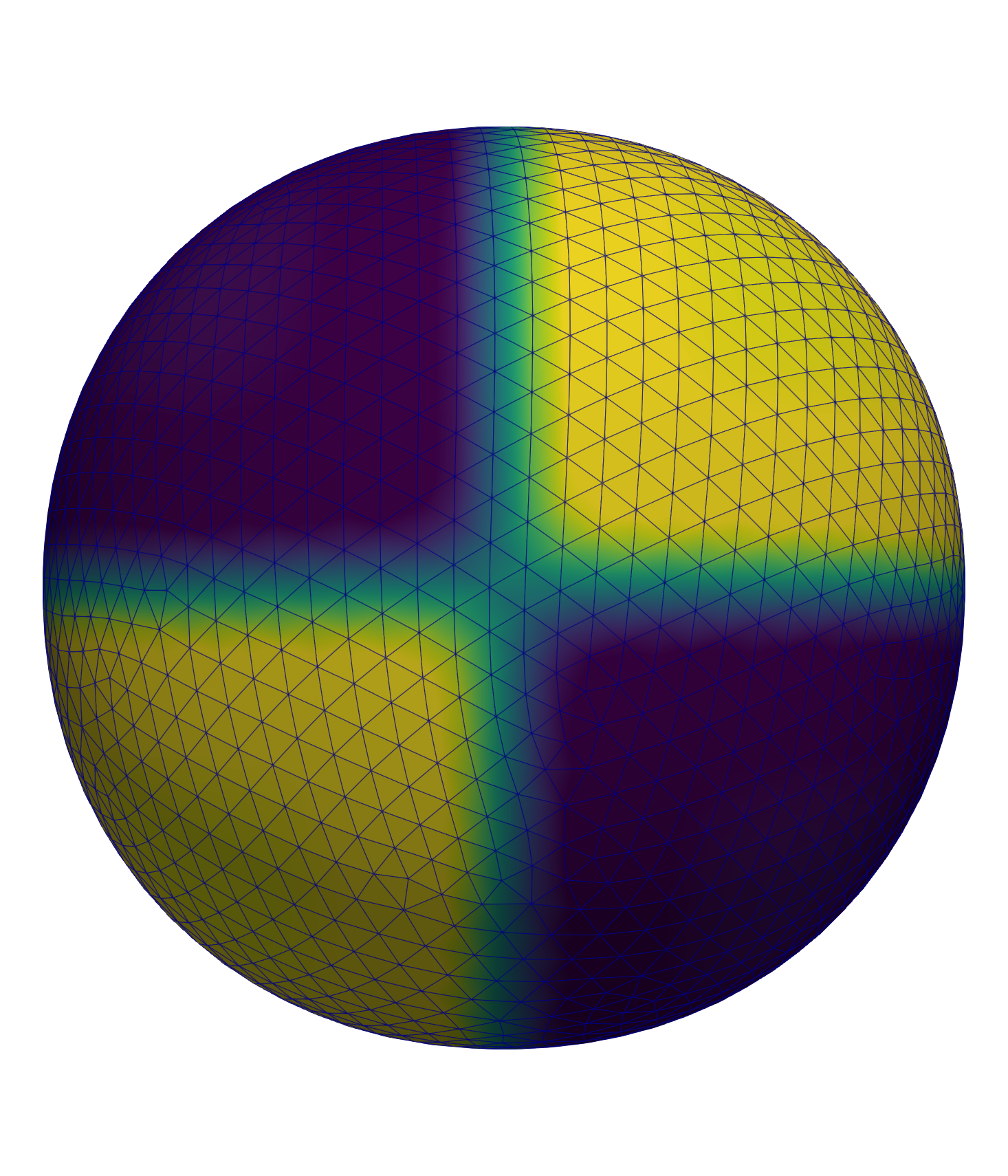}
        \caption{$\gamma_W=0.5, ~t=1$}
    \end{subfigure}

    \begin{subfigure}[tbhp!]{0.48\textwidth}
        \centering
        \includegraphics[clip, trim = {2cm 6cm 2cm 6cm}, width = 0.78\textwidth]{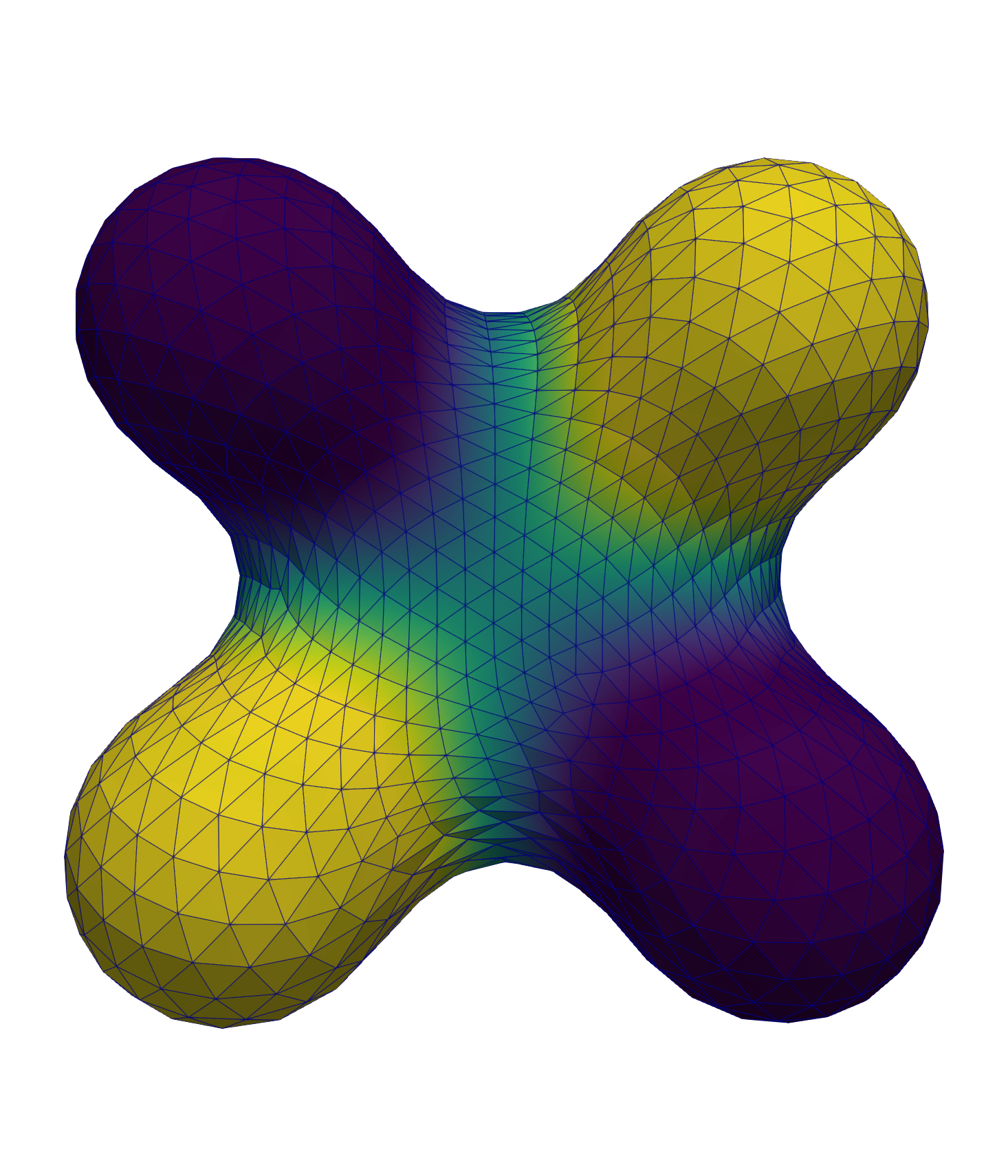}
        \caption{$\gamma_W=0.1, ~t=0.45$}
        \label{fig:mdr_bachini_trig_01}
    \end{subfigure}%
    ~
    \begin{subfigure}[tbhp!]{0.48\textwidth}
        \centering
        \includegraphics[clip, trim = {2cm 6cm 2cm 6cm}, width = 0.78\textwidth]{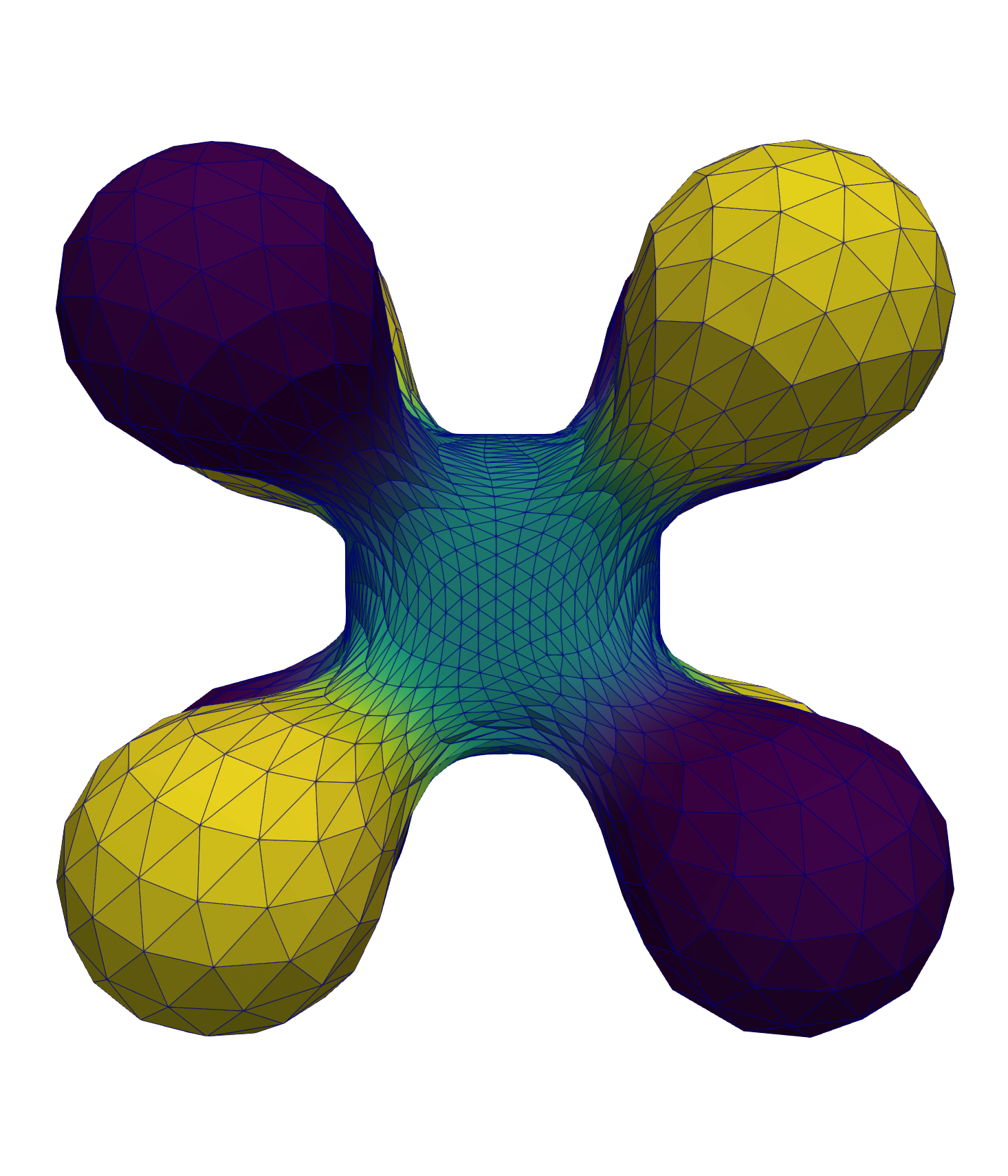}
        \caption{$\gamma_W=0.02, t=0.2375$}
    \end{subfigure}
    \caption{(a) Initial setup for the problem \ref{eq:wmch_bachini_u0}. (b)-(d) Mesh for different elasticity modulus. The simulation is shown at end time $t=1$ or at last timestep before inextensibility threshold is violated.}
    \label{fig:mdr_bachini_trig}
\end{figure}

We proceed now to test the coupling for more complicated systems and geometry evolution to explore the potential of our framework. The overdamped limit model for fluid deformable two-component membranes presented in \cite[p.A41-33]{BachiniKrauseNitschkeEtAl2023} is considered. The discrete scheme \eqref{eq:discrete_willmore} is augmented with a Lagrange multiplier to maintain the inextensibility constraint $\nabla_{\Gamma_h^n}\cdot\bv_h^n=0$. The modified version reads: Given $\by^{n}_h, \bkappa^{n}_h\in [V_h(\Gamma_h^n)]^d$ find $\bd^{n+1}_h\in [V_{h0}(\Gamma_h^n)]^d$, $\by^{n+1}_h\in [V_h(\Gamma_h^n)]^d$ and $\lambda_h^{n+1}\in V_{h}(\Gamma_h^n)$  such that
\begin{subequations}
\begin{align}
    &\inner[auto]{\frac{\bd_h^{n+1}}{\tau}}{\bphi_h}_{\Gamma_h^n}^h -\inner{\nabla_{\Gamma_h^n} \by^{n+1}_h}{\nabla_{\Gamma_h^n} \bphi_h}_{\Gamma_h^n} + \inner{\lambda_h^{n+1}\bn_{\Gamma_h^n}}{\bphi_h} \\
    & \qquad=~  p_h\inner{\by_h^n}{\bphi_h}, \nonumber \\
    &\frac{1}{\gamma_W}\inner{\by^{n+1}_h}{\bpsi_h}_{\Gamma_h^n}^h+ \inner{\nabla_{\Gamma_h^n}\bd_h^{n+1}}{\nabla_{\Gamma_h^n} \bpsi_h}_{\Gamma_h^n} \\
    & \qquad= -\inner{\bbP_{\Gamma_h^n}}{\nabla_{\Gamma_h^n} \bpsi_h}_{\Gamma_h^n}+ \inner{\kappa_0\bn_{\Gamma_h^n}}{\bpsi_h}_{\Gamma_h^n}^h + \inner{\bmu}{\bpsi_h}_{\partial\Gamma^n},\nonumber \\
    & \inner{\nabla_{\Gamma_h^n}\lambda_h^{n+1}}{\nabla_{\Gamma_h^n} \mu_h}_{\Gamma_h^n} + \inner{|\bkappa_h^n|^2\lambda_h^{n+1}}{\mu_h}_{\Gamma_h^n} \\
    & \qquad = \inner{\inner{\nabla_{\Gamma_h^n} \by^{n}_h}{\nabla_{\Gamma_h^n} \bkappa^n_h}_{\Gamma_h^n}}{\mu_h} + \inner{p_h\inner{\by_h^n}{\bkappa_h^n}}{\mu_h}, \nonumber
\end{align}
\label{eq:discrete_willmore_inex}
\end{subequations}
for all $\bphi_h\in [V_{h0}(\Gamma_h^n)]^d,~\bpsi_h\in [V_h(\Gamma_h^n)]^d,~ \mu_h \in V_{h}(\Gamma_h^n)$
where
\begin{align}
    &p_h\inner{\by_h^n}{\bphi_h} =\inner{\nabla_{\Gamma_h^n}\cdot \by^n_h}{\nabla_{\Gamma_h^n}\cdot \bphi_h}_{\Gamma_h^n}-\inner{(\nabla_{\Gamma_h^n}\by^n_h)^T}{\cD(\bphi_h)\;\bbP_{\Gamma_h^n}}_{\Gamma_h^n} \nonumber \\
    &\qquad- \kappa_0\inner{\bkappa^n_h}{(\nabla_{\Gamma_h^n}\bphi_h)^T \bn_{\Gamma_h^n}}_{\Gamma_h^n}^h -\frac{1}{2}\inner{\gamma_W|\bkappa^n_h - \kappa_0\bn_{\Gamma_h^n}|^2\;\bbP_{\Gamma_h^n}}{\nabla_{\Gamma_h^n} \bphi_h}_{\Gamma_h^n}^h \nonumber \\
    &\qquad + \inner{(\by^n_h\cdot\bkappa^n_h)\;\bbP_{\Gamma_h^n}}{\nabla_{\Gamma_h^n} \bphi_h}_{\Gamma_h^n}^h + \inner{\bf_h^{n+1}}{ \bphi_h}_{\Gamma_h^n}^h.
\label{eq:willmore_inex_auxiliary}
\end{align}
In \eqref{eq:willmore_inex_auxiliary}, $\bf_h^{n+1}$ is a general right-hand side. The system is composed of the following three steps:
\begin{enumerate}
    \item Solver \eqref{eq:discrete_willmore_inex} with the right-hand side (recall \eqref{eq:cahn_hillard_energy}): 
    \begin{equation}
        \bf_h^{n+1}=\cE_{CH}(u_h^{n+1})\bkappa_h^{n}-\sigma\epsilon \nabla_{\Gamma_h^n }u_h^{n+1}\bbH_h^n\nabla_{\Gamma_h^n} u_h^{n+1},
    \end{equation}
    where $\bomega_{\Gamma_h^n}$ is a certain vertex normal, see \cite{BarrettGarckeNurnberg2007} for details, and $\bbH_h^n = -\nabla_{\Gamma_h^n}\bomega_{\Gamma_h^n}$ is the discrete extended Weingarten map.
    \item Solver \eqref{eq:duanli_discrete} where the prescribed material velocity is set to
    \begin{equation}
        \bv_h^\top = -\nabla_{\Gamma_h^n}\lambda_h^{n+1}, \quad \bv_h^\perp = \bd_h^{n+1}/\tau.
    \end{equation}
    \item Solver \eqref{eq:cahn_hilliard_fully_discrete} with tangential advective velocity $\bv_h^{n+1}=\bv_h^\top$ and polynomial potential $F_2$ \eqref{eq:ch_potentials}.
\end{enumerate}
The test considers a sphere of radius one and the following initial parameters:
\vspace{2mm}
\begin{center}
\label{tab:wmch_bachini_params}
\begin{tabular}{||c|c|c|c|c||} 
\hline
$T$ & $h$ & $\sigma$ & $\varepsilon$ & $m$ \\ [0.5ex] 
\hline\hline
1 & 0.08 & $1.5\sqrt{2}$ & 0.1 & 0.001 \\ 
\hline
\end{tabular}
\end{center}
\vspace{2mm}
The common initial condition for the Cahn-Hilliard solvers is 
\begin{equation}
    u_0(\bx) = \cos(\pi x_1) \cos(\pi x_2) \cos(\pi x_3).
\label{eq:wmch_bachini_u0}
\end{equation}
The simulation is run with variable elasticity modulus $\gamma_W \in \{ 0.5, ~0.1, ~0.02 \}$ following the parameters in \cite{BachiniKrauseNitschkeEtAl2023}. The timestep is varied accordingly with values $\tau \in \{10^{-4},~10^{-4}, 2.5\cdot 10^{-5}\}$ respectively. The simulation is stopped once the area difference between initial and evolved surface differs by more than $10\%$. Results are presented in \autofigref{fig:bachini_trig}, where it can be seen that the characteristic bulging dynamics is correctly recovered. As already noted in \ref{subsec:helfrich_pinching} and \autofigref{fig:cigar}, the pearlings' mesh appears to be very coarse due to the extreme deformation. This is characteristic of the employed mesh redistribution algorithm presented in Section \ref{sec:mesh_redistribution}. As already noted in \cite{GarckeNurnbergZhao2025}, alternative choices are possible that might lead to better results \cite{Sauer2014, HuLi2022, PorrmannBartelsVoigt2025, Sauer2025, ZhangHuang2025, Duan2025}. Our setting is flexible enough to allow to switch between various options. 
We repeat the test in \ref{subsec:helfrich_pinching} substituting the surface mesh redistribution in Section \ref{sec:mesh_redistribution} \cite{DuanLi2024} with the minimal deformation rate (MDR) redistribution of \cite{HuLi2022} presented in \eqref{eq:huli_discrete}. The results are shown in \autofigref{fig:mdr_bachini_trig}. It can be seen how using the MDR technique leads, in this very specific case, to better results with respect to \autofigref{fig:bachini_trig}. This said, we found the performance of the different redistributions to be case-specific. For example, one can notice how, despite the enhanced node distribution, the inextensibility constaint for \autofigref{fig:mdr_bachini_trig_01} is violated at an earlier time with respect to \autofigref{fig:bachini_trig_01}.

\section{Conclusions and outlook}

Here, we presented a numerical framework to simulate moving boundary problems in biophysics. 
This work addresses a long-standing need for robust computational frameworks to tackle such problems. 
While other frameworks exist \cite{MokbelMokbelLieseEtAl2024, ArroyoDeSimone2009, BachiniKrauseNitschkeEtAl2023, MackenzieRowlattInsall2021, BarrettGarckeNurnberg2017}, our work targets flexibility of application while maintaining accuracy. Algorithms that enjoy stability and convergence properties are tied to non-invasive postprocessing techniques so to favor complex applications. The result is a stable pipeline with the capability of being scaled to more realistic scenarios.

Specifically, structure preservation is introduced on moving bulk and surface domains as a mean to allow for biophysically justified models and guarantee interpretability. 
We show  that the scheme is readily implementable while being fast and accurate.
The flexibility  of the algorithm, which does not require \emph{ad hoc} implementation, is proven by applying it to both advection-dominant ADR equations and phase-field models of the Cahn-Hilliard type. 
Convergence studies were performed for the bound- and mass-preserving scheme and the expected convergence rates were numerically verified. 
The framework also couples recent mesh redistribution techniques for gradient flows with the ALE method. 
This union is showed to be efficient in handling complex reshaping dynamics without the need for remeshing. 
As a result, the algorithm is tunable to the needs of the user and does not modify the underlying domain evolution. This is crucial in biophysical settings where the dynamics of the overall system is strongly coupled to the domains' shape.
Convergence studies were performed, demonstrating the accuracy of the scheme. Additionally, complex benchmark tests were conducted, showing good agreement with previous results in the literature. 
Building on previous results, we show that  staggering is an effective method to simulate increasingly complex physics such as those that occur in biophysics.

The framework developed here has strong potential to simulate problems with increasing biological complexity. 
To achieve this goal, we foresee the need to incorporate large systems of ADR equations, that are not discussed in the present work \cite{FrancisLaughlinDokkenEtAl2024, MackenzieRowlattInsall2021, MacDonaldMackenzieNolanEtAl2016}.
Furthermore, fluid equations have also been shown to play an important role at certain timescales \cite{Seifert1997, ShinBrangwynne2017, AlbertiHyman2021}. Including Navier-Stokes type equations on both interfaces and bulk \cite{MokbelMokbelLieseEtAl2024,BachiniKrauseNitschkeEtAl2023} will broaden the applicability of the computational schemes presented here.  

In summary, we presented a finite element framework able to handle a rich set of reshaping dynamics arising from bulk-surface coupled PDEs in biophysics. This framework poses itself as a major stepping stone on the path to physical simulations of moving boundary problems in cell biology. 

\section*{Acknowledgments}
We acknowledge Dr. Emmet Francis for his valuable input on biologically relevant models for computational biology and for sharing the lessons learned from his programming experience \cite{FrancisLaughlinDokkenEtAl2024}. All the simulations are performed using the software NgSolve \cite{Schoberl2014, GanglSturmNeunteufelEtAl2021}.

\section*{Declaration of Interests}
P.R. is a consultant for Simula Research Laboratories in Oslo, Norway and receives income. The terms of this arrangement have been reviewed and approved by the University of California, San Diego in accordance with its conflict-of-interest policies.

\section*{Funding sources}
P.R. supported by an NIH grant R35 GM158446.





\bibliographystyle{elsarticle-num}
\bibliography{references_contri,references_massing}

@article{AlandEgererLowengrubEtAl2014,
  title = {Diffuse Interface Models of Locally Inextensible Vesicles in a Viscous Fluid},
  author = {Aland, Sebastian and Egerer, Sabine and Lowengrub, John and Voigt, Axel},
  year = 2014,
  month = nov,
  journal = {Journal of Computational Physics},
  volume = {277},
  pages = {32--47},
  issn = {0021-9991},
  doi = {10.1016/j.jcp.2014.08.016},
  urldate = {2025-10-16}
}

@article{AlbertiHyman2021,
  title = {Biomolecular Condensates at the Nexus of Cellular Stress, Protein Aggregation Disease and Ageing},
  author = {Alberti, Simon and Hyman, Anthony A.},
  year = 2021,
  month = mar,
  journal = {Nat Rev Mol Cell Biol},
  volume = {22},
  number = {3},
  pages = {196--213},
  publisher = {Nature Publishing Group},
  issn = {1471-0080},
  doi = {10.1038/s41580-020-00326-6},
  urldate = {2025-10-16},
  copyright = {2021 Springer Nature Limited},
  langid = {english}
}

@article{AlphonseCaetanoDjurdjevacEtAl2023,
  title = {Function Spaces, Time Derivatives and Compactness for Evolving Families of {{Banach}} Spaces with Applications to {{PDEs}}},
  author = {Alphonse, Amal and Caetano, Diogo and Djurdjevac, Ana and Elliott, Charles M.},
  year = 2023,
  month = apr,
  journal = {Journal of Differential Equations},
  volume = {353},
  eprint = {2105.07908},
  primaryclass = {math},
  pages = {268--338},
  issn = {00220396},
  doi = {10/grp75n},
  urldate = {2025-04-25},
  archiveprefix = {arXiv}
}

@article{AlphonseElliottStinner2015,
  title = {An Abstract Framework for Parabolic {{PDEs}} on Evolving Spaces},
  author = {Alphonse, Amal and Elliott, Charles M. and Stinner, Bj{\"o}rn},
  year = 2015,
  month = mar,
  journal = {Portugaliae Mathematica},
  volume = {72},
  number = {1},
  pages = {1--46},
  issn = {0032-5155},
  doi = {10.4171/pm/1955},
  urldate = {2025-10-06},
  langid = {english}
}

@article{AmbroggioCostaNavarroPerezSocasEtAl2021,
  title = {Dengue and {{Zika}} Virus Capsid Proteins Bind to Membranes and Self-Assemble into Liquid Droplets with Nucleic Acids},
  author = {Ambroggio, Ernesto E. and Costa Navarro, Guadalupe S. and P{\'e}rez Socas, Luis Benito and Bagatolli, Luis A. and Gamarnik, Andrea V.},
  year = 2021,
  month = sep,
  journal = {J Biol Chem},
  volume = {297},
  number = {3},
  pages = {101059},
  issn = {1083-351X},
  doi = {10.1016/j.jbc.2021.101059},
  langid = {english},
  pmcid = {PMC8397897},
  pmid = {34375636}
}

@article{ArroyoDeSimone2009,
  title = {Relaxation Dynamics of Fluid Membranes},
  author = {Arroyo, Marino and DeSimone, Antonio},
  year = 2009,
  month = mar,
  journal = {Phys. Rev. E},
  volume = {79},
  number = {3},
  pages = {031915},
  publisher = {American Physical Society},
  doi = {10.1103/PhysRevE.79.031915},
  urldate = {2024-04-03}
}

@article{AuddyaZhangGulatiEtAl2021,
  title = {Biomembranes Undergo Complex, Non-Axisymmetric Deformations Governed by {{Kirchhoff}}--{{Love}} Kinematicsand Revealed by a Three-Dimensional Computational Framework},
  author = {Auddya, Debabrata and Zhang, Xiaoxuan and Gulati, Rahul and Vasan, Ritvik and Garikipati, Krishna and Rangamani, Padmini and Rudraraju, Shiva},
  year = 2021,
  month = nov,
  journal = {Proceedings of the Royal Society A: Mathematical, Physical and Engineering Sciences},
  volume = {477},
  number = {2255},
  pages = {20210246},
  publisher = {Royal Society},
  doi = {10.1098/rspa.2021.0246},
  urldate = {2025-10-29}
}

@article{BachiniKrauseNitschkeEtAl2023,
  title = {Derivation and Simulation of a Two-Phase Fluid Deformable Surface Model},
  author = {Bachini, Elena and Krause, Veit and Nitschke, Ingo and Voigt, Axel},
  year = 2023,
  journal = {Journal of Fluid Mechanics},
  volume = {977},
  pages = {A41},
  doi = {10.1017/jfm.2023.943}
}

@article{BadiaCodina2006,
  title = {Analysis of a Stabilized Finite Element Approximation of the Transient Convection-Diffusion Equation Using an Ale Framework},
  author = {Badia, Santiago and Codina, Ramon},
  year = 2006,
  month = jan,
  urldate = {2024-10-22},
  langid = {english}
}

@article{BaiHuLi2024,
  title = {A {{Convergent Evolving Finite Element Method}} with {{Artificial Tangential Motion}} for {{Surface Evolution}} under a {{Prescribed Velocity Field}}},
  author = {Bai, Genming and Hu, Jiashun and Li, Buyang},
  year = 2024,
  month = oct,
  journal = {SIAM J. Numer. Anal.},
  volume = {62},
  number = {5},
  pages = {2172--2195},
  publisher = {{Society for Industrial and Applied Mathematics}},
  issn = {0036-1429},
  doi = {10.1137/23M156968X},
  urldate = {2025-07-16}
}

@article{BananiLeeHymanEtAl2017,
  title = {Biomolecular Condensates: Organizers of Cellular Biochemistry},
  shorttitle = {Biomolecular Condensates},
  author = {Banani, Salman F. and Lee, Hyun O. and Hyman, Anthony A. and Rosen, Michael K.},
  year = 2017,
  month = may,
  journal = {Nat Rev Mol Cell Biol},
  volume = {18},
  number = {5},
  pages = {285--298},
  issn = {1471-0080},
  doi = {10.1038/nrm.2017.7},
  langid = {english},
  pmcid = {PMC7434221},
  pmid = {28225081}
}

@article{BarreiraElliottMadzvamuse2011,
  title = {The Surface Finite Element Method for Pattern Formation on Evolving Biological Surfaces},
  author = {Barreira, R. and Elliott, C. M. and Madzvamuse, A.},
  year = 2011,
  month = dec,
  journal = {J. Math. Biol.},
  volume = {63},
  number = {6},
  pages = {1095--1119},
  issn = {1432-1416},
  doi = {10/b386pd},
  urldate = {2024-09-02},
  langid = {english}
}

@incollection{BarrettGarckeNurnberg_bookSection,
  title = {Chapter 4 - {{Parametric}} Finite Element Approximations of Curvature-Driven Interface Evolutions},
  booktitle = {Handbook of {{Numerical Analysis}}},
  author = {Barrett, John W. and Garcke, Harald and N{\"u}rnberg, Robert},
  editor = {Bonito, Andrea and Nochetto, Ricardo H.},
  year = 2020,
  month = jan,
  series = {Geometric {{Partial Differential Equations}} - {{Part I}}},
  volume = {21},
  pages = {275--423},
  publisher = {Elsevier},
  doi = {10.1016/bs.hna.2019.05.002},
  urldate = {2024-10-22}
}

@article{BarrettGarckeNurnberg2007,
  title = {A Parametric Finite Element Method for Fourth Order Geometric Evolution Equations},
  author = {Barrett, John W. and Garcke, Harald and N{\"u}rnberg, Robert},
  year = 2007,
  month = mar,
  journal = {Journal of Computational Physics},
  volume = {222},
  number = {1},
  pages = {441--467},
  issn = {0021-9991},
  doi = {10/cr3n7j},
  urldate = {2024-10-22}
}

@article{BarrettGarckeNurnberg2016,
  title = {Computational {{Parametric Willmore Flow}} with {{Spontaneous Curvature}} and {{Area Difference Elasticity Effects}}},
  author = {Barrett, John W. and Garcke, Harald and N{\"u}rnberg, Robert},
  year = 2016,
  month = jan,
  journal = {SIAM J. Numer. Anal.},
  volume = {54},
  number = {3},
  pages = {1732--1762},
  publisher = {{Society for Industrial and Applied Mathematics}},
  issn = {0036-1429},
  doi = {10/f865cb},
  urldate = {2024-10-09}
}

@article{BarrettGarckeNurnberg2017,
  title = {Finite {{Element Approximation}} for the {{Dynamics}} of {{Fluidic Two-Phase Biomembranes}}},
  author = {Barrett, John W. and Garcke, Harald and N{\"u}rnberg, Robert},
  year = 2017,
  month = nov,
  journal = {ESAIM: M2AN},
  volume = {51},
  number = {6},
  eprint = {1611.05343},
  primaryclass = {math},
  pages = {2319--2366},
  issn = {0764-583X, 1290-3841},
  doi = {10/gsj8sk},
  urldate = {2025-05-30},
  archiveprefix = {arXiv}
}

@article{BarrettGarckeNurnberg2017a,
  title = {Stable Variational Approximations of Boundary Value Problems for {{Willmore}} Flow with {{Gaussian}} Curvature},
  author = {Barrett, John W and Garcke, Harald and N{\"u}rnberg, Robert},
  year = 2017,
  month = oct,
  journal = {IMA Journal of Numerical Analysis},
  volume = {37},
  number = {4},
  pages = {1657--1709},
  issn = {0272-4979},
  doi = {10/gb4tvx},
  urldate = {2024-10-09}
}

@article{BaumgartHessWebb2003,
  title = {Imaging Coexisting Fluid Domains in Biomembrane Models Coupling Curvature and Line Tension},
  author = {Baumgart, Tobias and Hess, Samuel T. and Webb, Watt W.},
  year = 2003,
  month = oct,
  journal = {Nature},
  volume = {425},
  number = {6960},
  pages = {821--824},
  publisher = {Nature Publishing Group},
  issn = {1476-4687},
  doi = {10.1038/nature02013},
  urldate = {2025-06-30},
  copyright = {2003 Macmillan Magazines Ltd.},
  langid = {english}
}

@article{BeutelMaraspiniPombo-GarciaEtAl2019,
  title = {Phase {{Separation}} of {{Zonula Occludens Proteins Drives Formation}} of {{Tight Junctions}}},
  author = {Beutel, Oliver and Maraspini, Riccardo and {Pombo-Garc{\'i}a}, Karina and {Martin-Lemaitre}, C{\'e}cilie and Honigmann, Alf},
  year = 2019,
  month = oct,
  journal = {Cell},
  volume = {179},
  number = {4},
  pages = {923-936.e11},
  issn = {1097-4172},
  doi = {10.1016/j.cell.2019.10.011},
  langid = {english},
  pmid = {31675499}
}

@article{Bonilla-QuintanaRangamani2024,
  title = {Biophysical {{Modeling}} of {{Actin-Mediated Structural Plasticity Reveals Mechanical Adaptation}} in {{Dendritic Spines}}},
  author = {{Bonilla-Quintana}, Mayte and Rangamani, Padmini},
  year = 2024,
  month = mar,
  journal = {eNeuro},
  volume = {11},
  number = {3},
  publisher = {Society for Neuroscience},
  issn = {2373-2822},
  doi = {10.1523/ENEURO.0497-23.2024},
  urldate = {2024-04-01},
  chapter = {Research Article: New Research},
  copyright = {Copyright \copyright{} 2024 Bonilla-Quintana and Rangamani. This is an open-access article distributed under the terms of the Creative Commons Attribution 4.0 International license, which permits unrestricted use, distribution and reproduction in any medium provided that the original work is properly attributed.},
  langid = {english},
  pmid = {38383589}
}

@article{BonitoKyzaNochetto2013,
  title = {Time-Discrete Higher Order {{ALE}} Formulations: A Priori Error Analysis},
  shorttitle = {Time-Discrete Higher Order {{ALE}} Formulations},
  author = {Bonito, Andrea and Kyza, Irene and Nochetto, Ricardo H.},
  year = 2013,
  month = oct,
  journal = {Numer. Math.},
  volume = {125},
  number = {2},
  pages = {225--257},
  issn = {0945-3245},
  doi = {10/f49ctb},
  urldate = {2025-04-15},
  langid = {english}
}

@article{BrangwynneEckmannCoursonEtAl2009,
  title = {Germline {{P Granules Are Liquid Droplets That Localize}} by {{Controlled Dissolution}}/{{Condensation}}},
  author = {Brangwynne, Clifford P. and Eckmann, Christian R. and Courson, David S. and Rybarska, Agata and Hoege, Carsten and Gharakhani, J{\"o}bin and J{\"u}licher, Frank and Hyman, Anthony A.},
  year = 2009,
  month = jun,
  journal = {Science},
  volume = {324},
  number = {5935},
  pages = {1729--1732},
  publisher = {American Association for the Advancement of Science},
  doi = {10.1126/science.1172046},
  urldate = {2025-10-16}
}

@article{BurmanClausHansboEtAl2015,
  title = {{{CutFEM}}: {{Discretizing}} Geometry and Partial Differential Equations},
  shorttitle = {{{CutFEM}}},
  author = {Burman, Erik and Claus, Susanne and Hansbo, Peter and Larson, Mats G. and Massing, Andr{\'e}},
  year = 2015,
  journal = {International Journal for Numerical Methods in Engineering},
  volume = {104},
  number = {7},
  pages = {472--501},
  issn = {1097-0207},
  doi = {10.1002/nme.4823},
  urldate = {2025-10-14},
  copyright = {Copyright \copyright{} 2014 John Wiley \& Sons, Ltd.},
  langid = {english}
}

@article{CaetanoElliott2021,
  title = {Cahn--{{Hilliard}} Equations on an Evolving Surface},
  author = {Caetano, D. and Elliott, C. M.},
  year = 2021,
  month = oct,
  journal = {European Journal of Applied Mathematics},
  volume = {32},
  number = {5},
  pages = {937--1000},
  issn = {0956-7925, 1469-4425},
  doi = {10.1017/S0956792521000176},
  urldate = {2025-10-06},
  langid = {english}
}

@article{CaetanoElliottGrasselliEtAl2023,
  title = {Regularization and {{Separation}} for {{Evolving Surface Cahn}}--{{Hilliard Equations}}},
  author = {Caetano, Diogo and Elliott, Charles M. and Grasselli, Maurizio and Poiatti, Andrea},
  year = 2023,
  month = dec,
  journal = {SIAM J. Math. Anal.},
  volume = {55},
  number = {6},
  pages = {6625--6675},
  publisher = {{Society for Industrial and Applied Mathematics}},
  issn = {0036-1410},
  doi = {10/g74q2x},
  urldate = {2025-04-17}
}

@article{CahnHilliard1958,
  title = {Free {{Energy}} of a {{Nonuniform System}}. {{I}}. {{Interfacial Free Energy}}},
  author = {Cahn, John W. and Hilliard, John E.},
  year = 1958,
  month = feb,
  journal = {J. Chem. Phys.},
  volume = {28},
  number = {2},
  pages = {258--267},
  issn = {0021-9606},
  doi = {10.1063/1.1744102},
  urldate = {2025-10-06}
}

@article{Canham1970,
  title = {The Minimum Energy of Bending as a Possible Explanation of the Biconcave Shape of the Human Red Blood Cell},
  author = {Canham, P. B.},
  year = 1970,
  month = jan,
  journal = {Journal of Theoretical Biology},
  volume = {26},
  number = {1},
  pages = {61--81},
  issn = {0022-5193},
  doi = {10/frc472},
  urldate = {2024-09-23}
}

@article{ChengShen2022,
  title = {A {{New Lagrange Multiplier Approach}} for {{Constructing Structure Preserving Schemes}}, {{II}}. {{Bound Preserving}}},
  author = {Cheng, Qing and Shen, Jie},
  year = 2022,
  month = jun,
  journal = {SIAM J. Numer. Anal.},
  volume = {60},
  number = {3},
  pages = {970--998},
  publisher = {{Society for Industrial and Applied Mathematics}},
  issn = {0036-1429},
  doi = {10/gr5qjs},
  urldate = {2025-02-05}
}

@article{ChengShen2022a,
  title = {A New {{Lagrange}} Multiplier Approach for Constructing Structure Preserving Schemes, {{I}}. {{Positivity}} Preserving},
  author = {Cheng, Qing and Shen, Jie},
  year = 2022,
  month = mar,
  journal = {Computer Methods in Applied Mechanics and Engineering},
  volume = {391},
  pages = {114585},
  issn = {0045-7825},
  doi = {10/gr5qjr},
  urldate = {2025-02-05}
}

@article{ChenWangWangEtAl2019,
  title = {Positivity-Preserving, Energy Stable Numerical Schemes for the {{Cahn-Hilliard}} Equation with Logarithmic Potential},
  author = {Chen, Wenbin and Wang, Cheng and Wang, Xiaoming and Wise, Steven M.},
  year = 2019,
  month = jun,
  journal = {Journal of Computational Physics: X},
  volume = {3},
  pages = {100031},
  issn = {2590-0552},
  doi = {10.1016/j.jcpx.2019.100031},
  urldate = {2025-07-16}
}

@article{DayelAkinLanderyouEtAl2009,
  title = {In {{Silico Reconstitution}} of {{Actin-Based Symmetry Breaking}} and {{Motility}}},
  author = {Dayel, Mark J. and Akin, Orkun and Landeryou, Mark and Risca, Viviana and Mogilner, Alex and Mullins, R. Dyche},
  year = 2009,
  month = sep,
  journal = {PLOS Biology},
  volume = {7},
  number = {9},
  pages = {e1000201},
  publisher = {Public Library of Science},
  issn = {1545-7885},
  doi = {10.1371/journal.pbio.1000201},
  urldate = {2025-10-21},
  langid = {english}
}

@article{DeckelnickDziukElliott2005,
  title = {Computation of Geometric Partial Differential Equations and Mean Curvature Flow},
  author = {Deckelnick, Klaus and Dziuk, Gerhard and Elliott, Charles M.},
  year = 2005,
  month = may,
  journal = {Acta Numerica},
  volume = {14},
  pages = {139--232},
  issn = {1474-0508, 0962-4929},
  doi = {10/c388gg},
  urldate = {2024-10-22},
  langid = {english}
}

@article{DoubrovinskiKruse2011,
  title = {Cell Motility Resulting from Spontaneous Polymerization Waves},
  author = {Doubrovinski, K. and Kruse, K.},
  year = 2011,
  month = dec,
  journal = {Phys Rev Lett},
  volume = {107},
  number = {25},
  pages = {258103},
  issn = {1079-7114},
  doi = {10.1103/PhysRevLett.107.258103},
  langid = {english},
  pmid = {22243118}
}

@article{Duan2025,
  title = {Mesh-{{Preserving}} and {{Energy-Stable Parametric FEM}} for {{Geometric Flows}} of {{Surfaces}}},
  author = {Duan, Beiping},
  year = 2025,
  month = apr,
  journal = {SIAM J. Numer. Anal.},
  volume = {63},
  number = {2},
  pages = {619--640},
  publisher = {{Society for Industrial and Applied Mathematics}},
  issn = {0036-1429},
  doi = {10.1137/24M1671542},
  urldate = {2025-06-16}
}

@article{DuanLi2024,
  title = {New {{Artificial Tangential Motions}} for {{Parametric Finite Element Approximation}} of {{Surface Evolution}}},
  author = {Duan, Beiping and Li, Buyang},
  year = 2024,
  month = feb,
  journal = {SIAM J. Sci. Comput.},
  volume = {46},
  number = {1},
  pages = {A587-A608},
  publisher = {{Society for Industrial and Applied Mathematics}},
  issn = {1064-8275},
  doi = {10/gtj6hx},
  urldate = {2024-05-02}
}

@article{DuJuLiEtAl2021,
  title = {Maximum {{Bound Principles}} for a {{Class}} of {{Semilinear Parabolic Equations}} and {{Exponential Time-Differencing Schemes}}},
  author = {Du, Qiang and Ju, Lili and Li, Xiao and Qiao, Zhonghua},
  year = 2021,
  month = jan,
  journal = {SIAM Rev.},
  volume = {63},
  number = {2},
  pages = {317--359},
  publisher = {{Society for Industrial and Applied Mathematics}},
  issn = {0036-1445},
  doi = {10/gp9fzj},
  urldate = {2025-03-26}
}

@article{DuJuTian2011,
  title = {Finite Element Approximation of the {{Cahn}}--{{Hilliard}} Equation on Surfaces},
  author = {Du, Qiang and Ju, Lili and Tian, Li},
  year = 2011,
  month = jul,
  journal = {Computer Methods in Applied Mechanics and Engineering},
  volume = {200},
  number = {29},
  pages = {2458--2470},
  issn = {0045-7825},
  doi = {10.1016/j.cma.2011.04.018},
  urldate = {2025-10-06}
}

@article{Dziuk2008,
  title = {Computational Parametric {{Willmore}} Flow},
  author = {Dziuk, Gerhard},
  year = 2008,
  month = nov,
  journal = {Numer. Math.},
  volume = {111},
  number = {1},
  pages = {55--80},
  issn = {0945-3245},
  doi = {10/cgpjqt},
  urldate = {2024-05-01},
  langid = {english}
}

@article{DziukElliott2007,
  title = {Finite Elements on Evolving Surfaces},
  author = {Dziuk, G. and Elliott, C. M.},
  year = 2007,
  month = apr,
  journal = {IMA Journal of Numerical Analysis},
  volume = {27},
  number = {2},
  pages = {262--292},
  issn = {0272-4979},
  doi = {10/c4x9hd},
  urldate = {2025-04-16}
}

@article{DziukElliott2013,
  title = {Finite Element Methods for Surface {{PDEs}}},
  author = {Dziuk, Gerhard and Elliott, Charles M.},
  year = 2013,
  month = may,
  journal = {Acta Numerica},
  volume = {22},
  pages = {289--396},
  issn = {0962-4929, 1474-0508},
  doi = {10/ggzzsq},
  urldate = {2024-09-23},
  langid = {english}
}

@article{DziukLubichMansour2012,
  title = {Runge--{{Kutta}} Time Discretization of Parabolic Differential Equations on Evolving Surfaces},
  author = {Dziuk, G. and Lubich, {\relax Ch}. and Mansour, D.},
  year = 2012,
  month = apr,
  journal = {IMA J Numer Anal},
  volume = {32},
  number = {2},
  pages = {394--416},
  issn = {0272-4979},
  doi = {10.1093/imanum/drr017},
  urldate = {2025-10-25}
}

@article{EilksElliott2008,
  title = {Numerical Simulation of Dealloying by Surface Dissolution via the Evolving Surface Finite Element Method},
  author = {Eilks, C. and Elliott, C. M.},
  year = 2008,
  month = dec,
  journal = {Journal of Computational Physics},
  volume = {227},
  number = {23},
  pages = {9727--9741},
  issn = {0021-9991},
  doi = {10.1016/j.jcp.2008.07.023},
  urldate = {2025-10-06}
}

@article{ElliottRanner2021,
  title = {A Unified Theory for Continuous-in-Time Evolving Finite Element Space Approximations to Partial Differential Equations in Evolving Domains},
  author = {Elliott, C M and Ranner, T},
  year = 2021,
  month = jul,
  journal = {IMA Journal of Numerical Analysis},
  volume = {41},
  number = {3},
  pages = {1696--1845},
  issn = {0272-4979},
  doi = {10/gr8rfz},
  urldate = {2024-10-10}
}

@article{ElliottStinner2010,
  title = {Modeling and Computation of Two Phase Geometric Biomembranes Using Surface Finite Elements},
  author = {Elliott, Charles M. and Stinner, Bj{\"o}rn},
  year = 2010,
  month = sep,
  journal = {Journal of Computational Physics},
  volume = {229},
  number = {18},
  pages = {6585--6612},
  issn = {0021-9991},
  doi = {10.1016/j.jcp.2010.05.014},
  urldate = {2025-10-06}
}

@article{ErlebacherAzizKarmaEtAl2001,
  title = {Evolution of Nanoporosity in Dealloying},
  author = {Erlebacher, Jonah and Aziz, Michael J. and Karma, Alain and Dimitrov, Nikolay and Sieradzki, Karl},
  year = 2001,
  month = mar,
  journal = {Nature},
  volume = {410},
  number = {6827},
  pages = {450--453},
  publisher = {Nature Publishing Group},
  issn = {1476-4687},
  doi = {10.1038/35068529},
  urldate = {2025-10-06},
  copyright = {2001 Macmillan Magazines Ltd.},
  langid = {english}
}

@book{ErnGuermond_book,
  title = {Finite {{Elements III}}: {{First-Order}} and {{Time-Dependent PDEs}}},
  shorttitle = {Finite {{Elements III}}},
  author = {Ern, Alexandre and Guermond, Jean-Luc},
  year = 2021,
  series = {Texts in {{Applied Mathematics}}},
  volume = {74},
  publisher = {Springer International Publishing},
  address = {Cham},
  doi = {10.1007/978-3-030-57348-5},
  urldate = {2024-10-03},
  copyright = {http://www.springer.com/tdm},
  isbn = {978-3-030-57347-8 978-3-030-57348-5},
  langid = {english}
}

@book{ErnGuermond_booka,
  title = {Finite {{Elements II}}: {{Galerkin Approximation}}, {{Elliptic}} and {{Mixed PDEs}}},
  shorttitle = {Finite {{Elements II}}},
  author = {Ern, Alexandre and Guermond, Jean-Luc},
  year = 2021,
  series = {Texts in {{Applied Mathematics}}},
  volume = {73},
  publisher = {Springer International Publishing},
  address = {Cham},
  doi = {10.1007/978-3-030-56923-5},
  urldate = {2024-10-03},
  copyright = {https://www.springer.com/tdm},
  isbn = {978-3-030-56922-8 978-3-030-56923-5},
  langid = {english}
}

@book{ErnGuermond_bookb,
  title = {Finite {{Elements I}}: {{Approximation}} and {{Interpolation}}},
  shorttitle = {Finite {{Elements I}}},
  author = {Ern, Alexandre and Guermond, Jean-Luc},
  year = 2021,
  series = {Texts in {{Applied Mathematics}}},
  volume = {72},
  publisher = {Springer International Publishing},
  address = {Cham},
  doi = {10.1007/978-3-030-56341-7},
  urldate = {2024-10-03},
  copyright = {http://www.springer.com/tdm},
  isbn = {978-3-030-56340-0 978-3-030-56341-7},
  langid = {english}
}

@article{Flory1942,
  title = {Thermodynamics of {{High Polymer Solutions}}},
  author = {Flory, Paul J.},
  year = 1942,
  month = jan,
  journal = {J. Chem. Phys.},
  volume = {10},
  number = {1},
  pages = {51--61},
  issn = {0021-9606},
  doi = {10.1063/1.1723621},
  urldate = {2025-10-21}
}

@article{FormaggiaNobile2004,
  title = {Stability Analysis of Second-Order Time Accurate Schemes for {{ALE}}--{{FEM}}},
  author = {Formaggia, Luca and Nobile, Fabio},
  year = 2004,
  month = oct,
  journal = {Computer Methods in Applied Mechanics and Engineering},
  series = {The {{Arbitrary Lagrangian-Eulerian Formulation}}},
  volume = {193},
  number = {39},
  pages = {4097--4116},
  issn = {0045-7825},
  doi = {10/cvf4fz},
  urldate = {2024-10-18}
}

@misc{FrancisLaughlinDokkenEtAl2024,
  title = {Spatial Modeling Algorithms for Reactions and Transport ({{SMART}}) in Biological Cells},
  author = {Francis, Emmet A. and Laughlin, Justin and Dokken, J{\o}rgen S. and Finsberg, Henrik and Lee, Christopher T. and Rognes, Marie E. and Rangamani, Padmini},
  year = 2024,
  month = may,
  primaryclass = {New Results},
  pages = {2024.05.23.595604},
  publisher = {bioRxiv},
  doi = {10.1101/2024.05.23.595604},
  urldate = {2024-10-09},
  archiveprefix = {bioRxiv},
  chapter = {New Results},
  copyright = {\copyright{} 2024, Posted by Cold Spring Harbor Laboratory. This pre-print is available under a Creative Commons License (Attribution 4.0 International), CC BY 4.0, as described at http://creativecommons.org/licenses/by/4.0/},
  langid = {english}
}

@article{GanglSturmNeunteufelEtAl2021,
  title = {Fully and Semi-Automated Shape Differentiation in {{NGSolve}}},
  author = {Gangl, Peter and Sturm, Kevin and Neunteufel, Michael and Sch{\"o}berl, Joachim},
  year = 2021,
  month = mar,
  journal = {Struct Multidisc Optim},
  volume = {63},
  number = {3},
  pages = {1579--1607},
  issn = {1615-1488},
  doi = {10.1007/s00158-020-02742-w},
  urldate = {2025-10-29},
  langid = {english}
}

@misc{GarckeNurnbergZhao2025,
  title = {An Energy-Stable Parametric Finite Element Method for {{Willmore}} Flow with Normal-Tangential Velocity Splitting},
  author = {Garcke, Harald and N{\"u}rnberg, Robert and Zhao, Quan},
  year = 2025,
  month = jun,
  number = {arXiv:2507.00193},
  eprint = {2507.00193},
  primaryclass = {math},
  publisher = {arXiv},
  doi = {10.48550/arXiv.2507.00193},
  urldate = {2025-07-09},
  archiveprefix = {arXiv}
}

@misc{GargOlshanskii2025,
  title = {A {{Stabilized Trace FEM}} for {{Surface Cahn--Hilliard Equations}}: {{Analysis}} and {{Simulations}}},
  shorttitle = {A {{Stabilized Trace FEM}} for {{Surface Cahn--Hilliard Equations}}},
  author = {Garg, Deepika and Olshanskii, Maxim},
  year = 2025,
  month = oct,
  number = {arXiv:2510.21662},
  eprint = {2510.21662},
  primaryclass = {math},
  publisher = {arXiv},
  doi = {10.48550/arXiv.2510.21662},
  urldate = {2025-10-27},
  archiveprefix = {arXiv}
}

@article{Gastaldi2001,
  title = {A Priori Error Estimates for the {{Arbitrary Lagrangian Eulerian}} Formulation with Finite Elements},
  author = {Gastaldi, L.},
  year = 2001,
  month = jun,
  volume = {9},
  number = {2},
  pages = {123--156},
  publisher = {De Gruyter},
  issn = {1569-3953},
  doi = {10.1515/JNMA.2001.123},
  urldate = {2025-10-26},
  chapter = {Journal of Numerical Mathematics},
  copyright = {De Gruyter expressly reserves the right to use all content for commercial text and data mining within the meaning of Section 44b of the German Copyright Act.},
  langid = {english}
}

@misc{GulatiRudraraju2024,
  title = {Electro-Diffusive Modeling and the Role of Spine Geometry on Action Potential Propagation in Neurons},
  author = {Gulati, Rahul and Rudraraju, Shiva},
  year = 2024,
  month = nov,
  number = {arXiv:2411.05329},
  eprint = {2411.05329},
  primaryclass = {q-bio},
  publisher = {arXiv},
  doi = {10.48550/arXiv.2411.05329},
  urldate = {2025-10-29},
  archiveprefix = {arXiv}
}

@article{Helfrich1973,
  title = {Elastic {{Properties}} of {{Lipid Bilayers}}: {{Theory}} and {{Possible Experiments}}},
  shorttitle = {Elastic {{Properties}} of {{Lipid Bilayers}}},
  author = {Helfrich, W.},
  year = 1973,
  month = dec,
  journal = {Zeitschrift f\"ur Naturforschung C},
  volume = {28},
  number = {11-12},
  pages = {693--703},
  publisher = {De Gruyter},
  issn = {1865-7125},
  doi = {10/gf3jnf},
  urldate = {2024-09-23},
  copyright = {De Gruyter expressly reserves the right to use all content for commercial text and data mining within the meaning of Section 44b of the German Copyright Act.},
  langid = {english}
}

@article{HirtAmsdenCook1974,
  title = {An Arbitrary {{Lagrangian-Eulerian}} Computing Method for All Flow Speeds},
  author = {Hirt, C. W and Amsden, A. A and Cook, J. L},
  year = 1974,
  month = mar,
  journal = {Journal of Computational Physics},
  volume = {14},
  number = {3},
  pages = {227--253},
  issn = {0021-9991},
  doi = {10/fdtb8s},
  urldate = {2025-04-16}
}

@article{HuangShen2021,
  title = {Bound/{{Positivity Preserving}} and {{Energy Stable Scalar}} Auxiliary {{Variable Schemes}} for {{Dissipative Systems}}: {{Applications}} to {{Keller--Segel}} and {{Poisson--Nernst--Planck Equations}}},
  shorttitle = {Bound/{{Positivity Preserving}} and {{Energy Stable Scalar}} Auxiliary {{Variable Schemes}} for {{Dissipative Systems}}},
  author = {Huang, Fukeng and Shen, Jie},
  year = 2021,
  month = jan,
  journal = {SIAM J. Sci. Comput.},
  volume = {43},
  number = {3},
  pages = {A1832-A1857},
  publisher = {{Society for Industrial and Applied Mathematics}},
  issn = {1064-8275},
  doi = {10/gsskvn},
  urldate = {2025-05-27}
}

@article{Huggins1941,
  title = {Solutions of {{Long Chain Compounds}}},
  author = {Huggins, Maurice L.},
  year = 1941,
  month = may,
  journal = {J. Chem. Phys.},
  volume = {9},
  number = {5},
  pages = {440},
  issn = {0021-9606},
  doi = {10.1063/1.1750930},
  urldate = {2025-10-21}
}

@article{HughesLiuZimmermann1981,
  title = {Lagrangian-{{Eulerian}} Finite Element Formulation for Incompressible Viscous Flows},
  author = {Hughes, Thomas J. R. and Liu, Wing Kam and Zimmermann, Thomas K.},
  year = 1981,
  month = dec,
  journal = {Computer Methods in Applied Mechanics and Engineering},
  volume = {29},
  number = {3},
  pages = {329--349},
  issn = {0045-7825},
  doi = {10/fjv3wd},
  urldate = {2024-10-22}
}

@article{HuHuang2020,
  title = {A Fully Discrete Positivity-Preserving and Energy-Dissipative Finite Difference Scheme for {{Poisson}}--{{Nernst}}--{{Planck}} Equations},
  author = {Hu, Jingwei and Huang, Xiaodong},
  year = 2020,
  month = may,
  journal = {Numer. Math.},
  volume = {145},
  number = {1},
  pages = {77--115},
  issn = {0945-3245},
  doi = {10/ggw36v},
  urldate = {2025-05-27},
  langid = {english}
}

@article{HuLi2022,
  title = {Evolving Finite Element Methods with an Artificial Tangential Velocity for Mean Curvature Flow and {{Willmore}} Flow},
  author = {Hu, Jiashun and Li, Buyang},
  year = 2022,
  month = sep,
  journal = {Numer. Math.},
  volume = {152},
  number = {1},
  pages = {127--181},
  issn = {0945-3245},
  doi = {10/gr8rf4},
  urldate = {2024-04-26},
  langid = {english}
}

@article{IyerGompperFedosov2023,
  title = {Dynamic Shapes of Floppy Vesicles Enclosing Active {{Brownian}} Particles with Membrane Adhesion},
  author = {Iyer, Priyanka and Gompper, Gerhard and Fedosov, Dmitry A.},
  year = 2023,
  month = may,
  journal = {Soft Matter},
  volume = {19},
  number = {19},
  pages = {3436--3449},
  publisher = {The Royal Society of Chemistry},
  issn = {1744-6848},
  doi = {10.1039/D3SM00004D},
  urldate = {2025-10-16},
  langid = {english}
}

@book{Karush_book,
  title = {Minima of {{Functions}} of {{Several Variables}} with {{Inequalities}} as {{Side Conditions}}},
  author = {Karush, William},
  year = 1939,
  publisher = {University of Chicago},
  googlebooks = {Z9I2zwEACAAJ},
  langid = {english}
}

@article{KerenPincusAllenEtAl2008,
  title = {Mechanism of Shape Determination in Motile Cells},
  author = {Keren, Kinneret and Pincus, Zachary and Allen, Greg M. and Barnhart, Erin L. and Marriott, Gerard and Mogilner, Alex and Theriot, Julie A.},
  year = 2008,
  month = may,
  journal = {Nature},
  volume = {453},
  number = {7194},
  pages = {475--480},
  issn = {1476-4687},
  doi = {10.1038/nature06952},
  langid = {english},
  pmcid = {PMC2877812},
  pmid = {18497816}
}

@article{KhanwaleLofquistSundarEtAl2020,
  title = {Simulating Two-Phase Flows with Thermodynamically Consistent Energy Stable {{Cahn-Hilliard Navier-Stokes}} Equations on Parallel Adaptive Octree Based Meshes},
  author = {Khanwale, Makrand A. and Lofquist, Alec D. and Sundar, Hari and Rossmanith, James A. and Ganapathysubramanian, Baskar},
  year = 2020,
  month = oct,
  journal = {Journal of Computational Physics},
  volume = {419},
  pages = {109674},
  issn = {0021-9991},
  doi = {10.1016/j.jcp.2020.109674},
  urldate = {2025-10-14}
}

@article{KhanwaleSaurabhIshiiEtAl2023,
  title = {A Projection-Based, Semi-Implicit Time-Stepping Approach for the {{Cahn-Hilliard Navier-Stokes}} Equations on Adaptive Octree Meshes},
  author = {Khanwale, Makrand A. and Saurabh, Kumar and Ishii, Masado and Sundar, Hari and Rossmanith, James A. and Ganapathysubramanian, Baskar},
  year = 2023,
  month = feb,
  journal = {Journal of Computational Physics},
  volume = {475},
  pages = {111874},
  issn = {0021-9991},
  doi = {10.1016/j.jcp.2022.111874},
  urldate = {2025-10-14}
}

@article{Kirchhoff1850,
  title = {Ueber Die {{Schwingungen}} Einer Kreisf\"ormigen Elastischen {{Scheibe}}},
  author = {Kirchhoff, G.},
  year = 1850,
  journal = {Annalen der Physik},
  volume = {157},
  number = {10},
  pages = {258--264},
  issn = {1521-3889},
  doi = {10.1002/andp.18501571005},
  urldate = {2025-06-30},
  copyright = {Copyright \copyright{} 1850 WILEY-VCH Verlag GmbH \& Co. KGaA, Weinheim},
  langid = {english}
}

@article{KovacsLiLubichEtAl2017,
  title = {Convergence of Finite Elements on an Evolving Surface Driven by Diffusion on the Surface},
  author = {Kov{\'a}cs, Bal{\'a}zs and Li, Buyang and Lubich, Christian and Power Guerra, Christian A.},
  year = 2017,
  month = nov,
  journal = {Numer. Math.},
  volume = {137},
  number = {3},
  pages = {643--689},
  issn = {0945-3245},
  doi = {10/gh4xvx},
  urldate = {2024-10-22},
  langid = {english}
}

@article{KovacsPowerGuerra2018,
  title = {Higher Order Time Discretizations with {{ALE}} Finite Elements for Parabolic Problems on Evolving Surfaces},
  author = {Kov{\'a}cs, Bal{\'a}zs and Power Guerra, Christian Andreas},
  year = 2018,
  month = jan,
  journal = {IMA Journal of Numerical Analysis},
  volume = {38},
  number = {1},
  pages = {460--494},
  issn = {0272-4979},
  doi = {10/gcz3xc},
  urldate = {2025-04-15}
}

@article{KrugerVarnikRaabe2011,
  title = {Efficient and Accurate Simulations of Deformable Particles Immersed in a Fluid Using a Combined Immersed Boundary Lattice {{Boltzmann}} Finite Element Method},
  author = {Kr{\"u}ger, T. and Varnik, F. and Raabe, D.},
  year = 2011,
  month = jun,
  journal = {Computers \& Mathematics with Applications},
  series = {Mesoscopic {{Methods}} for {{Engineering}} and {{Science}} --- {{Proceedings}} of {{ICMMES-09}}},
  volume = {61},
  number = {12},
  pages = {3485--3505},
  issn = {0898-1221},
  doi = {10.1016/j.camwa.2010.03.057},
  urldate = {2025-10-16}
}

@incollection{KuhnTucker_bookSection,
  title = {Nonlinear {{Programming}}},
  booktitle = {Proceedings of the {{Second Berkeley Symposium}} on {{Mathematical Statistics}} and {{Probability}}},
  author = {Kuhn, H. W. and Tucker, A. W.},
  year = 1951,
  month = jan,
  volume = {2},
  pages = {481--493},
  publisher = {University of California Press},
  urldate = {2025-10-12}
}

@article{LaadhariSaramitoMisbah2014,
  title = {Computing the Dynamics of Biomembranes by Combining Conservative Level Set and Adaptive Finite Element Methods},
  author = {Laadhari, Aymen and Saramito, Pierre and Misbah, Chaouqi},
  year = 2014,
  month = apr,
  journal = {Journal of Computational Physics},
  volume = {263},
  pages = {328--352},
  issn = {0021-9991},
  doi = {10.1016/j.jcp.2013.12.032},
  urldate = {2025-10-16}
}

@article{LeeCatheyWuEtAl2020,
  title = {Endoplasmic Reticulum Contact Sites Regulate the Dynamics of Membraneless Organelles},
  author = {Lee, Jason E. and Cathey, Peter I. and Wu, Haoxi and Parker, Roy and Voeltz, Gia K.},
  year = 2020,
  month = jan,
  journal = {Science},
  volume = {367},
  number = {6477},
  pages = {eaay7108},
  issn = {1095-9203},
  doi = {10.1126/science.aay7108},
  langid = {english},
  pmcid = {PMC10088059},
  pmid = {32001628}
}

@article{LeeLaughlinAnglivieldeLaBeaumelleEtAl2020,
  title = {{{3D}} Mesh Processing Using {{GAMer}} 2 to Enable Reaction-Diffusion Simulations in Realistic Cellular Geometries},
  author = {Lee, Christopher T. and Laughlin, Justin G. and {Angliviel de La Beaumelle}, Nils and Amaro, Rommie E. and McCammon, J. Andrew and Ramamoorthi, Ravi and Holst, Michael and Rangamani, Padmini},
  year = 2020,
  month = apr,
  journal = {PLoS Comput Biol},
  volume = {16},
  number = {4},
  pages = {e1007756},
  issn = {1553-7358},
  doi = {10.1371/journal.pcbi.1007756},
  langid = {english},
  pmcid = {PMC7162555},
  pmid = {32251448}
}

@article{LiuHuangXiaoEtAl2024,
  title = {Mathematical Modeling and Numerical Simulation of the {{N-component Cahn-Hilliard}} Model on Evolving Surfaces},
  author = {Liu, Lulu and Huang, Shijie and Xiao, Xufeng and Feng, Xinlong},
  year = 2024,
  month = sep,
  journal = {Journal of Computational Physics},
  volume = {513},
  pages = {113189},
  issn = {0021-9991},
  doi = {10.1016/j.jcp.2024.113189},
  urldate = {2025-10-06}
}

@article{LiuWangZhou2018,
  title = {Positivity-Preserving and Asymptotic Preserving Method for {{2D Keller-Segal}} Equations},
  author = {Liu, Jian-Guo and Wang, Li and Zhou, Zhennan},
  year = 2018,
  month = may,
  journal = {Math. Comp.},
  volume = {87},
  number = {311},
  pages = {1165--1189},
  issn = {0025-5718, 1088-6842},
  doi = {10.1090/mcom/3250},
  urldate = {2025-10-06},
  langid = {english}
}

@article{LiYangZhou2020,
  title = {Arbitrarily {{High-Order Exponential Cut-Off Methods}} for {{Preserving Maximum Principle}} of {{Parabolic Equations}}},
  author = {Li, Buyang and Yang, Jiang and Zhou, Zhi},
  year = 2020,
  month = jan,
  journal = {SIAM J. Sci. Comput.},
  volume = {42},
  number = {6},
  pages = {A3957-A3978},
  publisher = {{Society for Industrial and Applied Mathematics}},
  issn = {1064-8275},
  doi = {10/g5c5sd},
  urldate = {2025-05-27}
}

@article{LiZhang2020,
  title = {On the Monotonicity and Discrete Maximum Principle of the Finite Difference Implementation of \$\${{C}}\textasciicircum 0\$\$-\$\${{Q}}\textasciicircum 2\$\$finite Element Method},
  author = {Li, Hao and Zhang, Xiangxiong},
  year = 2020,
  month = jun,
  journal = {Numer. Math.},
  volume = {145},
  number = {2},
  pages = {437--472},
  issn = {0945-3245},
  doi = {10.1007/s00211-020-01110-6},
  urldate = {2025-07-16},
  langid = {english}
}

@article{LomakinLeeHanEtAl2015,
  title = {Competition for Actin between Two Distinct {{F-actin}} Networks Defines a Bistable Switch for Cell Polarization},
  author = {Lomakin, Alexis J. and Lee, Kun-Chun and Han, Sangyoon J. and Bui, Duyen A. and Davidson, Michael and Mogilner, Alex and Danuser, Gaudenz},
  year = 2015,
  month = nov,
  journal = {Nat Cell Biol},
  volume = {17},
  number = {11},
  pages = {1435--1445},
  issn = {1476-4679},
  doi = {10.1038/ncb3246},
  langid = {english},
  pmcid = {PMC4628555},
  pmid = {26414403}
}

@article{LubichMansourVenkataraman2013,
  title = {Backward Difference Time Discretization of Parabolic Differential Equations on Evolving Surfaces},
  author = {Lubich, Christian and Mansour, Dhia and Venkataraman, Chandrasekhar},
  year = 2013,
  month = oct,
  journal = {IMA Journal of Numerical Analysis},
  volume = {33},
  number = {4},
  pages = {1365--1385},
  issn = {0272-4979},
  doi = {10/mh2},
  urldate = {2024-06-13}
}

@article{LuHuangVanVleck2013,
  title = {The Cutoff Method for the Numerical Computation of Nonnegative Solutions of Parabolic {{PDEs}} with Application to Anisotropic Diffusion and {{Lubrication-type}} Equations},
  author = {Lu, Changna and Huang, Weizhang and Van Vleck, Erik S.},
  year = 2013,
  month = jun,
  journal = {Journal of Computational Physics},
  volume = {242},
  pages = {24--36},
  issn = {0021-9991},
  doi = {10/f4x4pn},
  urldate = {2025-05-27}
}

@article{M.ElliottFritz2017,
  title = {On Approximations of the Curve Shortening Flow and of the Mean Curvature Flow Based on the {{DeTurck}} Trick},
  author = {M. Elliott, Charles and Fritz, Hans},
  year = 2017,
  month = apr,
  journal = {IMA Journal of Numerical Analysis},
  volume = {37},
  number = {2},
  pages = {543--603},
  issn = {0272-4979},
  doi = {10/gmbnnc},
  urldate = {2025-06-04}
}

@article{MacDonaldMackenzieNolanEtAl2016,
  title = {A Computational Method for the Coupled Solution of Reaction-Diffusion Equations on Evolving Domains and Manifolds},
  author = {MacDonald, G. and Mackenzie, J.A. and Nolan, M. and Insall, R.H.},
  year = 2016,
  month = mar,
  journal = {J. Comput. Phys.},
  volume = {309},
  number = {C},
  pages = {207--226},
  issn = {0021-9991},
  doi = {10.1016/j.jcp.2015.12.038},
  urldate = {2025-06-30}
}

@article{MackenzieRowlattInsall2021,
  title = {A {{Conservative Finite Element ALE Scheme}} for {{Mass-Conservative Reaction-Diffusion Equations}} on {{Evolving Two-Dimensional Domains}}},
  author = {Mackenzie, John and Rowlatt, Christopher and Insall, Robert},
  year = 2021,
  month = jan,
  journal = {SIAM J. Sci. Comput.},
  volume = {43},
  number = {1},
  pages = {B132-B166},
  publisher = {{Society for Industrial and Applied Mathematics}},
  issn = {1064-8275},
  doi = {10/gh4xq8},
  urldate = {2024-06-06}
}

@article{MerckerMarciniak-CzochraRichterEtAl2013,
  title = {Modeling and {{Computing}} of {{Deformation Dynamics}} of {{Inhomogeneous Biological Surfaces}}},
  author = {Mercker, M. and {Marciniak-Czochra}, A. and Richter, T. and Hartmann, D.},
  year = 2013,
  month = jan,
  journal = {SIAM J. Appl. Math.},
  volume = {73},
  number = {5},
  pages = {1768--1792},
  publisher = {{Society for Industrial and Applied Mathematics}},
  issn = {0036-1399},
  doi = {10.1137/120885553},
  urldate = {2025-09-09}
}

@article{MokbelMokbelLieseEtAl2024,
  title = {A {{Simulation Method}} for the {{Wetting Dynamics}} of {{Liquid Droplets}} on {{Deformable Membranes}}},
  author = {Mokbel, Marcel and Mokbel, Dominic and Liese, Susanne and Weber, Christoph and Aland, Sebastian},
  year = 2024,
  month = dec,
  journal = {SIAM J. Sci. Comput.},
  volume = {46},
  number = {6},
  pages = {B806-B829},
  publisher = {{Society for Industrial and Applied Mathematics}},
  issn = {1064-8275},
  doi = {10.1137/24M1641142},
  urldate = {2025-10-06}
}

@book{Mullins_book,
  title = {The Biogenesis of Cellular Organelles},
  author = {Mullins, Chris},
  year = 2005,
  series = {Molecular Biology Intelligence Unit},
  publisher = {Landes Bioscience/Eurekah.com},
  address = {Georgetown, Tex},
  isbn = {978-0-306-47990-8},
  langid = {english},
  lccn = {2005 G-529}
}

@article{Nitsche1993,
  title = {Boundary Value Problems for Variational Integrals Involving Surface Curvatures},
  author = {Nitsche, Johannes C. C.},
  year = 1993,
  journal = {Quart. Appl. Math.},
  volume = {51},
  number = {2},
  pages = {363--387},
  issn = {0033-569X, 1552-4485},
  doi = {10.1090/qam/1218374},
  urldate = {2025-10-06},
  langid = {english}
}

@phdthesis{Nobile2001,
  title = {Numerical Approximation of Fluid-Structure Interaction Problems with Application to Haemodynamics},
  author = {Nobile, Fabio},
  year = 2001,
  doi = {10.5075/epfl-thesis-2458},
  urldate = {2024-10-22},
  langid = {english},
  school = {EPFL}
}

@article{Noguchi2004,
  title = {Fluid {{Vesicles}} with {{Viscous Membranes}} in {{Shear Flow}}},
  author = {Noguchi, Hiroshi},
  year = 2004,
  journal = {Phys. Rev. Lett.},
  volume = {93},
  number = {25},
  doi = {10.1103/PhysRevLett.93.258102}
}

@article{OlshanskiiPalzhanovQuaini2022,
  title = {A {{Comparison}} of {{Cahn}}--{{Hilliard}} and {{Navier}}--{{Stokes}}--{{Cahn}}--{{Hilliard Models}} on {{Manifolds}}},
  author = {Olshanskii, Maxim and Palzhanov, Yerbol and Quaini, Annalisa},
  year = 2022,
  month = oct,
  journal = {Vietnam J. Math.},
  volume = {50},
  number = {4},
  pages = {929--945},
  issn = {2305-2228},
  doi = {10.1007/s10013-022-00564-5},
  urldate = {2025-10-29},
  langid = {english}
}

@article{OngLai2020,
  title = {An Immersed Boundary Projection Method for Simulating the Inextensible Vesicle Dynamics},
  author = {Ong, Kian Chuan and Lai, Ming-Chih},
  year = 2020,
  month = may,
  journal = {Journal of Computational Physics},
  volume = {408},
  pages = {109277},
  issn = {0021-9991},
  doi = {10.1016/j.jcp.2020.109277},
  urldate = {2025-10-16}
}

@article{PalzhanovZhiliakovQuainiEtAl2021,
  title = {A Decoupled, Stable, and Linear {{FEM}} for a Phase-Field Model of Variable Density Two-Phase Incompressible Surface Flow},
  author = {Palzhanov, Yerbol and Zhiliakov, Alexander and Quaini, Annalisa and Olshanskii, Maxim},
  year = 2021,
  month = dec,
  journal = {Computer Methods in Applied Mechanics and Engineering},
  volume = {387},
  pages = {114167},
  issn = {0045-7825},
  doi = {10.1016/j.cma.2021.114167},
  urldate = {2025-10-29}
}

@misc{PorrmannBartelsVoigt2025,
  title = {Self-Avoiding Fluid Deformable Surfaces},
  author = {Porrmann, Maik and Bartels, S{\"o}ren and Voigt, Axel},
  year = 2025,
  month = sep,
  number = {arXiv:2509.22110},
  eprint = {2509.22110},
  primaryclass = {physics},
  publisher = {arXiv},
  doi = {10.48550/arXiv.2509.22110},
  urldate = {2025-10-05},
  archiveprefix = {arXiv}
}

@article{RosolenPecoArroyo2013,
  title = {An Adaptive Meshfree Method for Phase-Field Models of Biomembranes. {{Part I}}: {{Approximation}} with Maximum-Entropy Basis Functions},
  shorttitle = {An Adaptive Meshfree Method for Phase-Field Models of Biomembranes. {{Part I}}},
  author = {Rosolen, A. and Peco, C. and Arroyo, M.},
  year = 2013,
  month = sep,
  journal = {Journal of Computational Physics},
  volume = {249},
  pages = {303--319},
  issn = {0021-9991},
  doi = {10.1016/j.jcp.2013.04.046},
  urldate = {2025-10-16}
}

@article{Sauer2014,
  title = {Stabilized Finite Element Formulations for Liquid Membranes and Their Application to Droplet Contact},
  author = {Sauer, Roger A.},
  year = 2014,
  journal = {International Journal for Numerical Methods in Fluids},
  volume = {75},
  number = {7},
  pages = {519--545},
  issn = {1097-0363},
  doi = {10.1002/fld.3905},
  urldate = {2025-10-22},
  copyright = {Copyright \copyright{} 2014 John Wiley \& Sons, Ltd.},
  langid = {english}
}

@article{Sauer2025,
  title = {A Curvilinear Surface {{ALE}} Formulation for Self-Evolving {{Navier-Stokes}} Manifolds -- Stabilized Finite Element Formulation},
  author = {Sauer, Roger A.},
  year = 2025,
  month = dec,
  journal = {Computer Methods in Applied Mechanics and Engineering},
  volume = {447},
  pages = {118331},
  issn = {0045-7825},
  doi = {10.1016/j.cma.2025.118331},
  urldate = {2025-09-16}
}

@article{Schoberl2014,
  title = {C++11 {{Implementation}} of {{Finite Elements}} in {{NGSolve}}},
  author = {Sch{\"o}berl, Joachim},
  year = 2014,
  pages = {1--23},
  publisher = {{Institute of Analysis and Scientific Computing, TU Wien}},
  urldate = {2025-10-29},
  langid = {english}
}

@article{Seifert1997,
  title = {Configurations of Fluid Membranes and Vesicles},
  author = {Seifert, Udo},
  year = 1997,
  month = feb,
  journal = {Advances in Physics},
  volume = {46},
  number = {1},
  pages = {13--137},
  publisher = {Taylor \& Francis},
  issn = {0001-8732},
  doi = {10/dpkxvq},
  urldate = {2024-09-23}
}

@article{Sens2020,
  title = {Stick-Slip Model for Actin-Driven Cell Protrusions, Cell Polarization, and Crawling},
  author = {Sens, Pierre},
  year = 2020,
  month = oct,
  journal = {Proc Natl Acad Sci U S A},
  volume = {117},
  number = {40},
  pages = {24670--24678},
  issn = {1091-6490},
  doi = {10.1073/pnas.2011785117},
  langid = {english},
  pmcid = {PMC7547270},
  pmid = {32958682}
}

@article{ShinBrangwynne2017,
  title = {Liquid Phase Condensation in Cell Physiology and Disease},
  author = {Shin, Yongdae and Brangwynne, Clifford P.},
  year = 2017,
  month = sep,
  journal = {Science},
  volume = {357},
  number = {6357},
  pages = {eaaf4382},
  publisher = {American Association for the Advancement of Science},
  doi = {10.1126/science.aaf4382},
  urldate = {2025-10-16}
}

@article{SneadJalihalGerbichEtAl2022,
  title = {Membrane Surfaces Regulate Assembly of Ribonucleoprotein Condensates},
  author = {Snead, Wilton T. and Jalihal, Ameya P. and Gerbich, Therese M. and Seim, Ian and Hu, Zhongxiu and Gladfelter, Amy S.},
  year = 2022,
  month = apr,
  journal = {Nat Cell Biol},
  volume = {24},
  number = {4},
  pages = {461--470},
  issn = {1476-4679},
  doi = {10.1038/s41556-022-00882-3},
  langid = {english},
  pmcid = {PMC9035128},
  pmid = {35411085}
}

@book{Steigmann_book,
  title = {The {{Role}} of {{Mechanics}} in the {{Study}} of {{Lipid Bilayers}}},
  editor = {Steigmann, David J.},
  year = 2018,
  series = {{{CISM International Centre}} for {{Mechanical Sciences}}},
  volume = {577},
  publisher = {Springer International Publishing},
  address = {Cham},
  doi = {10.1007/978-3-319-56348-0},
  urldate = {2024-10-21},
  copyright = {http://www.springer.com/tdm},
  isbn = {978-3-319-56347-3 978-3-319-56348-0}
}

@article{WangPalzhanovQuainiEtAl2022,
  title = {Lipid Domain Coarsening and Fluidity in Multicomponent Lipid Vesicles: {{A}} Continuum Based Model and Its Experimental Validation},
  shorttitle = {Lipid Domain Coarsening and Fluidity in Multicomponent Lipid Vesicles},
  author = {Wang, Y. and Palzhanov, Y. and Quaini, A. and Olshanskii, M. and Majd, S.},
  year = 2022,
  month = jul,
  journal = {Biochimica et Biophysica Acta (BBA) - Biomembranes},
  volume = {1864},
  number = {7},
  pages = {183898},
  issn = {0005-2736},
  doi = {10.1016/j.bbamem.2022.183898},
  urldate = {2025-10-29}
}

@article{YangVenkataramanStylesEtAl2016,
  title = {A Computational Framework for Particle and Whole Cell Tracking Applied to a Real Biological Dataset},
  author = {Yang, Feng Wei and Venkataraman, Chandrasekhar and Styles, Vanessa and Kuttenberger, Verena and Horn, Elias and {von Guttenberg}, Zeno and Madzvamuse, Anotida},
  year = 2016,
  month = may,
  journal = {Journal of Biomechanics},
  series = {{{SI}}: {{Motility}} and Dynamics of Living Cells in Health, Disease and Healing},
  volume = {49},
  number = {8},
  pages = {1290--1304},
  issn = {0021-9290},
  doi = {10.1016/j.jbiomech.2016.02.008},
  urldate = {2025-10-29}
}

@misc{ZhangHuang2025,
  title = {A Unifying Moving Mesh Method for Curves, Surfaces, and Domains Based on Mesh Equidistribution and Alignment},
  author = {Zhang, Min and Huang, Weizhang},
  year = 2025,
  month = jan,
  number = {arXiv:2501.03086},
  eprint = {2501.03086},
  publisher = {arXiv},
  doi = {10.48550/arXiv.2501.03086},
  urldate = {2025-01-08},
  archiveprefix = {arXiv}
}

@article{ZhangShu2010,
  title = {On Positivity-Preserving High Order Discontinuous {{Galerkin}} Schemes for Compressible {{Euler}} Equations on Rectangular Meshes},
  author = {Zhang, Xiangxiong and Shu, Chi-Wang},
  year = 2010,
  month = nov,
  journal = {Journal of Computational Physics},
  volume = {229},
  number = {23},
  pages = {8918--8934},
  issn = {0021-9991},
  doi = {10.1016/j.jcp.2010.08.016},
  urldate = {2025-07-16}
}

@article{ZhangShu2011,
  title = {Maximum-Principle-Satisfying and Positivity-Preserving High-Order Schemes for Conservation Laws: Survey and New Developments},
  shorttitle = {Maximum-Principle-Satisfying and Positivity-Preserving High-Order Schemes for Conservation Laws},
  author = {Zhang, Xiangxiong and Shu, Chi-Wang},
  year = 2011,
  month = may,
  journal = {Proceedings of the Royal Society A: Mathematical, Physical and Engineering Sciences},
  volume = {467},
  number = {2134},
  pages = {2752--2776},
  publisher = {Royal Society},
  doi = {10/d4kh7c},
  urldate = {2024-10-22}
}

@article{ZhaoZhang2020,
  title = {Phase {{Separation}} in {{Membrane Biology}}: {{The Interplay}} between {{Membrane-Bound Organelles}} and {{Membraneless Condensates}}},
  shorttitle = {Phase {{Separation}} in {{Membrane Biology}}},
  author = {Zhao, Yan G. and Zhang, Hong},
  year = 2020,
  month = oct,
  journal = {Dev Cell},
  volume = {55},
  number = {1},
  pages = {30--44},
  issn = {1878-1551},
  doi = {10.1016/j.devcel.2020.06.033},
  langid = {english},
  pmid = {32726575}
}

@article{ZhiliakovWangQuainiEtAl2021,
  title = {Experimental Validation of a Phase-Field Model to Predict Coarsening Dynamics of Lipid Domains in Multicomponent Membranes},
  author = {Zhiliakov, A. and Wang, Y. and Quaini, A. and Olshanskii, M. and Majd, S.},
  year = 2021,
  month = jan,
  journal = {Biochimica et Biophysica Acta (BBA) - Biomembranes},
  volume = {1863},
  number = {1},
  pages = {183446},
  issn = {0005-2736},
  doi = {10.1016/j.bbamem.2020.183446},
  urldate = {2025-10-06}
}

@book{AndreaBonito2020,
  title = {Geometric Partial Differential Equations - Part {{I}}},
  editor = {Bonito, Andrea and Nochetto, Ricardo H.},
  year = 2020,
  series = {Handbook of Numerical Analysis},
  volume = {21}
}

@article{BurmanFernandez2009a,
  title = {Finite Element Methods with Symmetric Stabilization for the Transient Convection--Diffusion--Reaction Equation},
  author = {Burman, Erik and Fern{\'a}ndez, Miguel A},
  year = 2009,
  journal = {Comput. Methods Appl. Mech. Eng.},
  volume = {198},
  number = {33-36},
  pages = {2508--2519},
  publisher = {Elsevier},
  fjournal = {Computer Methods in Applied Mechanics and Engineering}
}

@article{BurmanHansbo2004,
  title = {Edge Stabilization for {{Galerkin}} Approximations of Convection--Diffusion--Reaction Problems},
  author = {Burman, E. and Hansbo, Peter},
  year = 2004,
  month = apr,
  journal = {Comput. Methods Appl. Mech. Engrg.},
  volume = {193},
  number = {15-16},
  pages = {1437--1453},
  doi = {10/fvsvh9}
}

@article{BurmanHansboLarsonEtAl2018b,
  title = {Stabilized {{CutFEM}} for the Convection Problem on Surfaces},
  author = {Burman, E. and Hansbo, P. and Larson, M. G. and Zahedi, S.},
  year = 2019,
  journal = {Numer. Math.},
  volume = {141},
  number = {1},
  doi = {10/gd8zcs}
}

@article{BurmanHansboLarsonEtAl2020,
  title = {A Stabilized Cut Streamline Diffusion Finite Element Method for Convection--Diffusion Problems on Surfaces},
  author = {Burman, Erik and Hansbo, Peter and Larson, Mats G. and Massing, Andr{\'e} and Zahedi, Sara},
  year = 2020,
  journal = {Comput. Methods Appl. Mech. Engrg.},
  volume = {358},
  pages = {112645},
  issn = {0045-7825},
  doi = {10/gf8vwp}
}

@article{DednerMadhavan2015,
  title = {Discontinuous {{Galerkin}} Methods for Hyperbolic and Advection-Dominated Problems on Surfaces},
  author = {Dedner, A. and Madhavan, P.},
  year = 2015,
  month = may,
  journal = {ArXiv e-prints},
  eprintclass = {math.NA},
  doi = {10/grm2g4},
  archiveprefix = {arXiv}
}

@incollection{GanesanHahnSimonEtAl2017,
  title = {{{ALE-FEM}} for Two-Phase and Free Surface Flows with Surfactants},
  booktitle = {Transport Processes at Fluidic Interfaces},
  author = {Ganesan, Sashikumaar and Hahn, Andreas and Simon, Kristin and Tobiska, Lutz},
  year = 2017,
  pages = {5--31},
  publisher = {Springer}
}

@book{GrossReusken2011,
  title = {Numerical {{Methods}} for {{Two-phase Incompressible Flows}}},
  author = {Gross, Sven and Reusken, Arnold},
  year = 2011,
  series = {Springer {{Series}} in {{Computational Mathematics}}},
  volume = {40},
  publisher = {Springer},
  address = {Berlin, Heidelberg},
  doi = {10.1007/978-3-642-19686-7},
  urldate = {2023-01-13},
  isbn = {978-3-642-19685-0 978-3-642-19686-7}
}

@article{HansboLarsonZahedi2015,
  title = {Characteristic Cut Finite Element Methods for Convection--Diffusion Problems on Time Dependent Surfaces},
  author = {Hansbo, Peter and Larson, Mats G. and Zahedi, Sara},
  year = 2015,
  month = aug,
  journal = {Computer Methods in Applied Mechanics and Engineering},
  volume = {293},
  pages = {431--461},
  issn = {0045-7825},
  doi = {10/f3pbxn},
  urldate = {2025-06-09}
}

@misc{NeivaTurlier2025,
  title = {Unfitted Finite Element Modelling of Surface-Bulk Viscous Flows in Animal Cells},
  author = {Neiva, Eric and Turlier, Herv{\'e}},
  year = 2025,
  month = may,
  number = {arXiv:2505.05723},
  eprint = {2505.05723},
  primaryclass = {cs},
  publisher = {arXiv},
  doi = {10.48550/arXiv.2505.05723},
  urldate = {2025-05-13},
  archiveprefix = {arXiv}
}

@article{Nitsche1971,
  title = {\"Uber Ein {{Variationsprinzip}} Zur {{L\"osung}} von {{Dirichlet-Problemen}} Bei {{Verwendung}} von {{Teilr\"aumen}}, Die Keinen {{Randbedingungen}} Unterworfen Sind},
  author = {Nitsche, J.},
  year = 1971,
  month = jul,
  journal = {Abhandlungen aus dem Mathematischen Seminar der Universit\"at Hamburg},
  volume = {36},
  number = {1},
  pages = {9--15},
  doi = {10/bps9r9}
}

@article{OlshanskiiReuskenXu2014,
  title = {A Stabilized Finite Element Method for Advection--Diffusion Equations on Surfaces},
  author = {Olshanskii, M. A. and Reusken, A. and Xu, X.},
  year = 2014,
  journal = {IMA J. Numer. Anal.},
  volume = {34},
  number = {2},
  pages = {732--758},
  publisher = {Oxford University Press},
  doi = {10/f5zf9x}
}

@article{SimonTobiska2019,
  title = {Local Projection Stabilization for Convection--Diffusion--Reaction Equations on Surfaces},
  author = {Simon, K and Tobiska, L},
  year = 2019,
  journal = {Comput. Methods Appl. Mech. Engrg.},
  volume = {344},
  pages = {34--53},
  publisher = {Elsevier},
  doi = {10/grm2gx}
}

@article{UlfsbyMassingSticko2023,
  title = {Stabilized Cut Discontinuous {{Galerkin}} Methods for Advection--Reaction Problems on Surfaces},
  author = {Ulfsby, Tale Bakken and Massing, Andr{\'e} and Sticko, Simon},
  year = 2023,
  month = aug,
  journal = {Computer Methods in Applied Mechanics and Engineering},
  volume = {413},
  pages = {116109},
  issn = {0045-7825},
  doi = {10/gsb7hg},
  urldate = {2023-06-12},
  langid = {english}
}

@article{XiaoZhaoFeng2020,
  title = {A Layers Capturing Type {{H-adaptive}} Finite Element Method for Convection--Diffusion--Reaction Equations on Surfaces},
  author = {Xiao, Xufeng and Zhao, Jianping and Feng, Xinlong},
  year = 2020,
  month = apr,
  journal = {Computer Methods in Applied Mechanics and Engineering},
  volume = {361},
  pages = {112792},
  issn = {0045-7825},
  doi = {10/g9prhn},
  urldate = {2025-06-09}
}

\end{document}